\documentclass{siamltex}
\usepackage{amsfonts}
\usepackage{amssymb,amsmath}
\usepackage{graphicx}
\usepackage{color}
\usepackage{hyperref}
\usepackage{subfigure}

\DeclareMathAlphabet{\itbf}{OML}{cmm}{b}{it}
\newcommand{\nc}{\newcommand}
\nc{\bx}{{\bf x}}
\nc{\bu}{{\bf u}}
\nc{\vx}{\vec{\bx}}
\nc{\ve}{\vec{\bf e}}
\nc{\vn}{\vec{\bf n}}
\nc{\bn}{{\bf n}}
\nc{\de}{\delta}
\nc{\la}{\lambda}
\nc{\eps}{\varepsilon}
\nc{\ep}{\varepsilon}
\nc{\om}{\omega}

\nc{\cD}{\mathcal D}
\nc{\cR}{\mathcal{R}}
\nc{\cA}{\mathcal A}
\nc{\cT}{\mathcal T}
\nc{\cS}{\mathcal S}
\nc{\cL}{\mathcal L}
\nc{\cW}{\mathcal W}
\nc{\cC}{\mathcal C}
\nc{\cU}{\mathcal U}
\nc{\cB}{\mathfrak B}
\nc{\brho}{{\boldsymbol{\rho}}}
\nc{\bpsi}{\boldsymbol \Psi}
\nc{\bphi}{{\boldsymbol{\varphi}}}
\nc{\bdeta}{{\boldsymbol{\eta}}}
\nc{\bM}{{\boldsymbol{\mathfrak{M}}}}
\nc{\bX}{{\boldsymbol{\mathfrak{X}}}}
\nc{\EE}{{\mathbb{E}}}
\nc{\PP}{{\mathbb{P}}}

\renewcommand{\hat}{\widehat}

\begin{document}

\title{Source estimation with incoherent waves in random waveguides}

\author{ Sebastian Acosta\footnotemark[1], Ricardo Alonso\footnotemark[2] and 
Liliana Borcea\footnotemark[3]}

\maketitle 

\renewcommand{\thefootnote}{\fnsymbol{footnote}}

\footnotetext[1]{Baylor College of Medicine, Houston, TX 77005. 
{\tt sacosta@bcm.edu}}
\footnotetext[2]{Departamento de Matem\'{a}tica, PUC--Rio, Brasil. {\tt
    ralonso@mat.puc-rio.br}} 
\footnotetext[3]{Department of
  Mathematics, University of Michigan, Ann Arbor, MI 48109. {\tt
    borcea@umich.edu}}

\begin{abstract}
We study an inverse source problem for the acoustic wave equation in a
random waveguide. The goal is to estimate the source of waves from
measurements of the acoustic pressure at a remote array of
sensors. The waveguide effect is due to boundaries that trap the waves
and guide them in a preferred (range) direction, the waveguide axis,
along which the medium is unbounded.  The random waveguide is a model
of perturbed ideal waveguides which have flat boundaries and are
filled with known media that do not change with range. The
perturbation consists of fluctuations of the boundary and of the wave
speed due to numerous small inhomogeneities in the medium.  The
fluctuations are uncertain in applications, which is why we model them
with random processes, and they cause significant cumulative
scattering at long ranges from the source. The scattering effect
manifests mathematically as an exponential decay of the expectation of
the acoustic pressure, the coherent part of the wave.  The incoherent
wave is modeled by the random fluctuations of the acoustic pressure, which
dominate the expectation at long ranges from the source. We use the
existing theory of wave propagation in random waveguides to analyze
the inverse problem of estimating the source from incoherent wave
recordings at remote arrays. We show how to obtain from the incoherent
measurements high fidelity estimates of the time resolved energy
carried by the waveguide modes, and study the invertibility of the
system of transport equations that model energy propagation in order
to estimate the source.
\end{abstract}
\renewcommand{\thefootnote}{\arabic{footnote}}


\section{Introduction}
\label{sect:intro}
We study an inverse problem for the scalar (acoustic) wave equation,
where we wish to estimate the source of waves from measurements of the
acoustic pressure field $p(t,\vx)$ at a remote array of receiver
sensors. The waves propagate in a waveguide, meaning that they are
trapped by boundaries and are guided in the range direction, the
waveguide axis, along which the medium is unbounded. Ideally the
boundaries are straight and the medium does not change with range. We
consider perturbed waveguides filled with heterogeneous media, where
the boundary and the wave speed have small fluctuations on scales
similar to the wavelength. These fluctuations have little effect in
the vicinity of the source, but they are important at long ranges
because they cause significant cumulative wave scattering. We suppose
that the array of receivers is far from the source, as is typical in
applications in underwater acoustics, sound propagation in corrugated
pipes, in tunnels, etc., and study how cumulative scattering impedes
the inversion.

In most setups the fluctuations are uncertain, which is why we
introduce a stochastic framework and model them with random
processes. The inversion is carried in only one perturbed waveguide,
meaning that the array measures one realization of the random pressure
field, the solution of the wave equation in that waveguide. The
stochastic framework allows us to study the chain of mappings from the
uncertainty in the waveguide to the uncertainty of the array
measurements and of the inversion results. The goal is to understand
how to process the uncertain data and quantify what can be estimated
about the source in a reliable (statistically stable)
manner. Statistical stability means that the estimates do not change
with the realization of the fluctuations of the waveguide, which are
unknown.

The problem of imaging (localizing) sources in waveguides has been
studied extensively in underwater acoustics
\cite{baggeroer1993overview,krolik1992matched,heaney1998very,abadi2012blind}.
Typical imaging approaches are matched field and related coherent
methods that match the measured $p(t,\vx)$ with its mathematical model
for search locations of the source. The model is based on wave
propagation in ideal waveguides and the imaging is successful when
$p(t,\vx)$ is mostly coherent. The coherent part of $p(t,\vx)$ is its
statistical expectation $\EE[p(t,\vx)]$ with respect to realizations
of the random waveguide, and the incoherent field is modeled by
$p(t,\vx)-\EE[p(t,\vx)]$. As the waves propagate in the random
waveguide they lose coherence due to scattering by the fluctuations of
the boundary and the inhomogeneities in the medium. This manifests as
an exponential decay in range of the expectation $\EE[p(t,\vx)]$, and
strengthening of the fluctuations $p(t,\vx)-\EE[p(t,\vx)]$. 

Detailed studies of the loss of coherence of sound waves due to
cumulative scattering are given in
\cite{kohler77,dozier,garnier_papa,gomez,garnier2008effective} for
waveguides filled with randomly heterogeneous media and in
\cite{ABG-12,gomez1} for waveguides with random boundaries. These
waveguides are two dimensional models of the ocean, and they may leak
(radiate) in the ocean floor. The problem is similar in three
dimensional acoustic waveguides with bounded cross-section. We refer
to \cite{BG-13} for wave propagation in three dimensional waveguide
models of the ocean which have unbounded cross-section and random
pressure release top boundary, and to \cite{AB_EM,marcuse91} for three
dimensional electromagnetic random waveguides.  In all cases the
analysis of loss of coherence is based on the decomposition of the
wave field in an infinite set of monochromatic waves called waveguide
modes, which are special solutions of the wave equation in the ideal
waveguide. Finitely many modes are propagating waves, and we may
associate them with plane waves that strike the boundary at different
angles of incidence and are reflected repeatedly.  The remaining
infinitely many modes are evanescent and/or radiating waves. The
cumulative scattering in the waveguide is modeled by fluctuations of
the amplitudes of the modes. When scattering is weak, as is the
case at moderate distances from the source, the amplitudes are
approximately constant in range, and they are determined solely by the
source excitation. Scattering builds up over long ranges and the mode
amplitudes become random fields with exponentially decaying
expectation on range scales called scattering mean free paths.

The mode dependence of the scattering mean free paths is analyzed in
\cite{ABG-12}. It turns out that the slow modes, which correspond to
plane waves that strike the boundary at almost normal incidence, are
most affected by scattering. These waves have long trajectories from
the source to the array, and thus interact more with the boundary and
medium fluctuations. We refer to \cite{BGT-14} for an adaptive
coherent imaging approach which detects which modes are incoherent and
filters them out from the measurements in order to achieve statistically
stable results. See also the results in
\cite{krolik1992matched,heaney1998very,yoo1999broadband}. However, 
when the array is farther from the source than the scattering mean free
paths of all the modes, the data is incoherent and coherent imaging
methods like matched field cannot work.  In this paper we assume that
this is the case and study an inversion approach based on a system of
transport equations that models the propagation of energy carried by
the modes. This system is derived in
\cite{kohler77,garnier_papa,ABG-12} and is used in \cite{BIT-10} to
estimate the location of a point source in random waveguides. Here we
study the inverse problem in more detail and answer the following
questions: (1) How can we obtain reliable estimates of the mode
energies from the incoherent pressure field measured at the array?
(2) What kind of information about the source can we recover from the
transport equations? (3) Can we quantify the deterioration of the
inversion results in terms of the range offset between the source and
the array? 

We begin in section \ref{sect:setup} with the mathematical formulation
of the inverse problem, and recall in section \ref{sect:stat} the
model of the random wave field $p(t,\vx)$ derived in
\cite{kohler77,garnier_papa,ABG-12,FRANCE}. The main results of the
paper are in sections \ref{sect:transpIm} and \ref{sect:inverse}. We
motivate there the inversion based on energy transport, and describe
the forward mapping from the source to the expectation of the time
resolved energy carried by the modes. We show how to calculate this
energy from the incoherent array data, and describe how to invert
approximately the transport equations.  The results quantify the
limited information that can be recovered about the source. We end
with a summary in section \ref{sect:summary}.

\section{Formulation of the problem}
\label{sect:setup}
\label{sect:2D}
\begin{figure}[t!]
\vspace{-0.6in}
     \hspace{-0.3in}\input{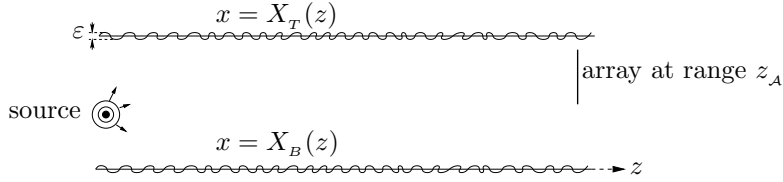}
   \vspace{-0.05in}
   \caption{Schematic of the problem setup. The source emits a signal
     in a waveguide and the wave field is recorded at a remote array.
     The perturbed waveguide has fluctuating boundaries and is filled
     with a medium with fluctuating wave speed.}
   \label{fig:Schem}
\end{figure}

We limit our study to two dimensional waveguides with reflecting
boundaries modeled by pressure release boundary conditions. This is
for simplicity, but the results extend to other boundary conditions
and to leaky and three dimensional waveguides, as discussed in section
\ref{sect:summary}.  We illustrate the setup in Figure
\ref{fig:Schem}, and introduce the system of coordinates $\vx = (x,z)$
with range $z$ originating from the center of the source.  The
waveguide occupies the domain
\[
\Omega = \left\{\vx = (x,z): ~ ~ x \in (X_{_B}(z), X_{_T}(z)), ~ ~ z
\in \mathbb{R} \right\},
\]
where the cross-range $x$ takes values between the bottom and top
boundaries modeled by $X_{_B}(z)$ and $X_{_T}(z)$. The source has an
unknown density $\rho(\vx)$ which is compactly supported in $\Omega$,
near $z = 0$, and emits a signal $F(t)$ which is a pulse $f(Bt)$ of
support of order $1/B$ around $t = 0$, modulated by an oscillatory
exponential
\begin{equation}
F(t) = e^{-i \om_o t} f(Bt).
\label{eq:fm1}
\end{equation}
We introduce the bandwidth $B$ in the argument of the pulse to
emphasize that the Fourier transform $\hat F(\om)$ of the signal is
supported in the interval $(\om_o - \pi B, \om_o+ \pi B)$ around the
central frequency $\om_o$,
\begin{equation}
\hat F(\om) = \int_{-\infty}^\infty d t\, F(t) e^{i \om t} =
\frac{1}{B} \hat f \left(\frac{\om-\om_o}{B} \right).
\end{equation}
The array is a collection of receivers that are placed close together
in the set
\[
{A} = \left\{\vx_{_\cA} = (x,z_{_\cA}): ~ ~ x \in \cA \subset
[X_{_B}(z), X_{_T}(z)] \right\},
\]
at range $z_{_\cA} > 0$ from the source, where $\cA$ is an interval
called the array aperture. The receivers record the acoustic pressure
field $p(t,\vx)$ modeled by the solution of the acoustic wave equation
\begin{equation}
\left[\partial_x^2 + \partial_z^2 - c^{-2}(\vx)
  \partial^2_t\right]p(t,\vx) = F(t) \rho(\vx), \quad \vx \in \Omega,
~ ~ t > 0,
\label{eq:fm2}
\end{equation}
with pressure release boundary conditions 
\begin{equation}
p(t,\vx) = 0, \quad t > 0, ~ ~ \vx \in \partial \Omega = \left\{\vx =
(x,z): ~ ~ x \in \{X_{_B}(z), X_{_T}(z)\}, ~ ~ z \in \mathbb{R}
\right\},
\label{eq:fm3}
\end{equation}
and initial
condition $ p(t,\vx) \equiv 0$ for $t \ll 0$. Here $c(\vx)$ is the
sound speed.

The inverse problem is to determine the source density $\rho(\vx)$
from the array data recordings $D(t,x)$. We model them by
\begin{equation}
D(t,x) = p(t,\vx_{_\cA}) 1_{_\cA}(x) \chi\left(\frac{t-t_o}{\cT}\right),
\quad \vx_{_\cA} = (x,z_{_\cA}),
\label{eq:fm4}
\end{equation}
using a recording time window $\chi$ centered at $t_o$ and of duration
$\cT$. We can take any continuous, compactly supported $\chi$, but we
assume henceforth that it equals one in the interval $(-1/2,1/2)$ and
tapers quickly to zero outside. We also approximate the array by a
continuum aperture in the interval $\cA$, and use the indicator
function $1_{_\cA}(x)$ which equals one when $x \in \cA$ and zero
otherwise.

\subsection{The random model of perturbed waveguides}
\label{sect:setuprand}
In ideal waveguides the sound speed is modeled by a function $c_o(x)$
that is independent of range and the boundaries are straight, meaning
that $X_{_B}(z) = 0$ and $X_{_T}(z) = X$, a constant. The sound speed
in the perturbed waveguide has fluctuations around $c_o$ and the
boundaries $X_{_B}$ and $X_{_T}$ fluctuate around $0$ and $X$. The
fluctuations are small, with amplitude quantified by a positive
dimensionless parameter $\ep \ll 1$.  It is used in
\cite{kohler77,garnier_papa,ABG-12,FRANCE} to analyze the pressure
field at properly scaled long ranges where scattering is significant,
in the asymptotic limit $\eps \to 0$.

We take $c_o$ constant for simplicity, to write explicitly the mode
decomposition, but the results extend easily to cross-range dependent
$c_o(x)$.  The perturbed sound speed $c(\vx)$ is modeled by
\begin{equation}
\frac{1}{c^2(\vx)} = \frac{1}{c_o^2} \left[ 1 + \ep_{c}\, \nu
  \left(\frac{\vx}{\ell}\right)\right],
\label{eq:fm5}
\end{equation}
where $\nu$ is a mean zero random process that is bounded almost
surely, so that the right hand side in (\ref{eq:fm5}) remains
positive.  We assume that $\nu$ is stationary and mixing in range,
meaning in particular that the auto-correlation
\begin{equation}
\cR_\nu(\xi,\xi',\eta) = \EE \left[ \nu(\xi,u) \nu(\xi',u + \eta)\right]
\label{eq:fm6}
\end{equation}
is absolutely integrable in the third argument over the real line. The
process $\nu$ is normalized by $\cR_\nu(0,0,0) = 1$ and
\[
\int_{-\infty}^\infty dz \,
\cR_\nu\left(\frac{x}{\ell},\frac{x'}{\ell},\frac{z}{\ell}\right) =
O(\ell),
\]
where $\ell$ is the correlation length, the range offset over which
the random fluctuations become statistically decorrelated. It compares
to the central wavelength $\la_o$ as $\ell \gtrsim \la_o$. The scaling
by the same $\ell$ of the cross-range in (\ref{eq:fm5}) means that the
heterogeneous medium is isotropic, but we could have $\ell_X =
O(\ell)$ as well, without changing the conclusions. The amplitude of
the fluctuations is scaled by $\ep_c$ which equals $\ep$ for a random
medium and zero for a homogeneous medium.

We model similarly the boundary fluctuations
\begin{equation}
X_{_B}(z) = \ep_{_B} \mu_{_B}\left(\frac{z}{\ell}\right), \qquad 
X_{_T}(z) = X \left[ 1 + \ep_{_T} \mu_{_T}\left(\frac{z}{\ell}\right)\right], 
\label{eq:fm7}
\end{equation}
using two mean zero, stationary and mixing random processes $\mu_{_B}$
and $\mu_{_T}$, that are bounded almost surely and have integrable
autocorrelation $\cR_{_B}$ and $\cR_{_T}$. We assume that $\nu$,
$\mu_{_B}$ and $\mu_{_T}$ are independent\footnote{If the random
  processes are not independent the moment formulae in this paper must
  be modified. Their derivation is a straightforward extension of
  the analysis in \cite{kohler77,garnier_papa,ABG-12}.}  and use the same
correlation length $\ell$ to simplify notation, but the results hold
for any scales $\ell_{_B}$ and $\ell_{_T}$ of the order of $\ell$. For
technical reasons related to the method of analysis used in
\cite{ABG-12} we also assume that the processes $\mu_{_B}$ and
$\mu_{_T}$ have bounded first and second derivatives, almost
surely. Less smooth boundary fluctuations are considered in
\cite{gomez1}. The boundary fluctuations are scaled by $\ep_{_B}$ and
$\ep_{_T}$ which can be $O(\ep)$, or they may be set to zero to study
separately the scattering effects of the medium and the boundary.

The theory of wave propagation in waveguides with long range
correlations of the random fluctuations of $c(\vx)$ is being developed
\cite{Gomez2}, and our results are expected to extend (with
modifications) to such settings. The case of turning waveguides with
smooth and large variations of the boundaries, on scales that are 
comparable to $z_{_\cA}$, is much more difficult.  The analysis of
wave propagation in such waveguides is quite involved
\cite{ting1983wave,felix2002multimodal,ahluwalia1974asymptotic} and
the mapping of random fluctuations of the sound speed to $p(t,\vx)$ is
not understood in detail, although it is considered formally in
\cite{marcuse91}.
\section{Cumulative scattering effects in the random waveguide}
\label{sect:stat}
We write the solution of the wave equation
(\ref{eq:fm2})-(\ref{eq:fm3}) as
\begin{equation}
p(t,\vx) = \int_{\Omega_{\rho}} d \vx' \, \rho(\vx') p(t,\vx,\vx'),
\label{eq:st1}
\end{equation}
where $p(t,\vx,\vx')$ is the wave field due to a point source at $\vx'
= (\bx',z')$, emitting the signal $F(t)$ defined in (\ref{eq:fm1}),
and $\Omega_{\rho} \subset \Omega$ is the compact support of the
source, which lies near $z' = 0$. The points $\vx = (\bx,z)$ in
(\ref{eq:st1}) are at range $z > z'$, for all $\vx' \in
\Omega_\rho$. 

It follows from \cite{kohler77,garnier_papa,ABG-12,FRANCE} that
$p(t,\vx,\vx')$ is a linear superposition of propagating and
evanescent waves, called waveguide modes
\begin{align}
p(t,\vx,\vx') = \int_{-\infty}^\infty \frac{d \om}{2 \pi B} \hat
f\left( \frac{\om-\om_o}{B}\right) e^{-i \om t} \Big[ \sum_{j=1}^N
a_{j}^+(\om,z,\vx') \Psi_j^+(\om,x,z-z') + \nonumber
\\ \sum_{j=1}^N a_j^-(\om,z,\vx') \Psi_j^-(\om,x,z-z') + \sum_{j =
  N+1}^\infty a_j^{e}(\om,z,\vx') \Psi_j^{e}(\om,x,z-z') \Big].
\label{eq:st2}
\end{align}
The modes are special solutions of the wave equation in the ideal
waveguide, and can be obtained with separation of variables. There are
$2 N$ propagating modes
\begin{equation}
\Psi^\pm_j(\om,x,z-z') = \phi_j(x) e^{\pm i \beta_j(\om) (z-z')},
\quad j = 1, \ldots. N,
\label{eq:st3}
\end{equation}
with index $+$ denoting forward going and $-$ backward going, and
infinitely many evanescent modes
\begin{equation}
\Psi^e_j(\om,x,z-z') = \phi_j(x) e^{-\beta_j(\om) |z-z'|}, \quad j > N. 
\label{eq:st4}
\end{equation}
They are defined by the complete and orthonormal set $\{\phi_j(x)\}_{j
  \ge 1}$ of eigenfunctions of the symmetric linear operator
$\mathbb{L}_x = \partial_x^2 + k^2$ with homogeneous Dirichlet
boundary conditions at $x = 0$ and $x = X$, where $k = \om/c_o$.
Because $c_o$ is constant we can write
\begin{equation}
\phi_j(x) = \sqrt{\frac{2}{X}} \sin \left(\frac{\pi j x}{X}\right),  
\label{eq:st5}
\end{equation}
and note explicitly how $\Psi_j^{+}$ are associated with monochromatic
plane waves that travel in the direction of the slowness vectors $(\pm
\pi j/X,\beta_j)$ and strike the boundaries where they reflect
according to Snell's law. The mode wavenumbers are denoted by
$\beta_j(\om)$, and are determined by the square root of the
eigenvalues of the operator $\mathbb{L}_x$
\begin{equation}
\beta_j(\om) = \left| k^2 - \left({\pi j}/{X}\right)^2 \right|^{1/2}.
\label{eq:st6}
\end{equation}
The $\beta_j$ of the propagating modes correspond to the first $N$
eigenvalues which are positive, where
\begin{equation}
N(\om) = \left \lfloor {k X}/{\pi} \right \rfloor 
\label{eq:st7}
\end{equation}
and $\lfloor \cdot \rfloor$ denotes the integer part. 

We assume for simplicity that the bandwidth $B$ is not too
large\footnote{In applications of imaging in open environments large
  bandwidths are desired for improved range resolution. In ideal
  waveguides good images can be formed with small bandwidths because
  the modes give different angle views of the support of the
  source. In random waveguides we may benefit from a large bandwidth,
  as explained in section \ref{sect:I6}. Such
  bandwidths may be divided in smaller sub-bands to which we can apply
  the analysis in this paper.}, so that there is the same number of
propagating modes for all the frequencies of the pulse, and drop the
dependence of $N$ on $\omega$.  We also suppose that there are no
standing waves, meaning that $\beta_j$ are bounded below by a positive
constant, for all $j \ge 1$.

The cumulative scattering effects in the random waveguide are modeled
by the mode amplitudes $\{a_j^\pm(\om,z,\vx')\}_{1\le j\le N}$ and
$\{a_j^e(\om,z,\vx')\}_{j>N}$, which are random fields. In ideal
waveguides the amplitudes are constant in range for $z > z'$
\begin{align} 
a_{j,o}^+(\om,\vx') &= \frac{\phi_j(x')}{2 i \beta_j(\om)},
\quad a_{j,o}^-(\om,\vx') = 0, \qquad ~~ ~ j = 1, \ldots, N,
\label{eq:st8} \\
 a_{j,o}^{e}(\om,\vx') &= -\frac{\phi_j(x')}{2
  \beta_j(\om)}, \quad ~ ~j > N.
\label{eq:st10}
\end{align}
They depend on the cross-range $x'$ in the support of the source, and
the second equation in (\ref{eq:st8}) complies with the wave being
outgoing. In random waveguides the mode amplitudes satisfy a coupled
system of stochastic differential equations driven by the random
fluctuations $\nu$, $\mu_{_B}$ and $\mu_{_T}$. They are analyzed in
detail in \cite{kohler77,garnier_papa,ABG-12,FRANCE} and the result is
that they are approximately the same as (\ref{eq:st8})-(\ref{eq:st10})
for range offsets $z-z' \ll \ep^{-2} \la_o$. This motivates the long
range scaling
\begin{equation}
z_{_\cA} = \ep^{-2} Z_{_\cA}, \qquad Z_{_\cA} = O(\la_o),
\label{eq:st11}
\end{equation}
where cumulative scattering becomes significant. The evanescent modes
may be neglected at such ranges\footnote{Note that although the
  evanescent modes do not appear explicitly in (\ref{eq:st12}), they
  affect the amplitudes of the propagating modes. This amplitude
  coupling is taken into account in the analysis in
  \cite{kohler77,garnier_papa,ABG-12,FRANCE} and thus in the results
  of this paper.}, and we use a further approximation that neglects
the backward going waves to write
\begin{align} 
p(t,\vx,\vx') \approx \int_{-\infty}^\infty \frac{d \om}{2 \pi B} \hat
f\left( \frac{\om-\om_o}{B}\right) e^{-i \om t} \sum_{j=1}^N
a_{j}^+(\om,z,\vx') \phi_j(x) e^{i \beta_j(\om)(z-z')}.
\label{eq:st12}
\end{align}
The forward scattering approximation holds for $\ell \gtrsim
\la_o$, and is justified by the fact that the backward mode amplitudes
have very weak coupling with the forward ones, for autocorrelations of
the fluctuations that are smooth enough in $z$
\cite{kohler77,garnier_papa,ABG-12,FRANCE}. We refer to
\cite{garnier2008effective} for the analysis of wave propagation that
includes both the forward and backward going modes, but for the
purpose of this paper it suffices to use (\ref{eq:st12}).

Let us write the amplitudes $a_j^+$ using the random propagator
$\PP^\ep \in \mathbb{C}^{N\times N}$, which maps the amplitudes
(\ref{eq:st8}) near the source at range $z'$, to those at the array
\begin{equation}
a_{j}^+\left(\om,\frac{Z_{_\cA}}{\ep^2},\vx'\right) = \sum_{l=1}^N
\PP^\ep_{jl}(\om,Z_\cA, z') a_{l,o}^+(\om,\vx').
\label{eq:st13}
\end{equation}
The propagator is analyzed in \cite{kohler77,garnier_papa,ABG-12} in
the asymptotic limit $\ep \to 0$. It converges in distribution to a
Markov diffusion $\PP$ with generator computed explicitly in
terms of the autocorrelations of the random fluctuations. Thus,
we can rewrite (\ref{eq:st13}) as
\begin{equation}
a_{j}^+\left(\om,\frac{Z_{_\cA}}{\ep^2},\vx'\right) \sim \sum_{l=1}^N
\PP_{jl}(\om,Z_\cA,z') a_{l,o}^+(\om,\vx'), 
\label{eq:st14}
\end{equation}
with symbol $\sim$ denoting approximate in distribution. It means that
we can approximate the statistical moments of $a_j^+$ using the right
hand side in (\ref{eq:st14}), with an $o(1)$ error in the limit $\ep
\to 0$. 

\subsection{Data model}
\label{sect:data}
The data model follows from (\ref{eq:fm4}), (\ref{eq:st1}) and
(\ref{eq:st12})
\begin{align}
D(t,x) \approx & \int_{\Omega_{\rho}} d \vx' \rho(\vx') \int_{-\infty}^\infty
\frac{du}{2 \pi} \, \hat \chi(u) \, e^{i u \frac{t_o}{\cT} }
\int_{-\infty}^\infty \frac{d \om}{2 \pi B} \hat
f\left(\frac{\om-\om_o}{B} - \frac{u}{B \cT} \right) e^{-i \om t}\times \nonumber
\\& \sum_{j=1}^N 1_{_\cA}(x) \phi_j(x) \, a_j^+
\left(\om-\frac{u}{\cT},\frac{Z_{_\cA}}{\ep^2},\vx'\right) e^{i
  \beta_j\left(\om-\frac{u}{\cT}\right) \left(\frac{Z_{_\cA}}{\ep^2} -
  z'\right)},
\label{eq:model}
\end{align}
where $\hat \chi$ is the Fourier transform of the recording window and
$a_j^+$ is given by (\ref{eq:st13})-(\ref{eq:st14}).  We take
henceforth the bandwidth
\begin{equation}
B = \om_o \ep^{\alpha}, \quad 1 < \alpha < 2,
\label{eq:model2}
\end{equation}
which is small with respect to the center frequency. We ask that
$\alpha < 2$ because the travel time of the modes is of order
$\ep^{-2}$, and we need a pulse of much smaller temporal support in
order to distinguish the arrival time of different modes. That $\alpha
< 2$ is also needed for the statistical stability of the inversion, as
we explain later. The choice $\alpha > 1$ is for
convenience\footnote{For $\alpha < 1$ the results are similar, but
  higher powers of $(\om - \om_o - u/\cT)$ enter in the phase, and
  they change the shape of the pulse carried by the modes.}, because
it allows us to linearize the phase in (\ref{eq:model}) as
\begin{equation*}
\beta_j\left(\om-{u}/{\cT} \right) \left(\ep^{-2}{Z_{_\cA}} -
z'\right) \approx \left[ \beta_j\left(\om_o\right) + \left( \om-\om_o
  - {u}/{\cT} \right)
  \beta_j'(\om_o)\right]\left(\ep^{-2}Z_{_\cA} - z'\right),
\end{equation*}
with small error of order $\ep^{2(\alpha -1)}$. When we use this
approximation in (\ref{eq:model}) we see that in ideal waveguides
where $a_j^+ = a_{j,o}^+$ the modes propagate with range speed
\begin{equation}
{1}/{\beta_j'(\om_o)} = c_o {\beta_j(\om_o)}/{k}.
\label{eq:model4}
\end{equation}
In random waveguides only the expectation (coherent part) of $a_j^+$
propagates at speed (\ref{eq:model4}), but the energy of the mode is
transported at different speed, as described in section
\ref{sect:I1}. In any case, we note that the wavenumbers $\beta_j$
decrease monotonically with $j$, so the first modes are faster as
expected, because they take a more direct path from the source to the
array. For example, in the case $N = \lfloor k X/\pi \rfloor \gg 1$
the slowness vectors $(\pm \pi /X,\beta_1)$ of the plane waves
associated with the first mode are almost parallel to the range
direction, and the speed (\ref{eq:model4}) is approximately equal to
$c_o$.  For the last modes the slowness vectors $(\pm \pi N
/X,\beta_N)$ are almost orthogonal to the range direction and the
speed is much smaller than $c_o$.
 
It is natural to choose the duration $\cT$ of the recording window to
be much longer than that of the pulse $\cT \gg 1/B$. We shall see in
section \ref{sect:transpIm} that in fact we need $\cT$ to be at least
of the order of the travel time of the waves in order for the
incoherent imaging method to work. Thus, we let
\begin{equation}
\cT = \ep^{-2} T,
\label{eq:durationT}
\end{equation}
with $T \ge O(1/\om_o)$.  We also assume that $\hat f$ is a continuous
function to simplify (\ref{eq:model}) slightly using the approximation
\begin{equation} 
\label{eq:model1}
\hat f\left(\frac{\om-\om_o}{B} - \frac{u}{B \cT} \right) = \hat f
\left(\frac{\om-\om_o}{B} - \ep^{2-\alpha} \frac{u}{\om_o T}\right)\approx
\hat f\left(\frac{\om-\om_o}{B} \right).
\end{equation}

\subsection{Loss of coherence} 
\label{sect:moments} 
To compute the coherent part of the data model, we recall from
\cite{kohler77,garnier_papa,ABG-12,FRANCE} the expectation of the
limit propagator
\begin{equation} 
\EE \left[ \PP_{jl}(\om,Z_{\cA},z')\right] \approx \delta_{jl} \exp
\left[ -\frac{Z_{_\cA}}{\cS_j(\om)} + i \frac{Z_{_\cA}}{\cL_j(\om)}\right],
\label{eq:st15}
\end{equation}
where $\delta_{jl}$ is the Kronecker delta symbol and the
approximation is due to the fact that $z'$ is much smaller than
$O(\ep^{-2}\la_o)$.  Although the mean propagator is a diagonal matrix as
in ideal waveguides, where it is the identity, its entries are
exponentially damped in $Z_{_\cA}$ on scales $\cS_j$, the scattering
mean free path of the modes. There is also an anomalous phase
accumulated on the mode dependent scales $\cL_j$.

The scales $\cS_j$ and $\cL_j$ are defined in \cite[equations
  (3.19),(3.28),(3.31)]{BGT-14} and depend on the frequency and the
autocorrelations $\cR_\nu$, $\cR_{_B}$ and $\cR_{_T}$ of the
fluctuations. Of particular interest in this paper are the scattering
mean free paths because they give the range scale on which the modes
randomize. The magnitude of the expectation (coherent part) of the
mode amplitudes follows from (\ref{eq:st14}) and (\ref{eq:st15})
\begin{equation}
\left|\EE\left[a_j^+\left(\om,\ep^{-2}{Z_{_\cA}},
  \vx'\right)\right]\right| = e^{-\frac{Z_{_\cA}}{\cS_j}} \left|
a_{j,o}^+(\om,\vx') \right|, \qquad j = 1, \ldots, N, 
\label{eq:st16NN}
\end{equation}
where $a_{j,o}^+$ is the initial condition of $a_j^+(\om,z,\vx')$ at
$z = z'$, equal to the amplitude (\ref{eq:st8}) in ideal waveguides.
The exponential decay in (\ref{eq:st16NN}) is not caused by
attenuation in the medium. The wave equation conserves energy, and we
state in the next section that $\EE \left[|a_j^+|^2\right]$ does not
tend to zero. The meaning of the decay in (\ref{eq:st16NN}) is the
randomization (loss of coherence) of the $j-$th mode due to
scattering. It says that beyond scaled ranges $Z_{_\cA} > \cS_j$ the
mode becomes incoherent i.e., the random fluctuations of its amplitude
dominate its expectation.

The scattering mean free paths are given by 
\begin{equation}
\cS_j(\om) = \frac{2}{\sum_{q =1}^N \Gamma_{jq}^{(c)}(\om)},
\label{eq:defSj}
\end{equation}
in terms of 
\begin{align}
\Gamma_{jq}^{(c)}(\om) = & \frac{\pi^4 \ell (jq)^2}{
  \beta_j(\om)\beta_q(\om) X^4} \left\{ \hat \cR_{_B}\left[\ell
  \left(\beta_j(\om)-\beta_q(\om)\right)\right] + \hat
\cR_{_T}\left[\ell \left(\beta_j(\om)-\beta_q(\om)\right) \right]
\right\} + \nonumber \\ &\frac{k^4 \ell }{4 \beta_j(\om) \beta_q(\om)}
\hat \cR_{\nu_{jq}}\left[\ell
  \left(\beta_j(\om)-\beta_q(\om)\right)\right],
\label{eq:st22}
\end{align}
where $\cR_{\nu_{jq}}$ is the autocorrelation of the stationary
process
\begin{align}
\nu_{jq}\left(\frac{z}{\ell}\right) &= \int_0^X dx \,
\nu\left(\frac{x}{\ell},\frac{z}{\ell}\right) \phi_j(x) \phi_q(x) .
\end{align}
This is for $\ep_c = \ep_{_B} = \ep_{_T} = \ep$ in (\ref{eq:fm5}) and
(\ref{eq:fm7}), and for statistically independent random processes
$\nu$, $\mu_{_B}$ and $\mu_{_T}$. The hat denotes the Fourier
transform of the autocorrelations, which is non-negative by Bochner's
theorem.

To compare the scattering effects in the random medium with those at
the boundary, we plot $\cS_j$ with solid line in Figure
\ref{fig:MatrixU} for the case $\ep_{_B} = \ep_{_T} = 0$ and Figure
\ref{fig:MatrixUBound} for $\ep_c = \ep_{_B}= 0$.  In the first case
we keep only the last term in (\ref{eq:st22}) and in the second case
we keep the second term. The setup of the simulations is explained in
the numerics section \ref{sect:I3}. Here we note two important facts
displayed by the plots: The scales $\cS_j$ decrease monotonically with
$j$, and their mode dependence is much stronger in the random boundary
case.  This is intuitive once we recall that the first modes are waves
that travel along a more direct path from the source to the
array. These waves interact with the random boundary only once in a
while and thus randomize on longer range scales than the slower
modes. For example $\cS_1$ is more than a hundred times longer than
$\cS_{20}$ in Figure \ref{fig:MatrixUBound}. The slow modes are waves
that reflect repeatedly at the boundary and travel a long way in the
waveguide as they progress slowly in range. They randomize on small
range scales for both random boundary and medium scattering. However,
the medium scatter leads to more dramatic loss of coherence as
illustrated in Figure \ref{fig:MatrixU}, where all but the last modes
have similar scattering mean free paths which are shorter than in
Figure \ref{fig:MatrixUBound}.

The goal of this paper is to analyze what can be determined about the
source of waves from measurements made at ranges $Z_{_\cA} > \cS_1$,
where
\begin{equation}
\left|\EE\left[a_j^+\left(\om,\ep^{-2}{Z_{_\cA}},
  \vx'\right)\right]\right| \le e^{-\frac{Z_{_\cA}}{\cS_1}} \left|
a_{j,o}^+(\om,\vx') \right| \ll \left| a_{j,o}^+(\om,\vx') \right|
\label{eq:st16}
\end{equation}
for all $j = 1, \ldots, N$ i.e., all the modes are incoherent. No
coherent method can work in this regime, so we study an incoherent
inversion approach based on the transport of energy theory summarized in
the next two sections.

\subsection{Statistical decorrelation}
\label{sect:decorr}
Since the wave equation is not dissipative, we have the conservation
of energy relation \cite{kohler77,garnier_papa,ABG-12,FRANCE}
\begin{equation}
\sum_{j=1}^N \left|a_j^+\left(\om,\ep^{-2}{Z_{_\cA}},
\vx'\right)\right|^2 \approx \sum_{j=1}^N \left| a_{j,o}^+(\om,\vx')
\right|^2,
\label{eq:st17}
\end{equation}
where the approximation is with an $o(1)$ error as $\ep \to 0$, due to
the neglect of the backward going and evanescent waves. Thus, some
second moments of the mode amplitudes remain finite, and can be used
in inversion. To decide if we can estimate them reliably from the
incoherent data, we need to know how the waves
decorrelate. Statistical decorrelation means that the second moments
of the amplitudes are equal approximately to the product of their
expectations, which is negligible by (\ref{eq:st16}).
 
The two frequency analysis of the propagator $\PP^\ep$ is carried out
in \cite{garnier_papa,ABG-12}, and the result is that the waves are
decorrelated for frequency offsets $ |\om-\om'| \ge O(\ep^2 \om_o).  $
Such small offsets are enough to cause the waves to interact
differently with the random fluctuations over ranges $\ep^{-2}Z_\cA$,
thus giving the statistical decorrelation. This result is important
because it says that we can estimate those second moments of the
amplitudes that do not decay in range by cross-correlating the Fourier
transform of the data at nearby frequencies $\om$ and $\om -\ep^2 h$
and integrating over $\om \in (\om_o-\pi B,\om_o+\pi B)$ to obtain a
statistically stable result. The bandwidth $B$ is much larger than
$\ep^2 \om_o$ by assumption (\ref{eq:model2}), and the statistical
stability follows essentially from a law of large numbers, because we
sum a large number of terms that are uncorrelated.

The second moments of the propagator at nearby frequencies are
\begin{align}
\EE \left[ \PP_{jl}^\ep (\om,Z_{_\cA},z) \overline{\PP_{j'l'}^\ep
    (\om- \ep^2 h, Z_{_\cA},z') } \right] \approx \delta_{jl}
\delta_{j'l'} \frac{\beta_l(\om)}{\beta_j(\om)} \hat
\cW_j^{(l)}(\om,h,Z_{_\cA}) e^{-i \beta'_j(\om) h Z_{_\cA}} +
\nonumber \\ (1-\delta_{jj'}) \EE \left[\PP_{jl}^\ep
  (\om,Z_{_\cA},z)\right] \overline{\EE\left[\PP_{j'l'}^\ep
    (\om,Z_{_\cA},z')\right] } e^{Z_{_\cA}/\cL_{jj'}}, \qquad \quad
\label{eq:st18}
\end{align}
where the bar denotes complex conjugate, $\hat \cW_j^{(l)}$ is the
Fourier transform of the Wigner distribution described below, and the
scale $\cL_{jj'}$ is defined in terms of the autocorrelations
$\cR_\nu$, $\cR_{_B}$ and $\cR_{_T}$ (see \cite[equation
  (6.26)]{FRANCE}). These formulas follow from the calculations in
\cite{garnier_papa,ABG-12} which assume $z = z'$, and the law of
iterated expectation with conditioning at $z$, for $z < z'$.  Denoting
by $\EE_{z}$ the conditional expectation and using
\[\EE_z\left[\PP_{j'l'}^\ep(\om-\ep^{2}h,Z_{\cA},z')\right]
\approx \PP_{j'l'}^\ep(\om-\ep^{2}h,Z_{\cA},z), \] because $z'-z \ll
\ep^{-2} \la_o$,  we obtain
\begin{align*}
\EE \left[ \PP_{jl}^\ep (\om,Z_{_\cA},z) \overline{\PP_{j'l'}^\ep
    (\om- \ep^2 h, Z_{_\cA},z') } \right] &= \EE \left[ \PP_{jl}^\ep
  (\om,Z_{_\cA},z) \,
  \EE_{z} \hspace{-0.03in}\left[\overline{\PP_{j'l'}^\ep (\om- \ep^2
      h, Z_{_\cA},z') } \right] \right] \\&\approx \EE \left[
  \PP_{jl}^\ep (\om,Z_{_\cA},z) \overline{\PP_{j'l'}^\ep (\om- \ep^2
    h, Z_{_\cA},z) } \right]
\end{align*}
and (\ref{eq:st18}) follows from \cite{garnier_papa,ABG-12} and the
fact that in the support of the source $z \ll \ep^{-2}\la_o$.

The last term in (\ref{eq:st18}) corresponds to the coherent part of
the mode amplitudes and it is negligible in our regime with $Z_{_\cA}
> \cS_1$. This is by (\ref{eq:st15}) and 
\[
Z_{_\cA}\left[\frac{1}{\cS_j} +\frac{1}{\cS_j'} -
  \frac{1}{\cL_{jj'}}\right] \sim \frac{Z_{_\cA}}{\cS_1} >1.
\]
Recalling the expression (\ref{eq:st13}) of the mode amplitudes in
terms of the propagator, we see that (\ref{eq:st18}) states that the
amplitudes of different modes are essentially uncorrelated. Therefore,
the only second moments that remain large are the mean energies of the
modes, which is why we use them in inversion.

\subsection{The system of transport equations}
\label{sect:transpeq}
The Wigner distribution defines the expectation of the energy of the
$j-$th mode resolved over a time window of duration similar to the
travel time, when the initial excitation is in the $l-$th mode. It
satisfies the following system of transport equations derived in
\cite{kohler77,garnier_papa,ABG-12}
\begin{equation}
\left[\partial_Z + \beta_j'(\om) \partial_\tau\right]
\cW_j^{(l)}(\om,\tau,Z) = \sum_{q = 1}^N \Gamma_{jq}(\om)
\cW_q^{(l)} (\om,\tau,Z), \quad Z > 0,
\label{eq:st19}
\end{equation}
with initial condition 
\begin{equation}
\cW_j^{(l)}(\om,\tau, 0) = \delta_{jl} \delta (\tau),
\label{eq:st20}
\end{equation}
where $\delta (\tau)$ is the Dirac delta distribution. The Fourier
transform that appears in (\ref{eq:st18}) is defined by
\begin{equation}
\hat \cW_j^{(l)}(\om,h,Z) = \int_{-\infty}^\infty d \tau \,
\cW_j^{(l)}(\om,\tau,Z) e^{i h \tau} = \left[ e^{\left(i h \cB'(\om) +
    \Gamma(\om)\right) Z}\right]_{jl},
\label{eq:st20FT}
\end{equation}   
where $\cB'$ is the diagonal matrix 
\begin{equation}
\cB'(\om) = {\rm diag}\left(\beta_1'(\om), \ldots, \beta_N'(\om)\right).
\label{eq:cB}
\end{equation}

The matrix $\Gamma(\om)$ in (\ref{eq:st19}) models the transfer
of energy between the modes, due to scattering. Its off-diagonal
entries are defined in (\ref{eq:st22})
\begin{equation}
\Gamma_{jq}(\om) = \Gamma_{jq}^{(c)}(\om), \qquad j \ne q, 
\end{equation}
and are non-negative, meaning that there is an outflow of energy from
mode $j$ to the other modes.  The energy lost by this mode is
compensated by the gain of energy in the other modes, as stated by
\begin{equation}
\Gamma_{jj}(\om) = - \sum_{q \ne j} \Gamma_{jq}(\om),
\qquad \forall \, j = 1, \ldots, N.
\label{eq:st21}
\end{equation}

\section{Inversion based on energy transport equations}
\label{sect:transpIm}
We now use the results summarized above to formulate our inversion
approach.  We give in section \ref{sect:formod} the forward model
which maps the source density to the cross-correlations of the mode
amplitudes. These are defined in section \ref{sect:dataproc} and are
self-averaging with respect to different realizations of the random
waveguide. Therefore, we can relate them to the Wigner
distribution. The inversion method is studied in section
\ref{sect:inverse}.

\subsection{Data processing}
\label{sect:dataproc}
The first question that arises is how to relate the incoherent array
data to the moments (\ref{eq:st18}) of the propagator which are
defined by the Wigner distribution.  The answer lies in computing
cross-correlations of the data projected on the eigenfunctions
$\phi_j$, as we now explain.

We denote by $\hat D(\om,x)$ the Fourier transform of the measurements
and by $\hat D_j(\om)$ its projection on the eigenfunction $\phi_j$
 \begin{equation}
\hat D_j(\om) = \int_0^X dx \, \hat D(\om,x) \phi_j(x).
\label{eq:inv1}
\end{equation}
We are interested in its cross-correlation $\hat \cC_j(h)$ at lag
$\ep^2 h$ and its inverse Fourier transform $\cC_j(\tau)$. The latter
has the physical interpretation of energy carried by the $j-$th mode
over the duration of a time window which we model with a bump function
$\psi$ of dimensionless argument and order one support
\begin{equation}
\cC_j(\tau) = \frac{2 \pi H}{\ep^2} \int_{-\infty}^\infty dt \, \psi(Ht)
\left|D_j\left(\frac{\tau-t}{\ep^2}\right)\right|^2.
\label{eq:inv3}
\end{equation}
Here $H$ has units of frequency, satisfying $HT \gg 1$, so the
integrand is compactly supported in the recording window $\chi$. The
scaling by $\ep^{-2}$ of the argument of $D_j$, the inverse Fourier
transform of (\ref{eq:inv1}), is to be consistent with the $O(\ep^{-2}
Z_{_\cA}/c_o)$ travel time of the waves to the array, and the factors
in front of the integral are chosen to get an order one 
\begin{align}
\hat \cC_j(h) = \hat \psi \left(\frac{h}{H} \right)
\int_{-\infty}^\infty d \om \, \hat D_j(\om) \overline{\hat
  D_j(\om-\ep^2 h)}.
\label{eq:inv2}
\end{align}
This expression is obtained by taking the inverse Fourier transform of
(\ref{eq:inv3}), and the integral over $\om$ is restricted by the
support of $\hat D_j(\om)$ to $|\om-\om_o| \le \pi B$.

We relate below the expectation of $\cC_j(\tau)$ to the Wigner
distribution, and explain in Appendix \ref{ap:A} under which
conditions $\cC_j(\tau)$ is self-averaging, meaning that it is
approximately equal to its expectation. The self-averaging is due to
the rapid frequency decorrelation of $\hat D_j(\om)$ over intervals of
order $\ep^2 \om_o$, and the bandwidth assumption
(\ref{eq:model2}). When we divide the frequency interval $(\om_o - \pi
B, \om_o + \pi B)$ in smaller ones of order $\ep^2 \om_o$, we see that
in (\ref{eq:inv2}) we are summing a large number $
{B}/({\ep^2 \om_o}) = \ep^{\alpha-2} \gg 1 $ of uncorrelated random
variables. The self-averaging is basically by the law of large
numbers, as long as $\EE[\cC_j(\tau)]$ is large. This happens for
large enough arrays, for long recording times that scale as
(\ref{eq:durationT}), and for times $\tau$ near the peak $\tau_j$ of
$\cC_j$.

The role of the projection (\ref{eq:inv1}) is to isolate in the data
the effect of the $j-$th mode. We see from
(\ref{eq:model})-(\ref{eq:model1}) that
\begin{align}
\hat D_j(\om) \approx & \frac{1}{B} \hat
f\left(\frac{\om-\om_o}{B}\right) \sum_{q=1}^N Q_{jq}
\int_{-\infty}^\infty \frac{du}{2 \pi } \, \hat \chi (u) e^{i u
  \left[\ep^2 t_o - \beta_q'(\om_o)Z_{_\cA}\right]/{T}} \times
\nonumber \\ & \int_{\Omega_\rho} d \vx' \, \rho(\vx')
a_q^+\left(\om-\frac{\ep^2 u}{T},\frac{Z_{_\cA}}{\ep^2}, \vx'\right)
e^{i \left[\beta_q(\om_o) +
    (\om-\om_o)\beta_q'(\om_o)\right]\left(\frac{Z_{_\cA}}{\ep^2} -
  z'\right)},
\label{eq:inv4}
\end{align}
where we introduced the mode coupling matrix $Q \in \mathbb{R}^{N
  \times N}$ with entries
\begin{equation}
Q_{jq} = \int_0^X dx \, 1_{_\cA}(x) \phi_j(x) \phi_q(x). 
\label{eq:inv5}
\end{equation}
This coupling is an effect of the aperture of the array. The ideal
setup is for an array with full aperture $\cA = [0,X]$, because $Q$ is
the identity by the orthonormality of the eigenfunctions, and $\hat
D_j$ involves only the amplitude of the $j-$th mode.  However, all the
mode amplitudes enter the expression of $\hat D_j$ when the array has
partial aperture, and they are weighted by $Q_{jq}$. The coupling
matrix is diagonally dominant when the length of the aperture $|\cA|$
is not much smaller than the waveguide depth $X$. This can be seen for
example in the case of an array starting at the top boundary $\cA =
[X-|\cA|, X]$, where
\begin{equation}
Q_{jq} = \delta_{jq} - \left(1-\frac{|\cA|}{X}\right)
\left\{ \begin{array}{ll}1 - {\rm sinc}\left( \frac{2 \pi j
      (X-|\cA|)}{X}\right) , \quad & q = j ,
  \\ {\rm sinc}\left(
    \frac{\pi (j+q) (X-|\cA|)}{X}\right) - {\rm sinc}\left( \frac{\pi
      (j-q) (X-|\cA|)}{X}\right) & j \ne q.
\end{array} \right. 
\label{eq:inv6}
\end{equation}

We note in (\ref{eq:inv4}) that by choosing the support $\cT$ of the
recording window $\chi$ as in (\ref{eq:durationT}), we can relate
$\hat D_j(\om)$ to the mode amplitudes in a frequency interval of
order $\ep^2 \om_o$. This is important in the calculation of the
cross-correlations $\hat \cC_j(h)$, where the amplitudes must be
evaluated at nearby frequencies.  If we had a smaller $\cT$, the
cross-correlations would involve products of the amplitudes at
frequency offsets that exceed $\ep^2 \om_o$. Such amplitudes are
statistically uncorrelated and there is no benefit in calculating the
cross-correlation.

\subsection{The forward model} 
\label{sect:formod}
We show in Appendix \ref{ap:A} that 
\begin{align}
\EE \left[ \cC_j(\tau) \right] \approx \frac{\|f\|^2}{4 B}\left| \chi
\left( \frac{\tau-\ep^2 t_o}{T} \right) \right|^2 \sum_{q, l =1}^N
Q^2_{jq} \frac{\left| \hat \rho_l\left[ \beta_q(\om_o)\right]
  \right|^2}{\beta_l(\om_o) \beta_q(\om_o)} \times \nonumber \\ \int
\frac{d h}{2 \pi} \hat \psi \left(\frac{h}{H}\right) \left[ e^{(i h
    \cB'(\om_o) + \Gamma(\om_o))Z_{_\cA}}\right]_{ql} e^{-i h
  \tau},
\label{eq:A4}
\end{align}
where the diagonal matrix $\cB'$ defined in (\ref{eq:cB}) is evaluated
at $\om_o$, 
\begin{equation}
\hat \rho_{l}(\beta) = \int_{\Omega_\rho} d \vx \, \rho(\vx) \phi_l(x)
e^{-i \beta z}
\label{eq:A2}
\end{equation}
are the Fourier coefficients of the unknown source density, and $
\|f\|^2 = \int_{-\infty}^\infty du \, |\hat f(u)|^2$.  Because the
cross-correlations are self-averaging we can define the forward map
$\mathfrak{F}$ from $\rho$ to the vector
$\left(\cC_j(\tau)\right)_{1\le j \le N},$ using equation
(\ref{eq:A4}).  We write it as
\begin{align}
\label{eq:FMod}
\left[\mathfrak{F}(\rho)\right]_j(\tau) = \frac{\|f\|^2}{4 B}
 \sum_{q, l =1}^N Q^2_{jq} \frac{
   \left| \hat \rho_l\left[ \beta_q\right] \right|^2}{\beta_l \beta_q}
 \int \frac{d h}{2 \pi} \hat \psi \left(\frac{h}{H}\right) \left[
   e^{(i h \cB' + \Gamma)Z_{_\cA}}\right]_{ql} e^{-i h \tau},
\end{align}
which is a simplification of (\ref{eq:A4}) based on the assumption
that the recording window $\chi$ is well centered and sufficiently
long to equal one at the times of interest.  We also simplify the
notation by dropping the $\om_o$ argument of the wavenumbers
$\beta_q$, their derivatives $\beta_q'$ and $\Gamma$.  The
unknown source density appears in the model as the $N\times N$
matrix of absolute values of its Fourier coefficients
(\ref{eq:A2}). This is the most that we can expect to recover from the
inversion.

\section{Inversion}
\label{sect:inverse}
We have the following unknowns: the range $Z_{_\cA}$, the $N \times N$
matrix $\left(|\hat \rho_l(\beta_q)|\right)_{1\le q,l\le N}$, and
possibly the autocorrelations of the fluctuations. The question is
what can be recovered from $\left(\cC_j(\tau) \right)_{1\le j \le N}$
and how to carry the inversion.  The range $Z_{_\cA}$ and some
information about the autocorrelation of the fluctuations can be
determined from the measurements of the travel times $\tau_j$ of
$\cC_j(\tau)$. This is the easier part of the inversion and we discuss
it first. The estimation of $\rho$ is more delicate and requires
knowing $Z_{_\cA}$ and the autocorrelation of the fluctuations, so we
can calculate the matrix $\Gamma$. We discuss it in sections
\ref{sect:I2}-\ref{sect:I6}. We illustrate the results with numerical
simulations in section \ref{sect:I3}.

\subsection{Arrival time analysis} 
\label{sect:I1}
If there were no random scattering effects i.e., no matrix
$\Gamma$, the $h$ integral in (\ref{eq:FMod}) would equal
$\delta_{ql} H \psi \left[H (\tau - \beta_q' Z_{_\cA})\right]$. This
implies in particular that for an array with full aperture, where $Q$
equals the identity, the cross-correlation $\cC_j(\tau)$ would have a
single peak at the travel time $\tau = \beta_j' Z_{_\cA}$. In random
waveguides the transport speed is not $1/\beta_j'$. The matrices
$\cB'$ and $\Gamma$ in the exponential in (\ref{eq:FMod}) do not
commute, so there is anomalous dispersion due to scattering which must
be taken into account in inversion.

The range estimation based on arrival (peak) times of $\cC_j(\tau)$
was studied with numerical simulations in \cite[Section 6.1]{BIT-10}
for the case of a point source. The method there uses definition
(\ref{eq:st20FT}) of the Wigner transform for a search range
$Z_{_{\cA}}^{s}$, and estimates $Z_{_\cA}$ as the minimizer of the
misfit between the peak time of the theoretical model (\ref{eq:FMod})
and the calculated $\left(\cC_j(\tau)\right)_{1\le j\le N}$ from the
data. It is observed in \cite{BIT-10} that the range estimation is not
sensitive to knowing the source density and that the search for
$Z_{_\cA}$ can be done in conjunction with the estimation of the
autocorrelation of the fluctuations, in case it is unknown. The method
in \cite{BIT-10} has been tested extensively with numerical
simulations for both large and small arrays in waveguides with random
wave speed.  The conclusion is that the estimation of $Z_{_\cA}$ is
very robust, but the success of the estimation of $\cR_\nu$ depends on
having the right model of the autocorrelation. For example, with a
Gaussian model of a Gaussian $\cR_\nu$, the optimization determines
correctly the correlation length $\ell$.  For another model the
optimization returns the wrong correlation length, but the range
$Z_{_\cA}$ is still well determined.  This is because the anomalous
dispersion depends on $\Gamma$, which is defined by (\ref{eq:st22}) in
terms of only a few Fourier coefficients of the autocorrelation
function. There are many functions that give the same Fourier
coefficients i.e., the same $\Gamma$, so to get the true correlation
length we need the true model of $\cR_\nu$.

Here we complement the results in \cite{BIT-10} with an explicit
arrival time analysis which can be carried out using perturbation
theory. We explain in Appendix \ref{ap:B} that in forward scattering
regimes, as assumed in this paper, the matrix $i h \cB'$ may be
treated as a perturbation of $\Gamma$. Thus, we can approximate the
matrix exponential in (\ref{eq:FMod}) using the perturbation of the
spectral decomposition of $\Gamma$. By definition $\Gamma$ is
symmetric, so it has real eigenvalues $\Lambda_j$ and eigenvectors
${\bf u}_j$ for $j = 1, \ldots, N$ that form an orthonormal basis of
$\mathbb{R}^N$. The eigenvalues satisfy $\Lambda_j \le 0$, otherwise
the energy would not be conserved (recall (\ref{eq:st17})), and the
null space of $\Gamma$ is nontrivial, since by (\ref{eq:st21})
\begin{equation}
\label{eq:defu1}
\Gamma {\bf u}_1 = 0, \qquad \mbox{where} ~ ~ {\bf u}_1 = (1, 1,
\ldots, 1)^T/\sqrt{N}.
\end{equation}
We count henceforth the eigenvalues in decreasing order, and suppose
they are distinct. This assumption is not needed for the inversion to
work, and we use it only in this section. It allows a simpler arrival
time analysis, because we can approximate the spectrum of $i h \cB' +
\Gamma$ with regular perturbation theory.

If we denote by $\Lambda_j(h)$ the eigenvalues and ${\bf u}_j(h)$ the
eigenvectors of  $i h \cB' + \Gamma$, we
have the standard results \cite{golub2012matrix}
\begin{equation}
\Lambda_j(h) \approx \Lambda_j + i h {\bf u}_j^T \cB' {\bf u}_j,
\label{eq:PertL}
\qquad {\bf u}_j(h) \approx {\bf u}_j + i h \sum_{q \ne j} \frac{{\bf
    u}_q^T \cB' {\bf u}_j}{\Lambda_j-\Lambda_q} {\bf u}_q.
\end{equation}
Thus, we approximate the matrix exponential by
\begin{equation}
e^{\left( i h \cB' + \Gamma\right) Z_{_\cA}} \approx \sum_{j =
  1}^N  e^{\left(\Lambda_j + i h {\bf u}_j^T \cB'
  {\bf u}_j\right) Z_{_\cA}}\, {\bf u}_j {\bf u}_j^T,
\label{eq:apM1}
\end{equation}
where we neglect the perturbation of the eigenvectors because it has
little influence on the arrival times.  Substituting (\ref{eq:apM1})
in the forward model (\ref{eq:FMod}), we obtain that
\begin{align}
\label{eq:FModPert}
\left[\mathfrak{F}(\rho)\right]_j(\tau)\approx \frac{H \|f\|^2}{4 B}
\sum_{r=1}^N e^{-|\Lambda_r| Z_{_\cA}} \, \Psi \left( H (\tau - Z_{_\cA}
    {\bf u}_r^T \cB' {\bf u}_r)\right) \sum_{q,l =1}^N
    Q^2_{jq} \frac{ \left| \hat \rho_l\left[ \beta_q\right]
      \right|^2}{\beta_l \beta_q} u_{qr} u_{lr},
\end{align}
where $u_{qr}$ is the $q$ component of the eigenvector ${\bf u}_r$.
This is a superposition of $N$ pulses (bumps) $\psi$ traveling at
transport speed
\begin{equation}
V_r = \frac{1}{{\bf u}_r^T \cB' {\bf u}_r} = \left[\sum_{q=1}^N
  \beta_q' u_{qr}^2\right]^{-1}.
\label{eq:AVGSpeed}
\end{equation}
Only the first term in (\ref{eq:FModPert}) does not decay in range,
and travels at speed\footnote{This equation is also derived in
  \cite[Section 20.6.2]{book07} using a probabilistic analysis of the
  transport equations (\ref{eq:st19}).}
\begin{equation}
V_1 = N \left(\sum_{q=1}^N \beta_q'\right)^{-1}.
\label{eq:AVGSpeed1}
\end{equation}
The other terms decay exponentially and their transport speeds $V_r$
are quite different than $1/\beta_r'$, unless the entries in ${\bf
  u}_r$ are concentrated around the $r-$th row.

We illustrate in Figure \ref{fig:TRANSP} the transport speeds
$\left(V_r\right)_{1 \le r \le N}$ calculated for two types of random
waveguides: filled with a random medium and with a random top
boundary. The setup is discussed in detail in the numerics section
\ref{sect:I3}, and the spectrum of $\Gamma$ is displayed in Figures
\ref{fig:MatrixU} and \ref{fig:MatrixUBound}.  Figure \ref{fig:TRANSP}
shows that the difference between $V_r$ and $1/\beta'_r$, which
quantifies the anomalous dispersion, depends on the ratio $\ell/\la_o$
and the type of scattering: in the medium or at the boundary.

The number of terms contributing in (\ref{eq:FModPert}) depends on the
array aperture via the coupling matrix $Q$, and the magnitude of the
entries in the eigenvectors ${\bf u}_r$.  We study in the next section
the structure of the matrix $(u_{qr})_{1 \le q,r \le N}$ and explain
that it has a nearly vanishing block in the upper right corner. This
is also illustrated in Figures \ref{fig:MatrixU} and
\ref{fig:MatrixUBound}. The implication is that the $r$ index of
summation in (\ref{eq:FModPert}) extends roughly up to $j$, so there
are more terms to sum for the slower modes than the fast ones. Thus,
at moderate ranges we expect a wider spread in $\tau$ of $\cC_j(\tau)$
for large $j$. As $Z_{_\cA}$ grows, only the first term $r = 1$ contributes,
and the arrival time becomes independent of $j$ 
\begin{align}
\label{eq:FModPertLim}
\left[\mathfrak{F}(\rho)\right]_j(\tau)\stackrel{\eps \to
 0}{\longrightarrow} \frac{H \|f\|^2}{4 B}
\Psi \left( H (\tau - Z_{_\cA} /V_1)\right) \hspace{-0.05in}\sum_{q,l
  =1}^N Q^2_{jq} \frac{ \left| \hat \rho_l\left[ \beta_q\right]
  \right|^2}{\beta_l \beta_q}.
\end{align}

Note from (\ref{eq:FModPert}) and (\ref{eq:FModPertLim}) that the
arrival (peak) time is mostly dependent on the spectral decomposition
of $\Gamma$, and not on the actual source density $\rho$, which
only changes the ``weights'' of the bump $\psi$.  The unlikely case
where the last sum in (\ref{eq:FModPert}) equals zero is taken into
account in \cite{BIT-10} by excluding from the optimization the modes
with small values of the calculated $\cC_j$. Consequently, the
estimation of $Z_{_\cA}$ is insensitive to the lack of knowledge of
$\rho$, as observed in \cite{BIT-10}.

\subsection{Estimation of the source density} 
\label{sect:I2}
We suppose henceforth that $Z_{_\cA}$ has been determined and that the
autocorrelations of the fluctuations are either known or have been
estimated as explained in the previous section in sufficient detail to
be able to approximate $\Gamma$.

Because the time $\tau$ does not appear in the $\rho$ dependent factor
in (\ref{eq:FMod}) or (\ref{eq:FModPert}), it suffices to consider the
peak values of $\cC_j(\tau_j)$ as the inversion data or alternatively,
to integrate $\cC_j(\tau)$ over $\tau$. We choose the latter because
it is more robust, and define the column vector $\bM \in \mathbb{R}^N$
of newly processed data with entries
\begin{align}
\mathfrak{M}_j: = \frac{4 B}{\|f\|^2 \hat \psi(0)}
\int_{-\infty}^\infty d \tau \,\cC_j(\tau) &\approx
\sum_{q, l =1}^N Q^2_{jq} \frac{ \left| \hat \rho_l\left[
    \beta_q\right] \right|^2}{\beta_l \beta_q} \left[
  e^{\Gamma Z_{_\cA}}\right]_{ql} \nonumber \\ &= \sum_{r=1}^N
e^{-|\Lambda_r| Z_{_\cA}}\sum_{q,l =1}^N Q^2_{jq} \frac{ \left| \hat
  \rho_l\left[ \beta_q\right] \right|^2}{\beta_l \beta_q} u_{qr}
u_{lr},
\label{eq:inve1}
\end{align}
where $\hat \psi(0) = H \int d\tau \, \psi(H \tau).$ We only have $N$
data so we cannot expect to determine uniquely the $N\times N$ matrix
with entries $|\hat \rho_l(\beta_q)|^2$, unless we have additional
assumptions on $\rho$. For example, in \cite{BIT-10} it is assumed
that the source has small, point-like support. Here we let instead
$\rho(\vx)$ be a separable function
\begin{equation}
\rho(\vx) = \xi(x) \zeta(z),
\label{eq:inve2}
\end{equation}
so that 
\begin{equation}
\hat \rho_l(\beta) = \hat \xi_l \, \hat \zeta(\beta), \qquad \hat
\xi_l = \int_0^X dx \, \xi(x) \phi_l(x), \qquad \hat \zeta(\beta) =
\int_{-\infty}^\infty dz \,  \zeta(z)  e^{-i \beta z},
\label{eq:inve3}
\end{equation}
and we can study separately the estimation of the range and
cross-range profiles of the source. Such separation is usual in
imaging, where the range is determined from the arrival time of the
waves and the cross-range from their direction of arrival. We used the
arrival times $\ep^{-2} \tau_j$ to determine the distance $\ep^{-2}
Z_{_\cA}$ from the source to the array.  We cannot get more
information from them because the cross-correlations are at $O(\ep^2
H)$ frequency lag, which means that the error in the arrival time
estimation is $O(\ep^{-2}/H)$.  If we do not know anything about
$\rho(\vx)$, we can only assume that the source is tightly supported
at distance $\ep^{-2} Z_{_\cA}$ from the array (i.e., let $\zeta(z) =
\delta(z)$), and estimate the cross-range profile $\xi(x)$. Only if we
know $\xi(x)$ we can estimate $\zeta(z)$.

Let us write (\ref{eq:inve1}) in vector form
\begin{equation}
\bM \approx \mathbb{Q} \, {\rm diag}\left(|\hat \zeta(\beta_1)|^2,
\ldots, |\hat \zeta(\beta_N)|^2\right) \cB^{-1}\sum_{r=1}^N {\bf u}_r
      {\bf u_r^T} \cB^{-1} \left( \begin{array}{c} |\hat \xi_1|^2
        \\ \vdots \\ |\hat \xi_N|^2
\end{array} \right) e^{-|\Lambda_r|Z_{_\cA}},
\label{eq:inve5}
\end{equation}
where $\mathbb{Q}$ is the matrix with entries $Q_{jq}^2$ and $\cB =
\diag( \beta_1, \ldots, \beta_N)$.  We have two cases:
\vspace{-0.05in}
\begin{enumerate}
\itemsep 0.05in
\item Invert for the range profile $\zeta(z)$ when $\xi(x)$ is known.
\item Invert for the cross-range profile $\xi(x)$ when $\zeta(z)$ is
  approximately $\delta(z)$.
\end{enumerate}
\vspace{0.05in} We analyze both cases under the assumption that
$\mathbb{Q}$ is strictly diagonally dominant and therefore
invertible. This holds for a large enough aperture $\cA$.

\vspace{0.1in}
\noindent \textbf{Case 1} When we know the cross-range profile
$\xi(x)$  we can calculate the vector
\begin{equation}
\bdeta = e^{\Gamma Z_{_\cA}} \cB^{-1} \left( \begin{array}{c}
  |\hat \xi_1|^2 \\ \vdots \\ |\hat \xi_N|^2
\end{array} \right) =  \sum_{r=1}^N  e^{-|\Lambda_r|Z_{_\cA}} {\bf u}_r
      {\bf u_r^T} \cB^{-1} \left( \begin{array}{c} |\hat \xi_1|^2
        \\ \vdots \\ |\hat \xi_N|^2
\end{array} \right),
\label{eq:inve7}
\end{equation}
to rewrite equation (\ref{eq:inve5}) as
\begin{equation}
\mathbb{Q}^{-1} \bM \approx {\rm diag}\left(|\hat \zeta(\beta_1)|^2,
\ldots, |\hat \zeta(\beta_N)|^2\right) \cB^{-1} \bdeta,
\label{eq:inve6}
\end{equation}
and invert it by
\begin{equation}
|\hat \zeta(\beta_j)|^2 \approx \frac{\beta_j \left(\mathbb{Q}^{-1}
  \bM\right)_j}{\eta_j}, \qquad {\rm if} ~ ~ \eta_j \ne 0.
\label{eq:inve8}
\end{equation}
We know that the matrix exponential has a trivial null space, 
so the vector $\bdeta$ cannot be zero, but can some of its components 
be zero or very small? 

To answer this question let us decompose $\bdeta$ in two orthogonal
parts: one that lies in ${\rm Null}(\Gamma)$ and is constant in
range, and the other that lies in ${\mathbb R}^N \setminus {\rm
  Null}(\Gamma)$ and decays exponentially in range. To be more
precise, suppose henceforth that the null space is one dimensional
\begin{equation}
{\rm Null}(\Gamma) = {\rm span}\{{\bf u}_1\}, 
\label{eq:null1d}
\end{equation}
and therefore $\Lambda_2 < 0$. A sufficient (not necessary) condition
for this to hold is that all the off-diagonal entries of
$\Gamma$ are strictly positive, which happens for 
autocorrelation functions like Gaussians for example. Then
$\Gamma$ is a matrix of Perron-Frobenius type, and its largest
eigenvalue $\Lambda_1$ is simple. Equation (\ref{eq:inve7}) gives 
\begin{align}
\bdeta &= {\bf u}_1 {\bf u}_1^T \cB^{-1} \left(\begin{array}{c} |\hat
  \xi_1|^2 \\ \vdots \\ |\hat \xi_N|^2
\end{array} \right) + {\boldsymbol{\mathcal E}} = \frac{1}{N} 
\left(\sum_{j=1}^N \frac{|\hat \xi_j|^2}{\beta_j} \right)
\left(\begin{array}{c} 1 \\ \vdots
  \\ 1 \end{array} \right) + {\boldsymbol{\mathcal E}}
\end{align}
with residual vector ${\boldsymbol{\mathcal E}}$ that decays in range
like $\exp(-|\Lambda_2| Z_{_\cA})$. Thus, all the components of
$\bdeta$ are bounded below by a positive constant as $Z_{_\cA}$ grows,
and the calculation (\ref{eq:inve8}) is well-posed.

\vspace{0.1in} \noindent \textbf{Case 2} When the source has
point-like support in range we let $\hat \zeta(\beta) \approx 1$ in
(\ref{eq:inve5}) and invert the system as
\begin{equation}
\left(\begin{array}{c}|\hat \xi_1|^2 \\ \vdots \\ |\hat
  \xi_N|^2 \end{array} \right) \approx \cB \bX 
\label{eq:inve9}
\end{equation}
where 
\begin{equation}
\bX = e^{-\Gamma Z_{_\cA}} \cB \mathbb{Q}^{-1} \bM = \sum_{j=1}^N
e^{|\Lambda_j| Z_{_\cA}}\left({\bf u_j^T} \cB
\mathbb{Q}^{-1} \bM\right)  {\bf u}_j .
\label{eq:inve10}
\end{equation}
However, this calculation is ill-posed due to the exponential growth
in $Z_{_\cA}$ of the right hand side, so we need regularization. There
are many ways to regularize, and the inversion can be improved with
prior information about $\xi(x)$.  Here we discuss a spectral cut-off
regularization which uses the first $J$ terms in (\ref{eq:inve10})
\begin{equation}
\bX_J = \sum_{j=1}^J e^{|\Lambda_j| Z_{_\cA}} \left( {\bf u}_j^T \cB
\mathbb{Q}^{-1} \bM \right){\bf u}_j.
\label{eq:XJ}
\end{equation}
This is the orthogonal projection of $\bX$ on the subspace spanned by
$\{{\bf u}_1, \ldots, {\bf u}_J\}$ or, equivalently, the minimum
Euclidian norm vector that gives a misfit of order ${\rm exp}
\left(-|\Lambda_{J+1}|Z_{_\cA}\right)$ between the data
(\ref{eq:inve5}) and the model.

But in what sense does $\bX_J$ approximate $\bX$ and therefore the
vector of absolute values of the Fourier coefficients of $\xi$?  We
expect that it should be easier to estimate $|\hat \xi_j|$ for lower
indices $j$ that correspond to the fast modes which have less
interaction with the random fluctuations than the slow modes.  To see
if this is the case, note first from (\ref{eq:inve9}) that since $\cB$
is diagonal, it is sufficient to investigate if $\bX_J$ approximates
better the first components of $\bX$.  Let the orthogonal projector
operator be $\mathbb{U}_J$, so that $\bX_J = \mathbb{U}_J \bX$. The
error can be bounded as
\begin{equation}
\label{eq:errEst}
\frac{|(\bX-\bX_J)_j|}{\|\bX\|} = \frac{\|{\bf e}_j^T
  (I-\mathbb{U}_J)\bX \|}{\|\bX\|} \le \|(I -\mathbb{U}_J){\bf e}_j\|
= \sqrt{\sum_{q=J+1}^N u^2_{jq}},
\end{equation}
and it is guaranteed to be small for $1 \le j \lesssim J$ if the
eigenvectors ${\bf u}_q$ for $q \ge J+1$ have small entries in the
first $J$ rows.  Here $I$ is the $N \times N$ identity matrix and
${\bf e}_j$ are the vectors of the canonical basis in $\mathbb{R}^N$.
We demonstrate in sections \ref{sect:I3} and \ref{sect:I4} with
numerical simulations and with analysis that indeed, the matrix ${\bf
  U} = ({\bf u}_1, \ldots, {\bf u}_N)$ of eigenvectors of $\Gamma$ has
a nearly vanishing block in the upper right corner. Thus, we expect a
good approximation of the first $J$ entries in $\bX$ if $Z_{_\cA}
\lesssim 1/|\Lambda_J|$.

\subsection{Estimation of $\rho$ from the absolute value of its Fourier transform}
\label{sect:discuss}
Given that we can only estimate a few absolute values of the Fourier
coefficients of the cross-range (range) profile of the source, what
can we actually say about the source density? Clearly, it is
impossible to reconstruct $\rho$ in detail unless we have prior
knowledge. Otherwise we get limited information such as its support.
Here are a few examples:

\vspace{0.05in}
\noindent $\bullet$ \textbf{Point like source.} If we let $\rho(\vx) =
\delta(x-x_\star) \delta(z)$, it is enough to determine the absolute
value of the first Fourier coefficient
\[
| \hat \rho_1(\beta)| \approx |\phi_1(x_\star)|, \qquad \forall \beta.
\]
Since $|\phi(x)|$ is monotonically increasing for $x \in [0,X/2)$ and
  decreasing for $x \in (X/2,0]$, this gives the cross-range location
$x_\star$ up to a reflection with respect to the axis of the
waveguide. This reflection ambiguity cannot be resolved by estimating
higher order Fourier coefficients of $\rho$.  It is due to the
symmetric boundary conditions at $x=0$ and $x = X$. If we had
Dirichlet conditions at $x = X$ and Neumann at $x = 0$, $|\phi_1(x)|$
would be monotone in $(0,X)$ and $x_\star$ would be uniquely
determined by $|\hat \rho_1(\beta)|$. We discuss next a more robust
way of estimating the support of the source.

\vspace{0.05in}
\noindent $\bullet$ \textbf{Size of cross-range support.} Let us
denote by $\xi_{e}(x)$ the odd extension of the cross-range profile of
the source about $x = 0$, and define its autocorrelation
\begin{equation}\label{autocorrelation}
\mathfrak{R}_{_\xi}(x) = \int_{-X}^X dx' \, \xi_e(x') \xi_e(x'+x) = 2
\sum_{j=1}^\infty |\hat \xi_j|^2 \cos \left(\frac{\pi j x}{X}\right),
\end{equation}
where the last equality follows by direct calculation using the
Fourier sin series expansion of the real valued $\xi_e(x)$. Obviously,
we can approximate $\mathfrak{R}_{_\xi}(x)$ using the regularized
solution described in Case 2 of the previous section, if the Fourier
coefficients $\hat \xi_j$ are small for $j > J$. Otherwise, we get the
autocorrelation of a smoothed version of the source. 
To illustrate what we can expect, suppose that  
\[
\xi(x) = \mathcal{N}(x-x_o,\sigma), \quad \mbox{where} \quad
\mathcal{N}(x;\sigma)=\tfrac{1}{\sqrt{2\pi}\,\sigma}e^{-\frac{x^{2}}{2\sigma^{2}}},
\] 
and $\sigma \ll X$ so that the essential support of the Gaussian is
inside the interval $(0,X)$. Then
$\xi_{e}(x)=\mathcal{N}(x- x_o,\sigma) - \mathcal{N}(x+x_o,\sigma)$,
and the autocorrelation is given by 
\begin{equation}\label{gaussian_auto}
\mathfrak{R}_{_\xi}(x) \approx 2\,\mathcal{N}(x;\sqrt{2}\,\sigma) -
\mathcal{N}(x-2x_o;\sqrt{2}\,\sigma) -
\mathcal{N}(x+2x_o;\sqrt{2}\,\sigma).
\end{equation}
The first term in \eqref{gaussian_auto} is invariant to
translations of the source, and can be used to estimate the
cross-range support of the source (i.e., $\sigma$).  The remaining two
terms depend on the source location, and can be used to estimate
$x_o$.  Because the autocorrelation is a $2X$-periodic function, the
translation by $2x_o$ in \eqref{gaussian_auto} is understood modulo
$2X$.  Consequently, sources that are symmetrically located about the
center of the waveguide ($x=X/2$) produce the same
autocorrelation. That is to say, the location $x_m$ of the minimum of
the autocorrelation determines the center of the source up to a
reflection ambiguity: at $x_o = x_m/2$ or at its reflection $x_o = X -
x_m/2$.  We illustrate the estimation of $\xi(x)$ with numerical
simulations in Figure \ref{fig:inversion}.

\vspace{0.05in}
\noindent $\bullet$ \textbf{Size of range support.} The
autocorrelation of the range profile is 
\[
\mathfrak{R}_{_\zeta}(z) = \int_{-\infty}^\infty dz' \, \zeta(z')
\zeta(z'+z) = \frac{1}{\pi} \int_{0}^\infty d \beta \, |\hat
\zeta(\beta)|^2 \cos(\beta z),
\]
where we used that $\zeta(z)$ is real valued. We can approximate
$\mathfrak{R}_{_\zeta}$ from $\{|\hat \zeta(\beta_j)|\}_{1 \le j \le
  N}$ when $N \gg 1$, so that $\beta_j$ sample well the interval
$(0,k)$, and $|\hat \zeta(\beta)| \ll 1$ for $\beta > \beta_1 \approx
k$. We already know that the source is centered at $z = 0$, and the
size of the support of $\zeta(z)$ follows from that of
$\mathfrak{R}_{_\zeta}(z)$ as above.

\subsection{The equipartition regime and the benefit of a large bandwidth}
\label{sect:I6}
We saw in the previous sections that the accuracy of the cross-range
estimation depends on how $Z_{_\cA}$ compares to the scales
$1/|\Lambda_j|$. We refer to Figures \ref{fig:MatrixU} and
\ref{fig:MatrixUBound} for an illustration of these scales and note
that while in waveguides with random boundaries $\cS_1 \approx
1/|\Lambda_2|$, in waveguides filled with random media there is a gap
between $\cS_1$ and $1/|\Lambda_2|$ of at least one order of
magnitude. The importance of the scale $1/|\Lambda_2|$ is revealed
once we calculate from (\ref{eq:st20FT}) and (\ref{eq:defu1}) the mean
energy carried by a mode
\[
\int_{-\infty}^\infty d \tau \, \cW_j^{(l)}(\om_o,\tau,Z) =
\left[e^{\Gamma Z_{_\cA}}\right]_{jl} = \sum_{r=1}^N
e^{-|\Lambda_r| Z_{_\cA}} u_{jr} u_{lr} \approx \frac{1}{N},
\]
where the approximation is for $Z_{_\cA} > {1}/{|\Lambda_2|}$ and all
$j, l = 1, \ldots, N$.  Cumulative scattering distributes the energy
uniformly over the modes, which is why
\begin{equation}
\cL_{eq} = 1/|\Lambda_2|
\end{equation}
is called the equipartition distance. The waves forget their initial
direction when they travel further than $\ep^{-2} \cL_{eq}$, and the
processed data (\ref{eq:inve5}) becomes approximately
\begin{equation}
\bM \approx \frac{1}{N} \left[\sum_{j=1}^N \frac{|\hat
    \xi_j|^2}{\beta_j}\right] \mathbb{Q} \left( \begin{matrix} |\hat
  \zeta(\beta_1)|^2/\beta_1 \\ \vdots \\ |\hat
  \zeta(\beta_N)|^2/\beta_N \end{matrix} \right).
\end{equation}
It depends only on the weighted average of $\left(|\hat
\xi_j|^2\right)_{1\le j \le N}$, so the cross-range profile estimation 
(Case 1 in section \ref{sect:I2}) is impossible.

Because $\cS_1 \approx \cL_{eq}$ in waveguides with random boundaries,
coherent inversion with mode filtering as in \cite{BGT-14} is the best
approach for estimating the cross-range profile $\xi(x)$ of the
source. That method fails at ranges that exceed $\ep^{-2} \cS_1$,
where all the modes are incoherent, but since the waves are in the
equipartition regime, it is impossible to determine $\xi(x)$ with any
other method. In waveguides filled with random media there is a range
interval between $\ep^{-2} \cS_1$ and $\ep^{-2} \cL_{eq}$ where
incoherent inversion based on the cross-correlations $\cC_j$ can
determine approximately $\xi(x)$. Thus, we may say that the incoherent
method analyzed in this paper is more useful in these
waveguides. However, all this is for narrow bandwidths, scaled as in
(\ref{eq:model2}). For large bandwidths we may be able to improve the
inversion, as we now explain.

Assuming a large bandwidth of the signal emitted by the source, let
us divide it in smaller sub-bands scaled as in (\ref{eq:model2}),
centered at frequencies $\om_j$ listed in increasing order, for $j= 1,
\ldots, M$. Definition (\ref{eq:st22}) and the relation $\pi/X \approx
k/N$ show that the magnitude of $\Gamma$ grows with the frequency, so
we expect the least scattering effects in the lower frequency band
centered at $\om_1$.  If it is the case that $Z_{_\cA} \lesssim
1/|\Lambda_J(\om_1)|$ for some $J>1$, then we can invert as in section
\ref{sect:I2}, and recover roughly $|\hat \xi_j|$ for $1 \le j
\lesssim J$. However, when $Z_{_\cA} > 1/|\Lambda_2(\om_1)|$, the
waves are in the equipartition regime throughout the whole frequency
range, and all we can determine from each sub-band are the weighted
averages
\begin{equation}
\theta_j = \frac{1}{N_j} \sum_{q=1}^{N_j} \frac{|\hat
  \xi_q|^2}{\beta_q(\om_j)}, \qquad N_j:= N(\om_j), \quad j = 1, \ldots, M.
\label{eq:AVGxi}
\end{equation}
Combining the results we obtain the linear system
\begin{equation}
\mathbb{B} \left(\begin{matrix} |\hat \xi_1|^2 \\ \vdots \\ |\hat
  \xi_{_{N_{_M}}}|^2 \end{matrix} \right) = \left(\begin{matrix}
  N_1 \theta_1 \\ \vdots \\ N_{_M} \theta_{_M}\end{matrix} \right)
\label{eq:LINSYST}
\end{equation}
with $M \times N(\om_{_M})$ matrix $\mathbb{B}$ with rows equal to 
\begin{equation}
{\bf e}_j^T \mathbb{B} = \left\{ \begin{array}{ll}\left(
  1/\beta_1(\om_j) , \ldots, 1/\beta_{N_j}(\om_j), 0 \ldots,
  0\right), \qquad &1 \le j < M \\ \\\left( 1/\beta_1(\om_{_M}) , \ldots,
  1/\beta_{N_{_M}}(\om_{_M})\right), & j = M.
\end{array} \right.
\end{equation}
Direct calculation shows that most of the rows in $\mathbb{B}$ are
linearly independent at frequency separation $|\om_j - \om_q| =
O(\om_o)$ for $j \ne q$, so it is possible to improve the estimation
of the cross-range profile of the source for large enough $M$. In
particular, when $M = N_{_M}$ we can determine uniquely the solution
from (\ref{eq:LINSYST}).  We refer to Figure \ref{fig:estim} for a
numerical illustration of the improvement brought by a wide bandwidth
in the estimation of the cross-range location of a point-like source.

\subsection{Numerical simulations} 
\label{sect:I3}
To illustrate the theoretical results of the previous sections, we
present here numerical simulations for two types of random
waveguides. The first has flat boundaries and random wave speed with
Gaussian autocorrelation of the fluctuations $\nu$
\begin{equation}
\cR_\nu\left(\frac{x}{\ell},\frac{x'}{\ell},\frac{z}{\ell}\right) =
\EE\left[ \nu\left(\frac{x}{\ell},\frac{z}{\ell}\right)
  \nu\left(\frac{x'}{\ell},0\right)\right] = e^{-\frac{(x-x')^2}{2
    \ell^2} - \frac{z^2}{2 \ell^2}}.
\label{eq:ex1}
\end{equation}
The second is for a waveguide filled with a homogeneous medium and 
random top boundary with Gaussian autocorrelation of the fluctuations 
$\mu_{_T}$, 
\begin{equation}
\cR_{_T}\left(\frac{z}{\ell}\right) = \EE\left[
  \mu_{T}\left(\frac{z}{\ell}\right) \mu_{T}(0)\right] = e^{-
  \frac{z^2}{2 \ell^2}}.
\label{eq:ex1B}
\end{equation}
Using these in the definition (\ref{eq:st22}) we obtain that in the
first case
\begin{equation}
\Gamma_{jq} \approx \frac{\pi}{X} \frac{\ell^2 k_o^4}{\beta_j
  \beta_q} e^{-\frac{\ell^2}{2} (\beta_j-\beta_q)^2} \left[
  e^{-\frac{(k_o \ell)^2}{2} \frac{(j-q)^2}{N^2}} + e^{-\frac{(k_o
      \ell)^2}{2} \frac{(j+q)^2}{N^2}}\right], \qquad j \ne q, 
\label{eq:ex2}
\end{equation}
where $k_o = \om_o/c_o$, and the approximation is for $\ell \ll X$.
In the second case we have 
\begin{equation}
\Gamma_{jq} = \frac{\pi^4 \sqrt{2 \pi} \ell
  (jq)^2}{\beta_j \beta_q X^4} e^{-\frac{\ell^2}{2} 
  (\beta_j-\beta_q)^2}, \qquad j \ne q.
\label{eq:ex2B}
\end{equation}
We take $c_o = 1.5$km/s, the sound speed in water, the wavelength
$\la_o = 1.5$m corresponding to central frequency $1$kHz, $X = 20.3
\la_o$, so that $N = 40$, and three choices of the correlation length:
$\ell = \la_o$, $\ell = 3 \la_o$ and $\ell = 5 \la_o$.  

\begin{figure}[t]
\begin{center}
\includegraphics[width=0.43 \textwidth]{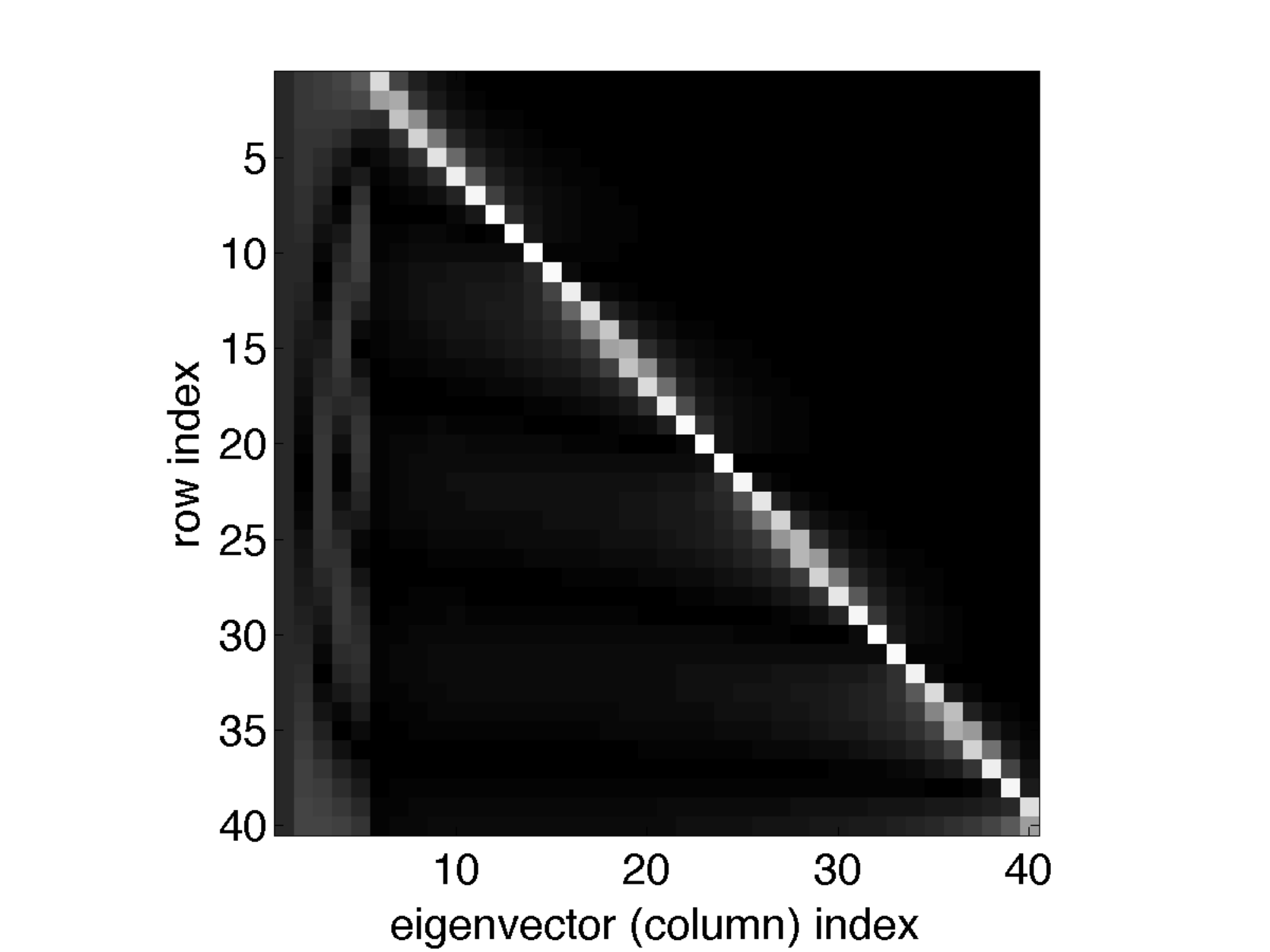} 
\raisebox{0.1in}{\includegraphics[width=0.35 \textwidth]{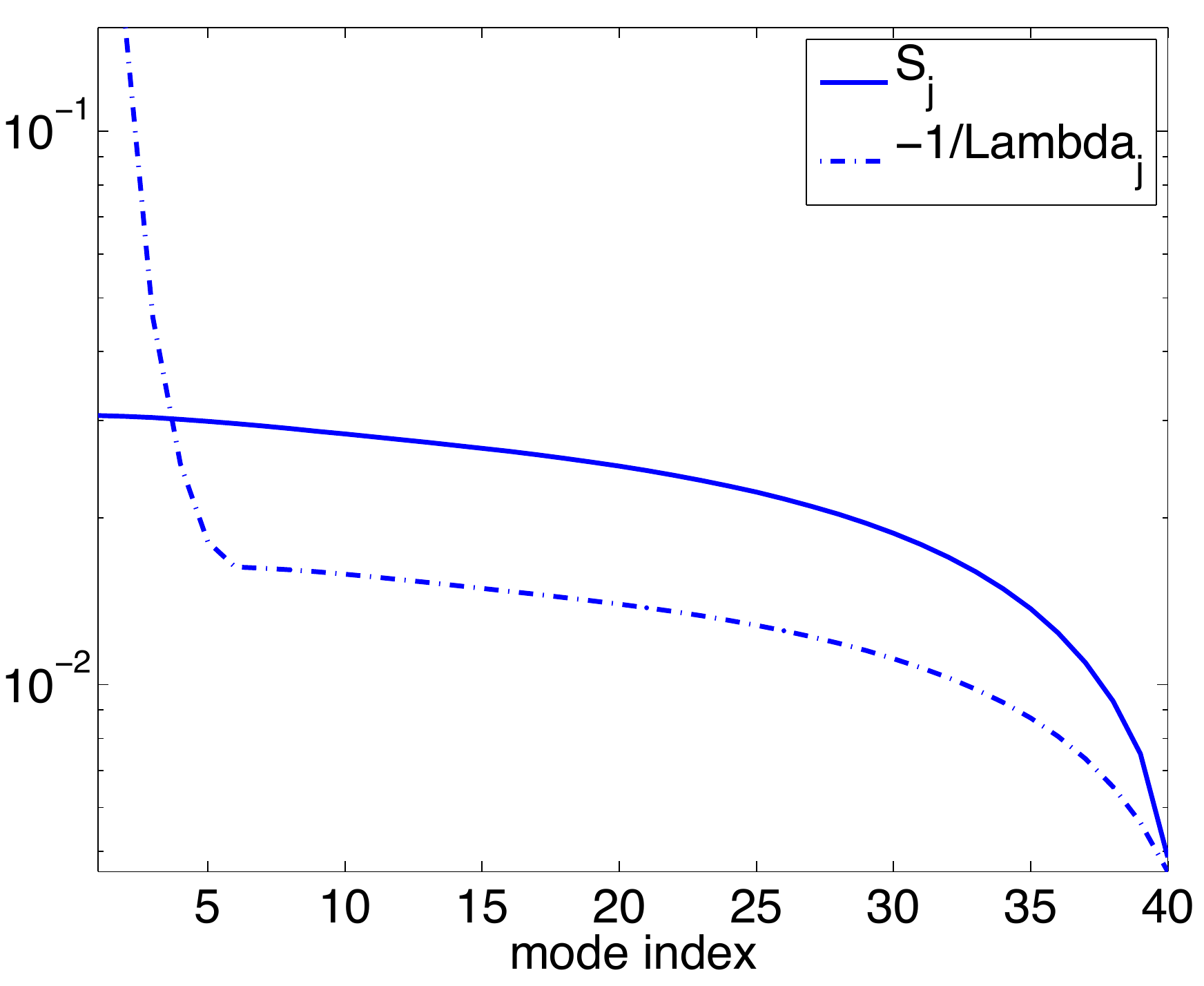}} 

\vspace{-0.03in}
\includegraphics[width=0.43 \textwidth]{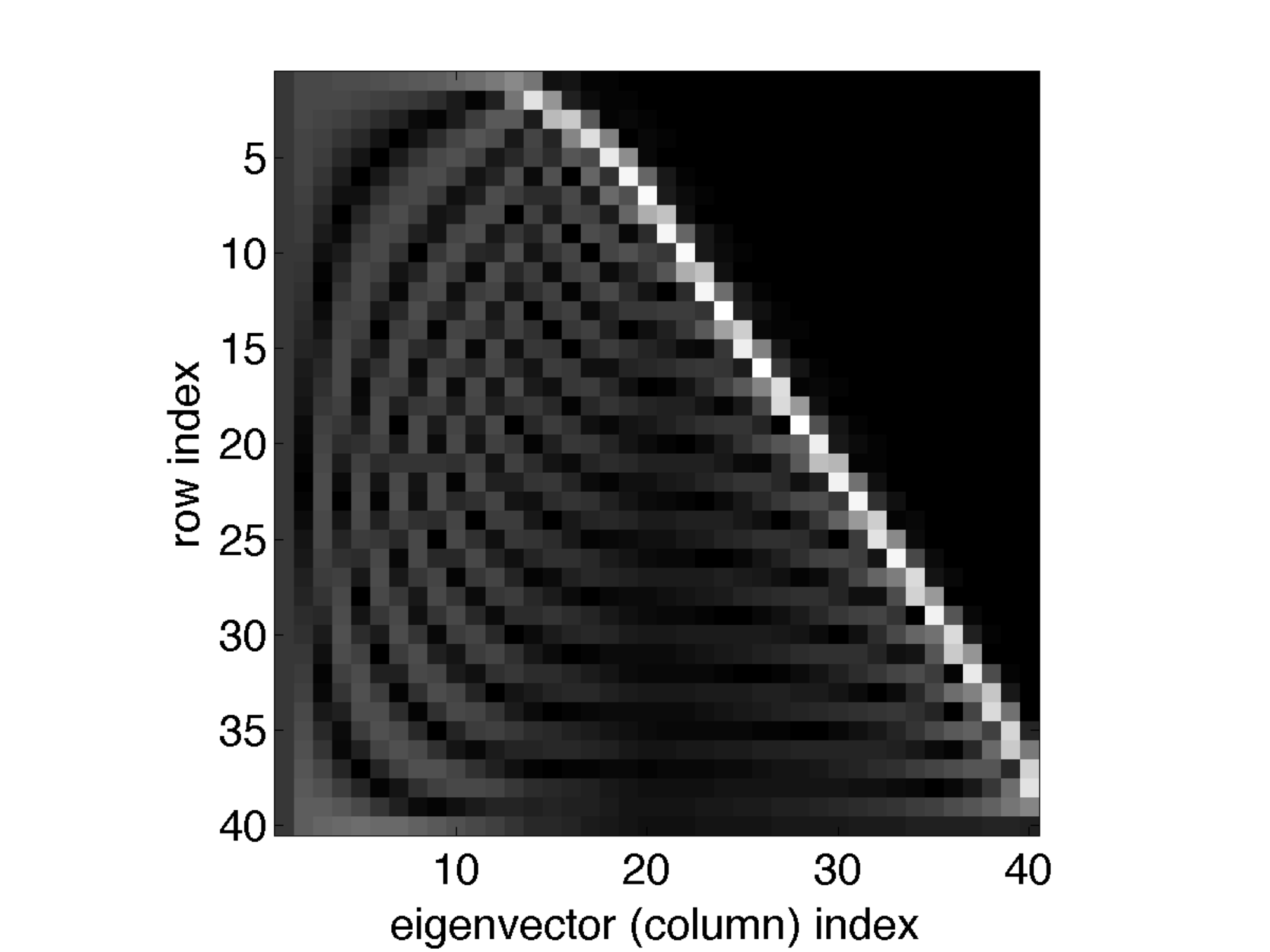} 
\raisebox{0.1in}{\includegraphics[width=0.35 \textwidth]{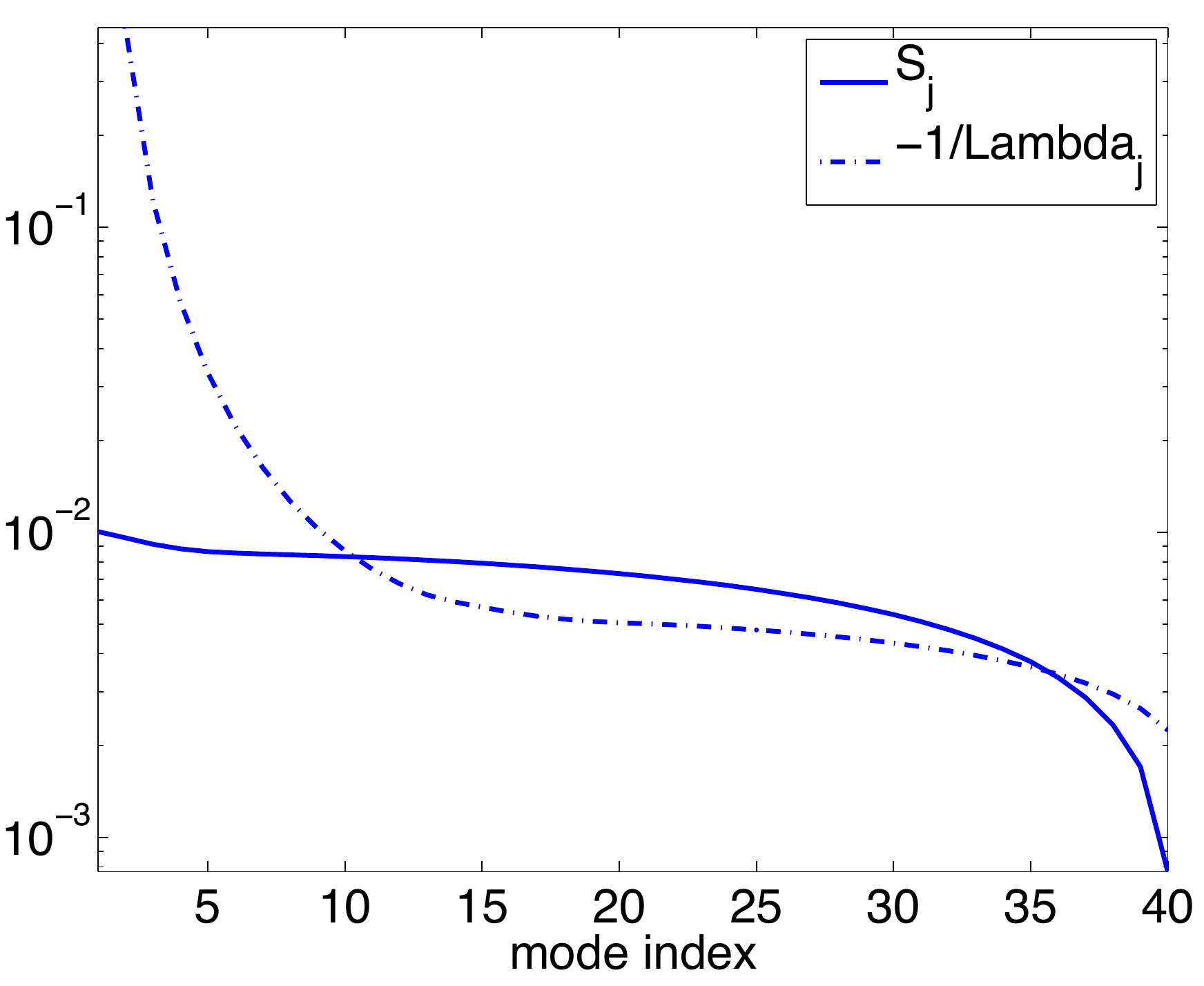}} 

\vspace{-0.03in}
\includegraphics[width=0.43 \textwidth]{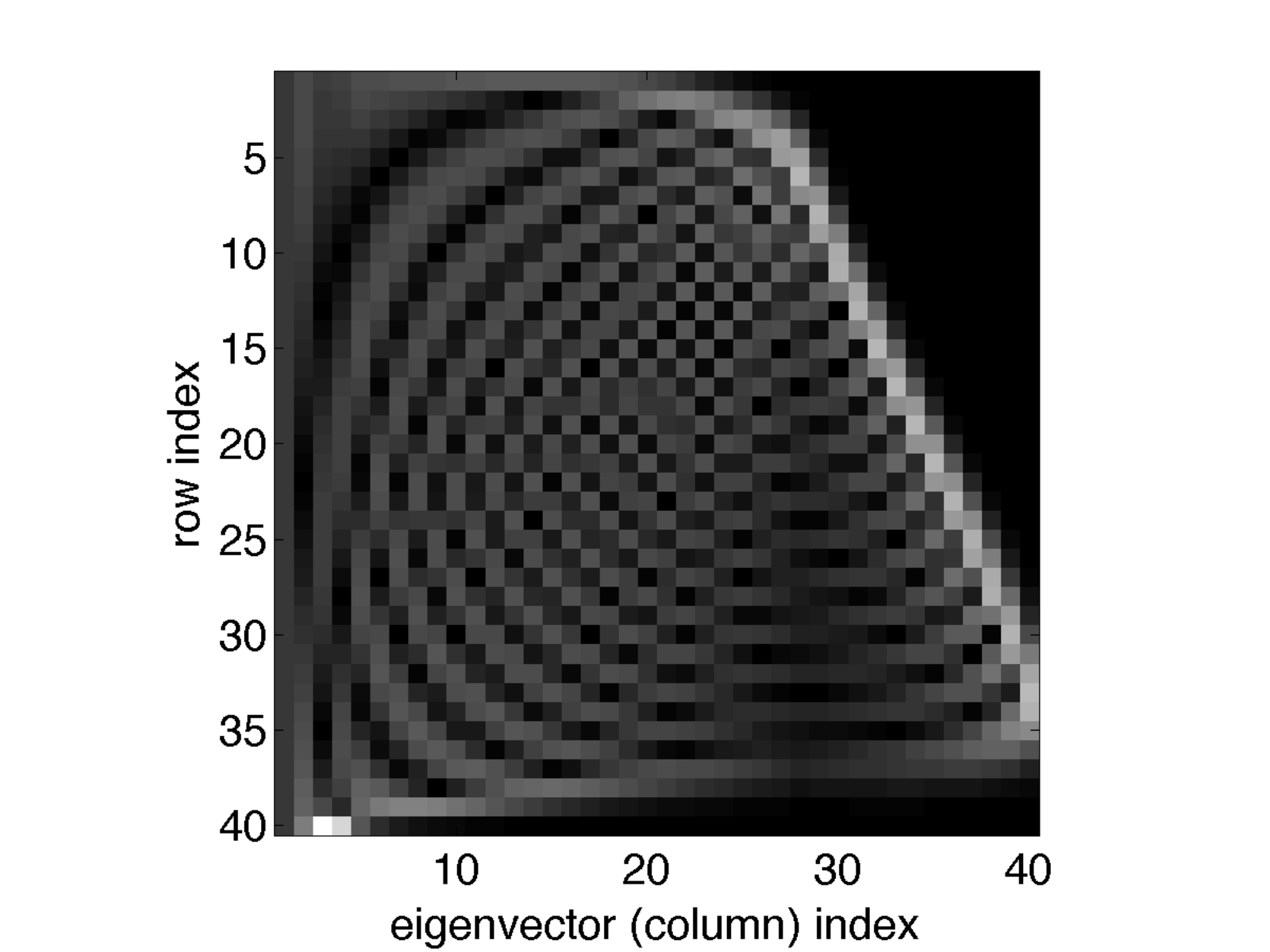} 
\raisebox{0.1in}{\includegraphics[width=0.35 \textwidth]{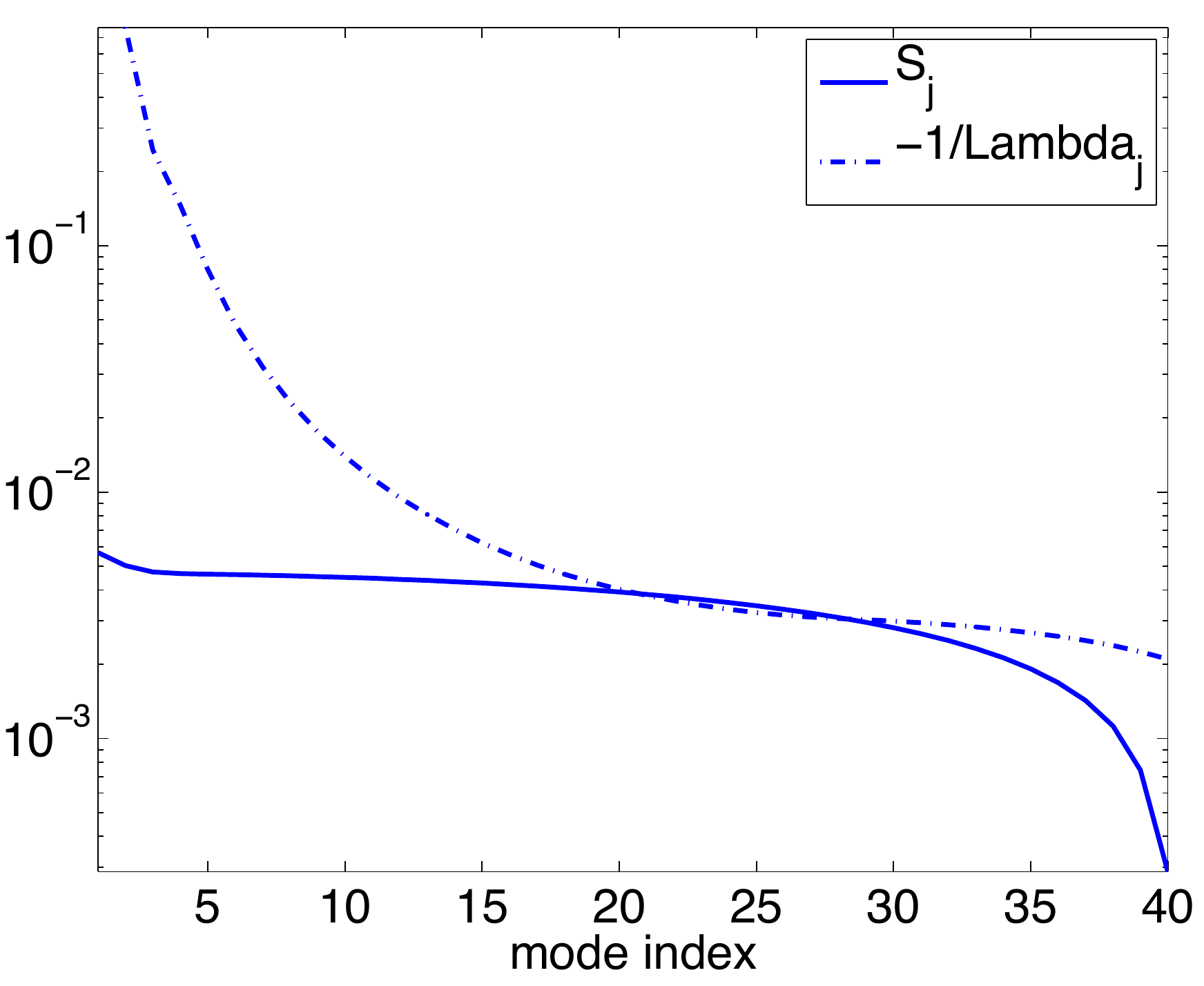}}
\end{center}
\vspace{-0.1in}
\caption{Waveguide filled with a random medium. Left plots: Absolute
  values of the entries of the matrix ${\bf U}$ of eigenvectors. Gray
  scale with lighter color indicates larger values and black indicates
  nearly zero. The column index is in the abscissa and the row index
  in the ordinate. Right plots: The scattering mean free path of the
  modes (full line) and the scales $-1/\Lambda_j$, for $j = 2, \ldots
  N$ (dotted line). In the top row $\ell = \la_o$, in the middle row
  $\ell = 3 \la_o$ and in the last row $\ell = 5 \la_o$. These scales
  should be multiplied by $\ep^{-2}$ e.g., in media with $1\%$
  fluctuations, the ordinate is in units of $\times 10^4$m.}
\label{fig:MatrixU}
\end{figure}
\begin{figure}[h]
\begin{center}
\includegraphics[width=0.43 \textwidth]{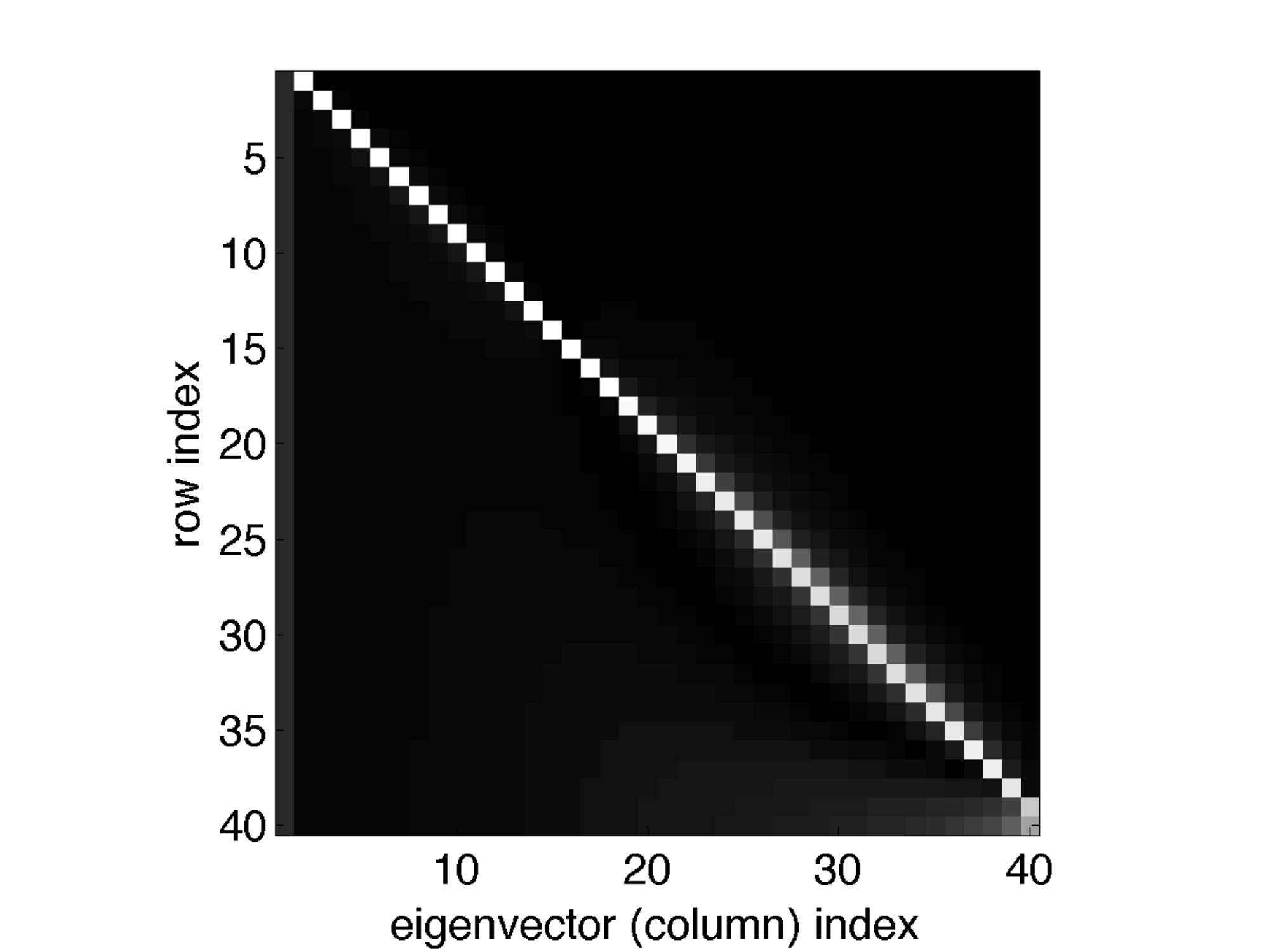} 
\raisebox{0.1in}{\includegraphics[width=0.35 \textwidth]{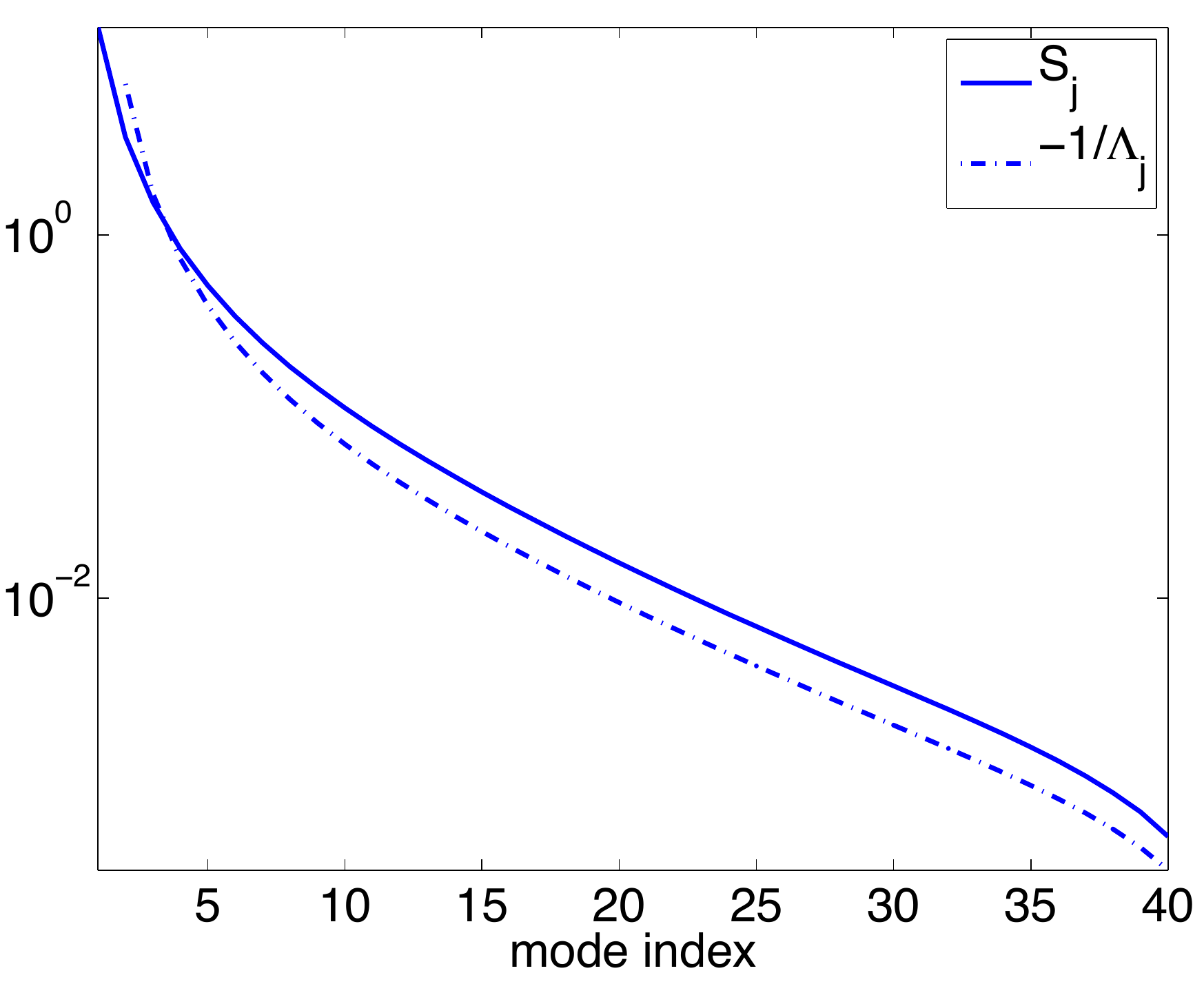}} 

\vspace{-0.03in}
\includegraphics[width=0.43 \textwidth]{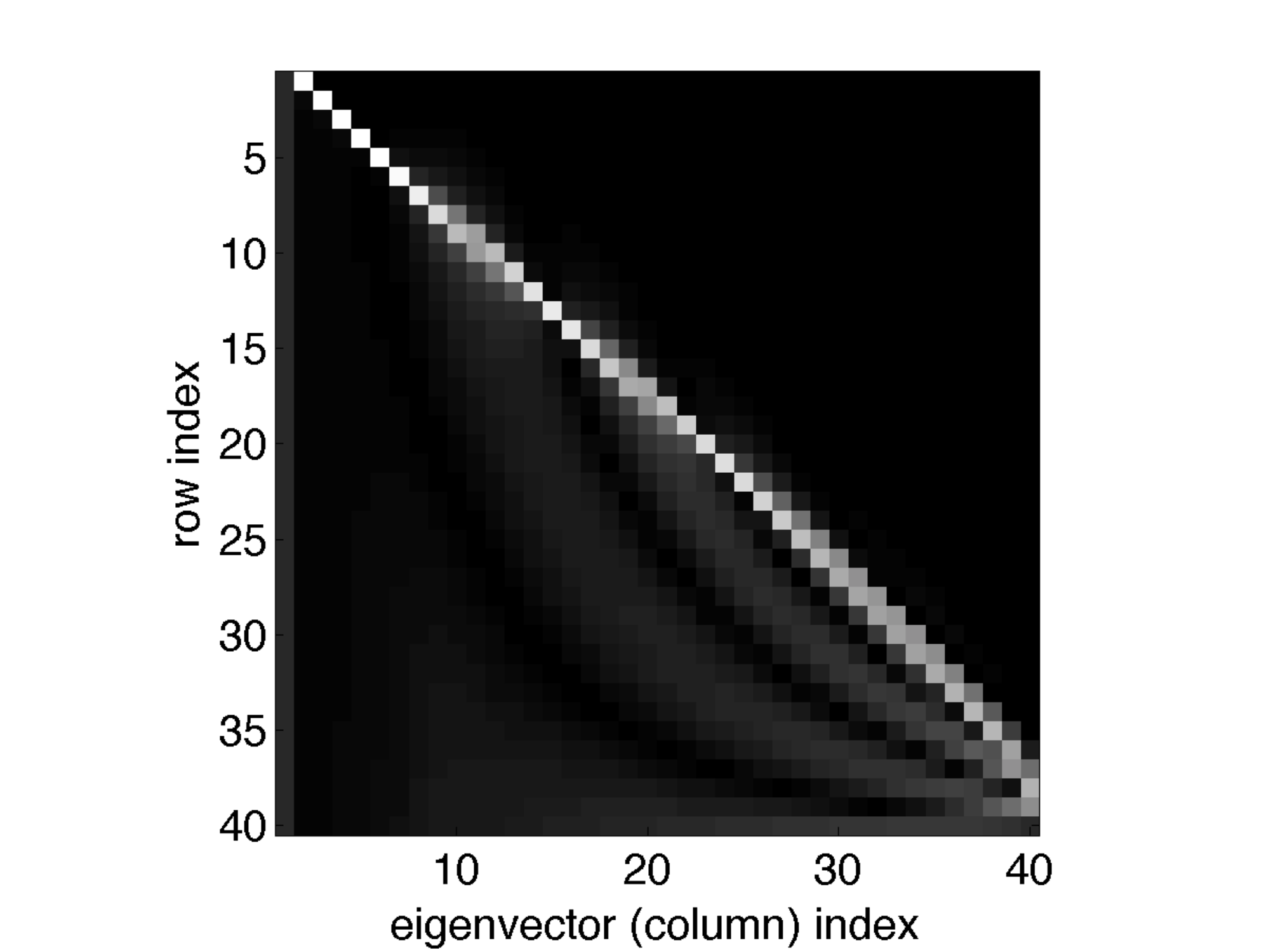} 
\raisebox{0.1in}{\includegraphics[width=0.35 \textwidth]{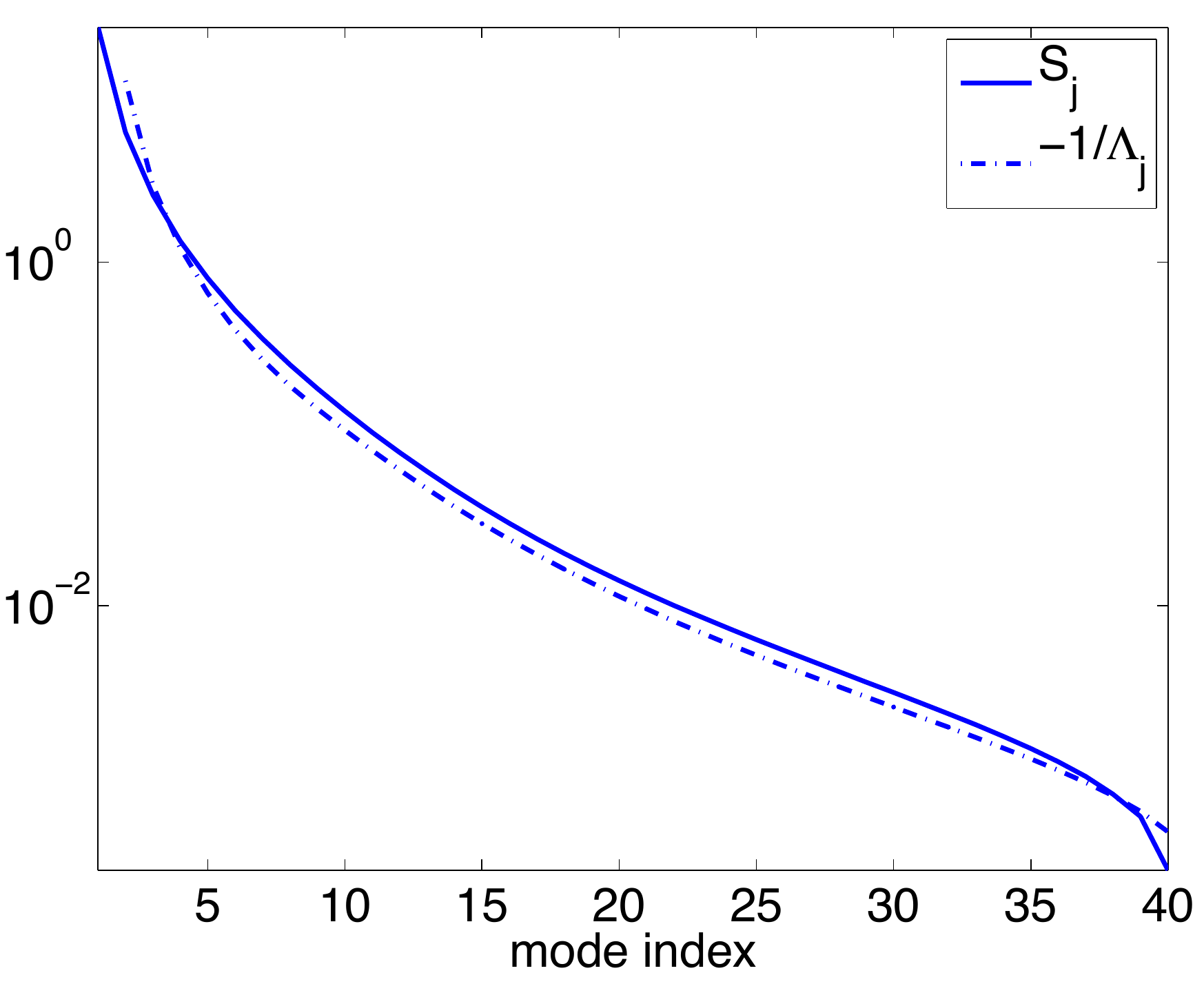}} 

\vspace{-0.03in}
\includegraphics[width=0.43 \textwidth]{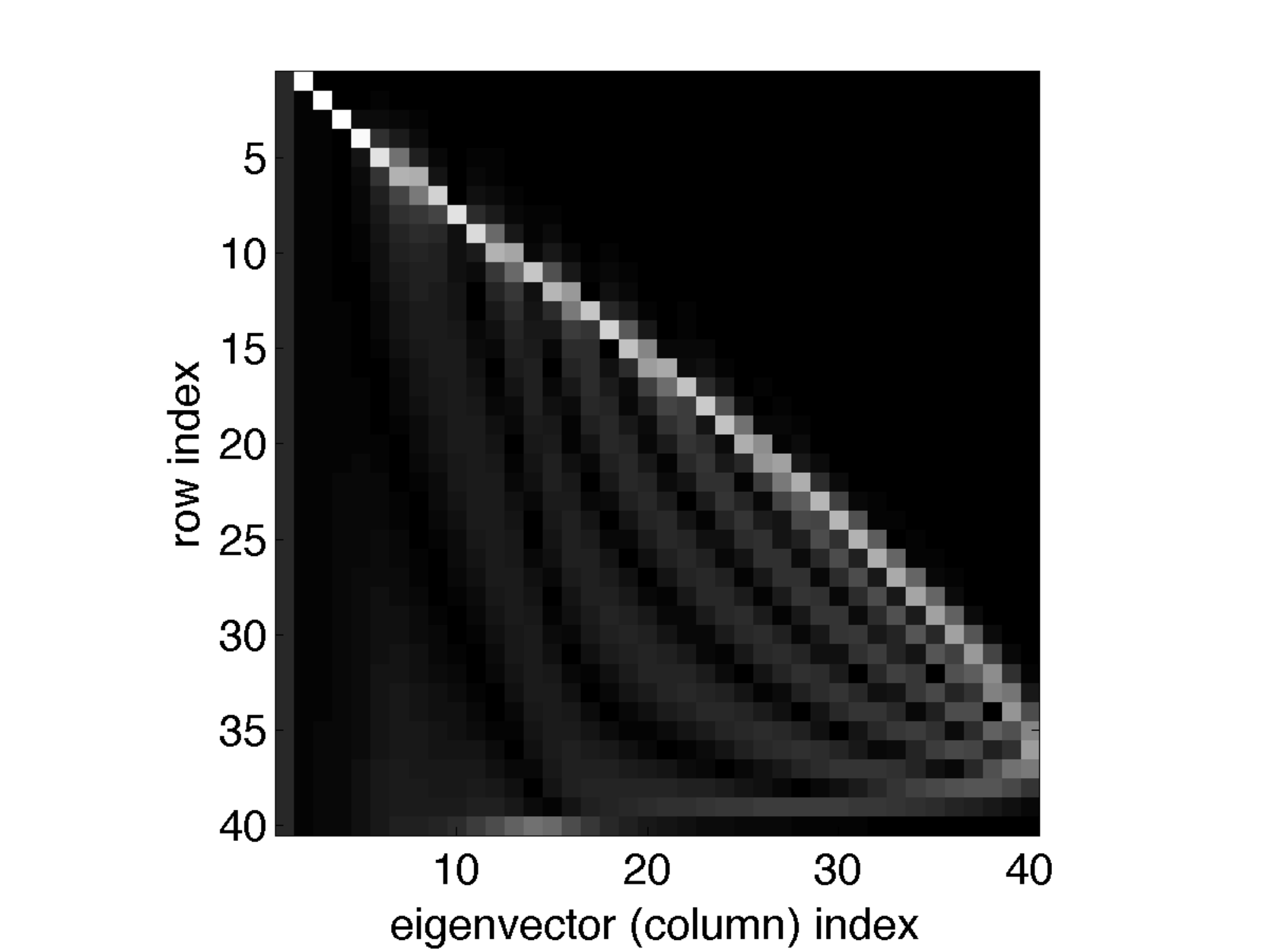} 
\raisebox{0.1in}{\includegraphics[width=0.35 \textwidth]{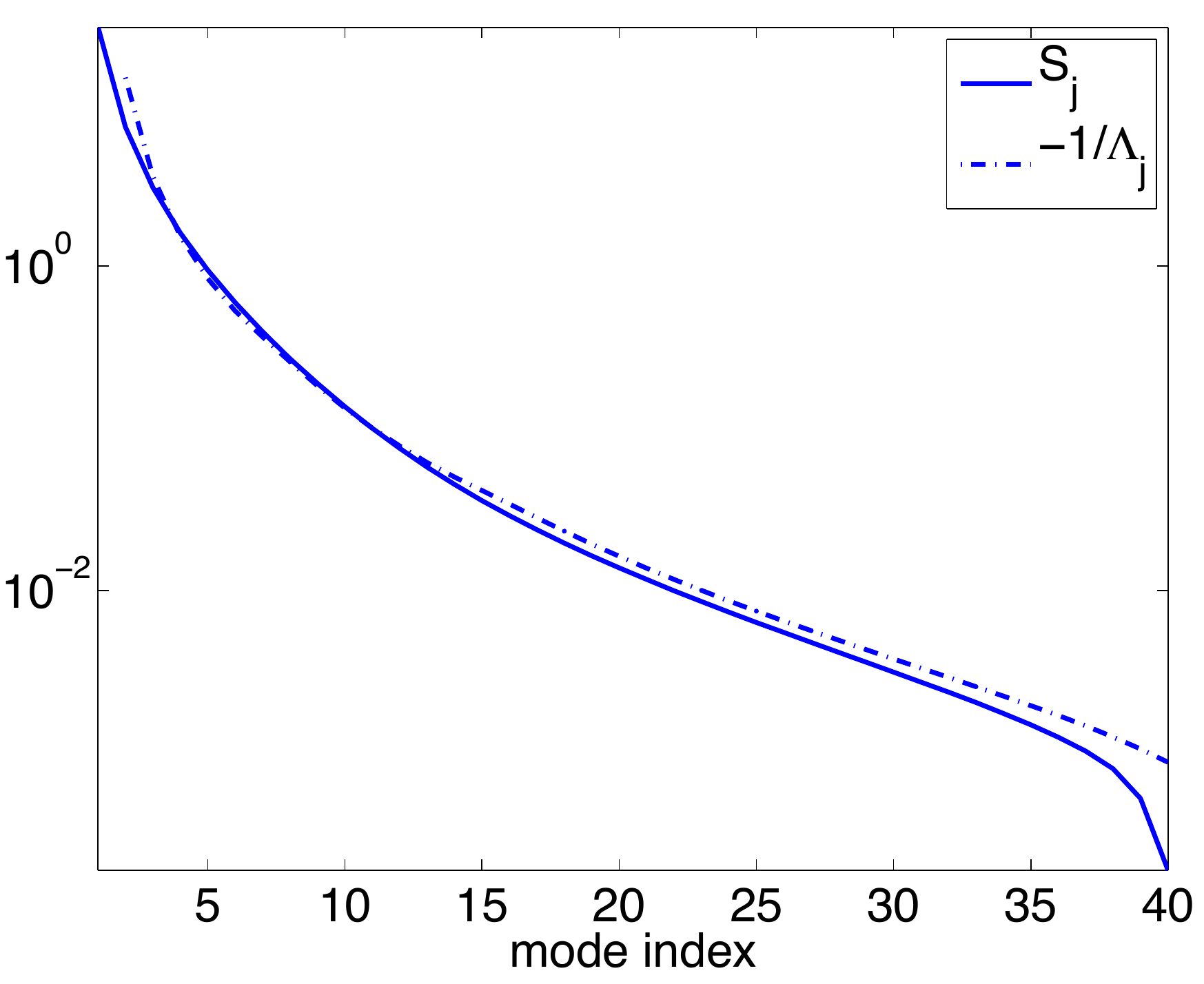}}
\end{center}
\vspace{-0.1in}
\caption{Waveguide with top random boundary. Left plots: Absolute
  values of the entries of the matrix ${\bf U}$ of eigenvectors.  Gray
  scale with lighter color indicates larger values and black indicates
  nearly zero. The column index is in the abscissa and the row index
  in the ordinate. Right plots: The scattering mean free path of the
  modes (full line) and the scales $-1/\Lambda_j$, for $j = 2, \ldots
  N$ (dotted line). In the top row $\ell = \la_o$, in the middle row
  $\ell = 3 \la_o$ and in the last row $\ell = 5 \la_o$. These scales
  should be multiplied by $\ep^{-2}$.}
\label{fig:MatrixUBound}
\end{figure}
\begin{figure}[h]
\begin{center}
\includegraphics[width=0.3 \textwidth]{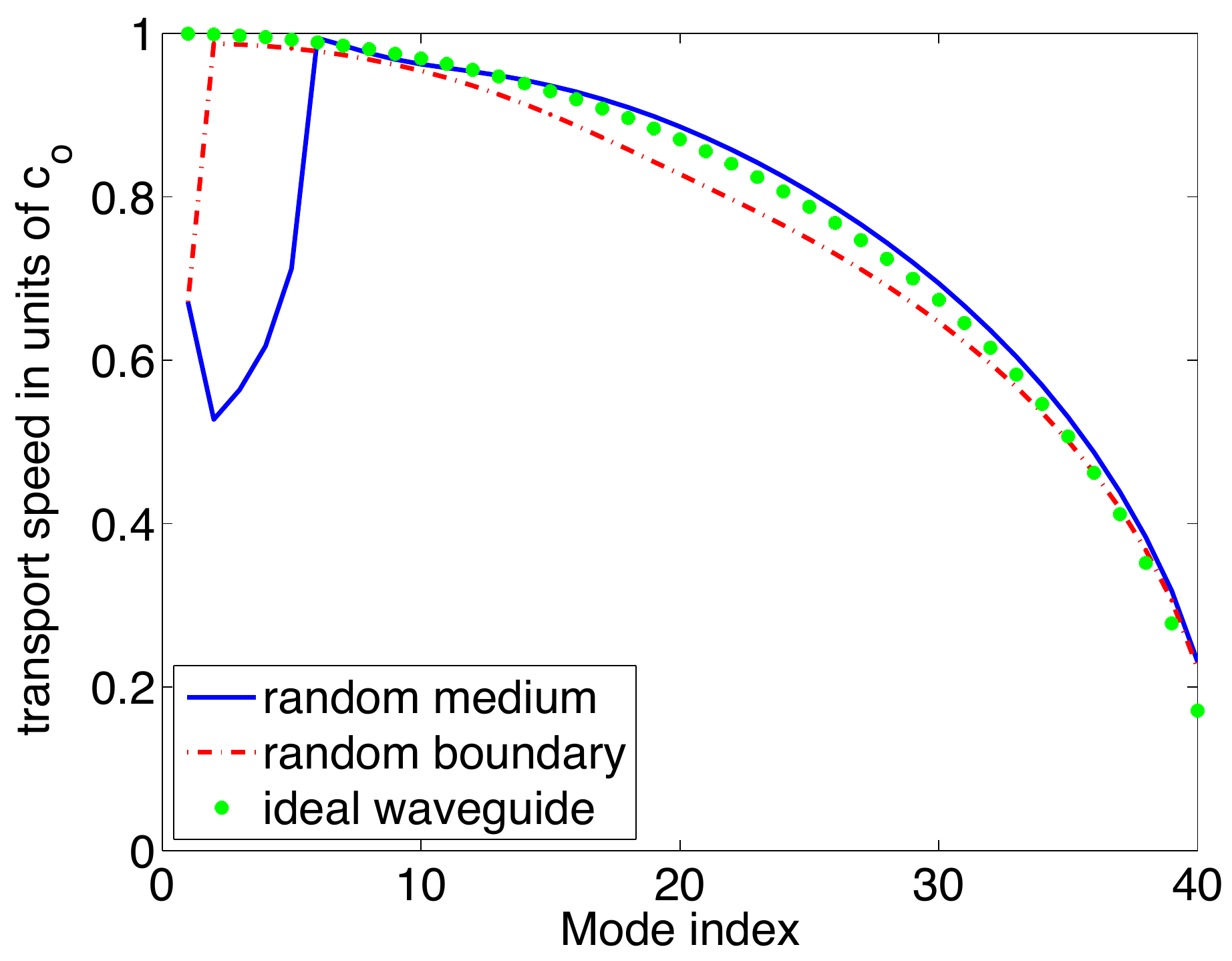} 
\includegraphics[width=0.3 \textwidth]{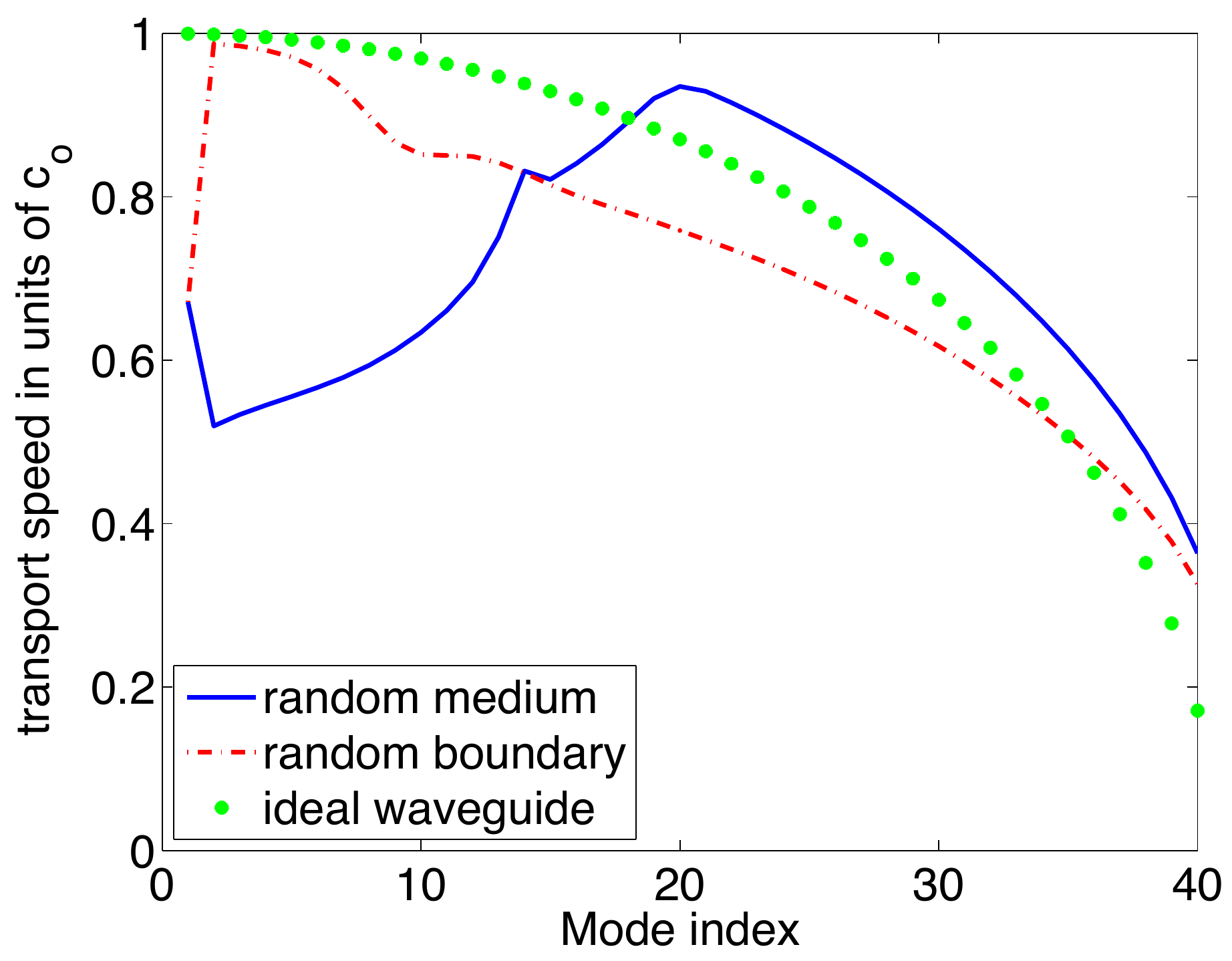}
\includegraphics[width=0.3 \textwidth]{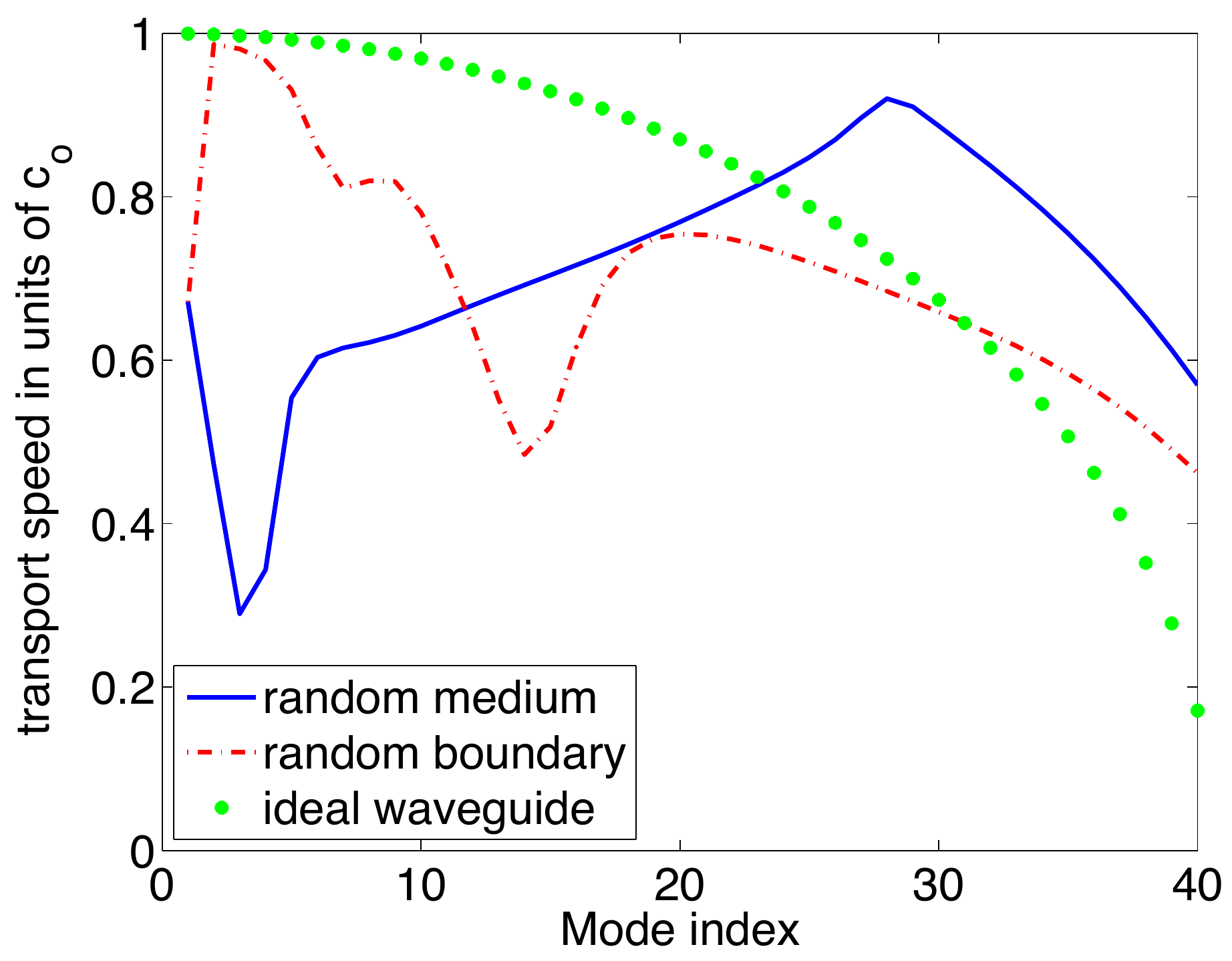}
\end{center}
\vspace{-0.1in}
\caption{Transport speeds in waveguide filled with a random medium
  (solid blue), and with random boundary (dotted red).  The speed in
  ideal waveguides (dotted green). From left to right: $\ell = \la_0$,
  $\ell = 3 \la_o$ and $\ell = 5 \la_o$. The abscissa is mode index
  and the ordinate is the speed scaled by $c_o$.}
\label{fig:TRANSP}
\end{figure}
We display in the left plots of Figures \ref{fig:MatrixU} and
\ref{fig:MatrixUBound} the absolute values of the entries in the
matrix ${\bf U}$ of the eigenvectors, and in the right plots the
scattering mean free paths of the modes and the range scales
$1/|\Lambda_j|$, for $j = 2, \ldots, N$.  We note that the matrix of
eigenvectors has a nearly vanishing block in the upper right
corner. Explicitly, there is an index $j_\star$ such that the first
entries of the eigenvectors ${\bf u}_j$ are negligible for $j >
j_\star$.  In the first simulation in Figure \ref{fig:MatrixU}
$j_\star \approx 5$, in the second $j_\star \approx 15$ and in the
last $j_\star \approx 25$. The effect is more pronounced in the case
of random boundaries where $j^\star \approx 1$ for all three
simulations.

The transport speeds are displayed in Figure \ref{fig:TRANSP}. They
are close to the deterministic ones for most of the modes in the case
$\ell = \la_o$, but they are very different when $\ell = 5
\la_o$. Thus, it is important to use the transport equations in the
range estimation, because the anomalous dispersion induced by
scattering may be significant.

\begin{figure}[h]
\begin{center}
\includegraphics[width=0.4 \textwidth]{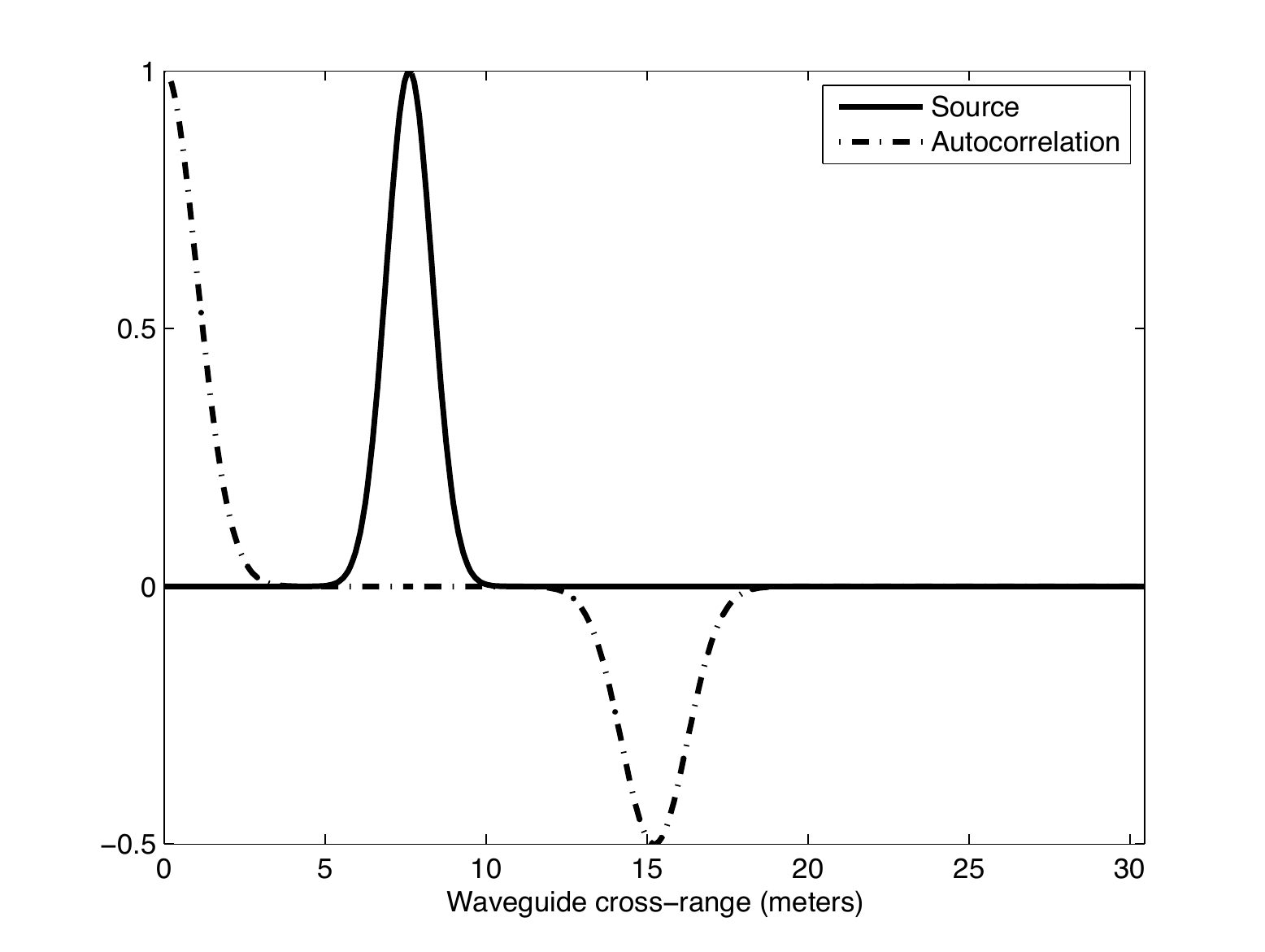}
\includegraphics[width=0.4 \textwidth]{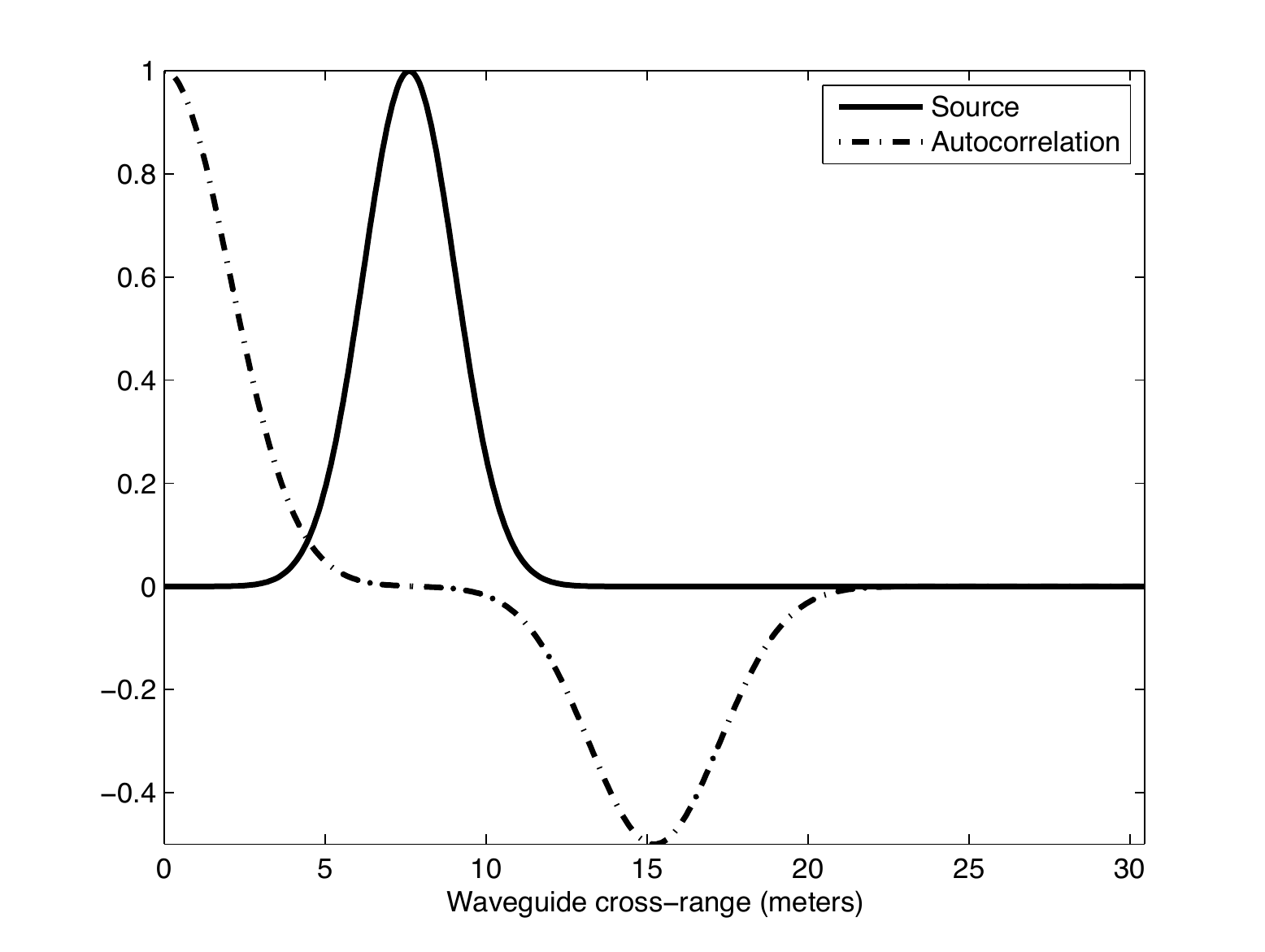}

\vspace{-0.05in} 
\includegraphics[width=0.4 \textwidth]{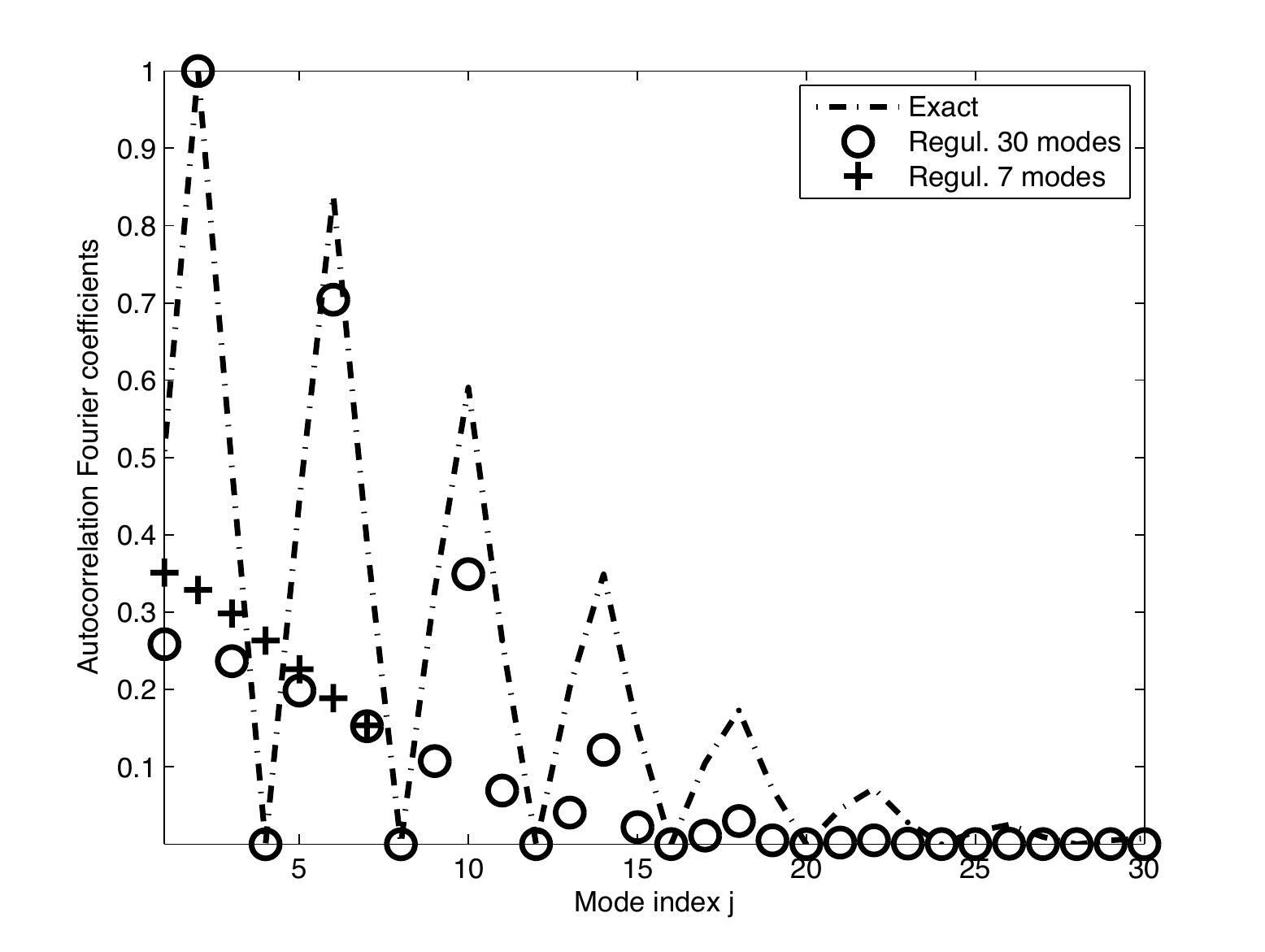} 
\includegraphics[width=0.4 \textwidth]{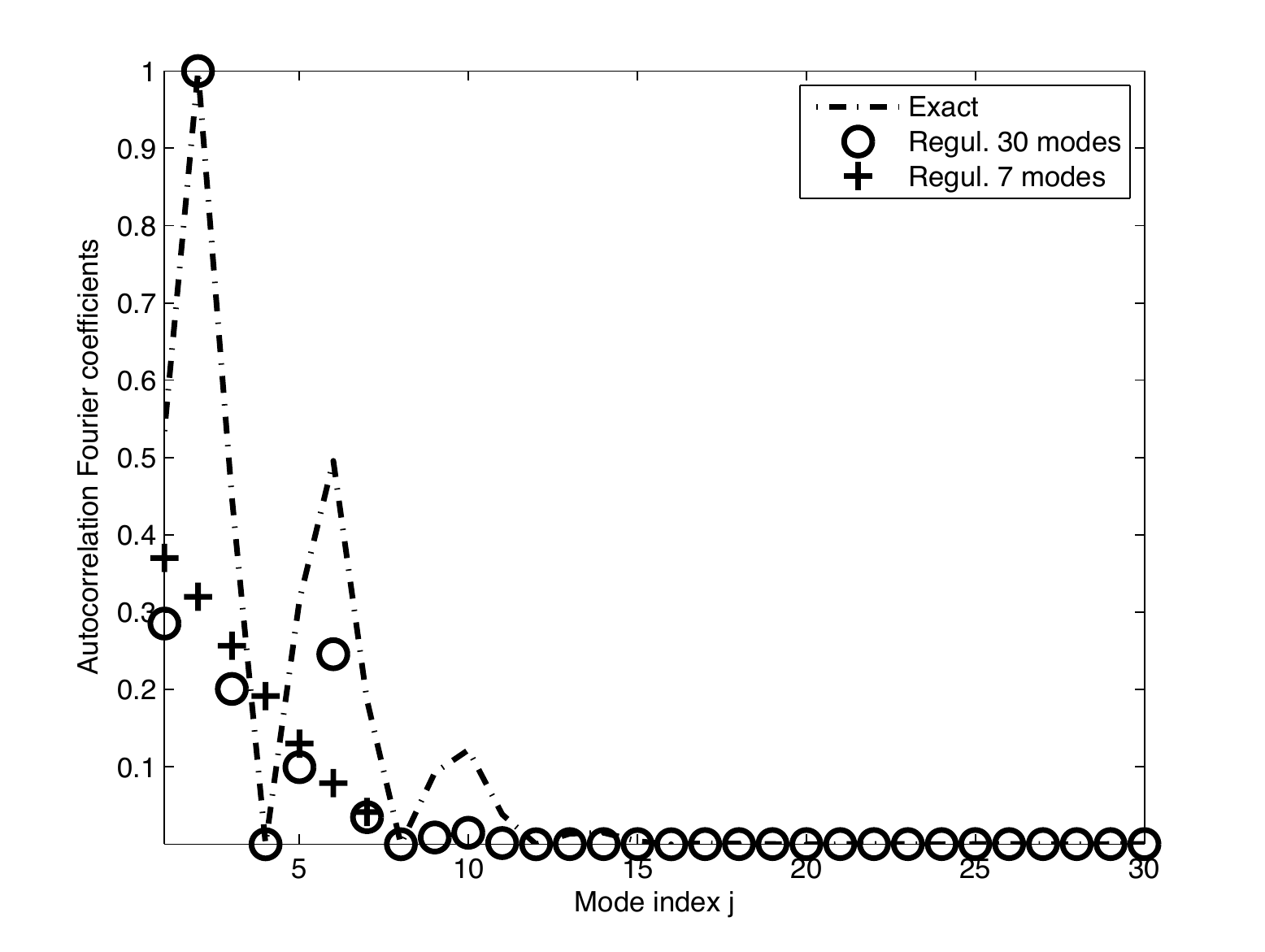} 

\vspace{-0.05in}
\includegraphics[width=0.4 \textwidth]{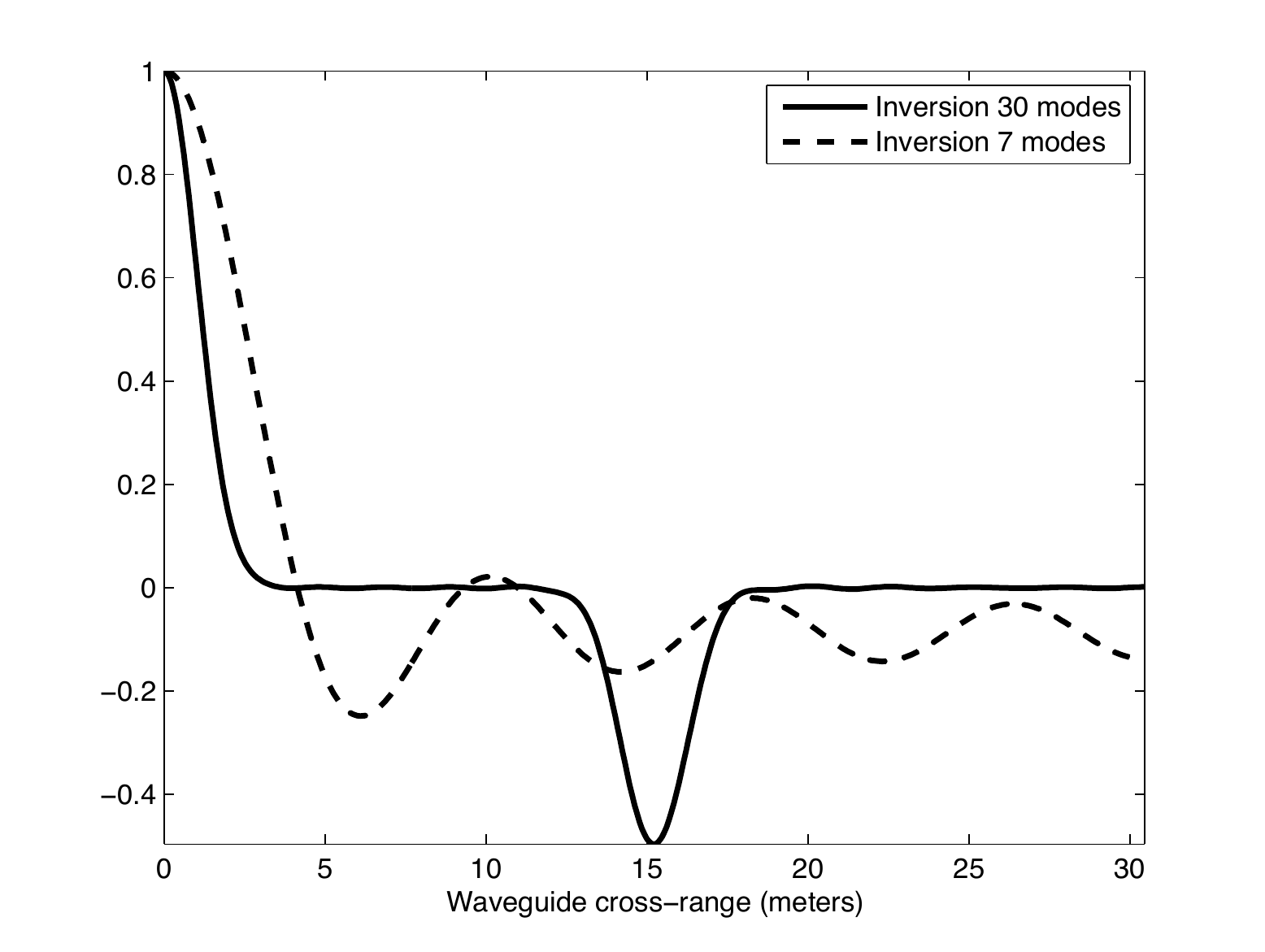}
\includegraphics[width=0.4 \textwidth]{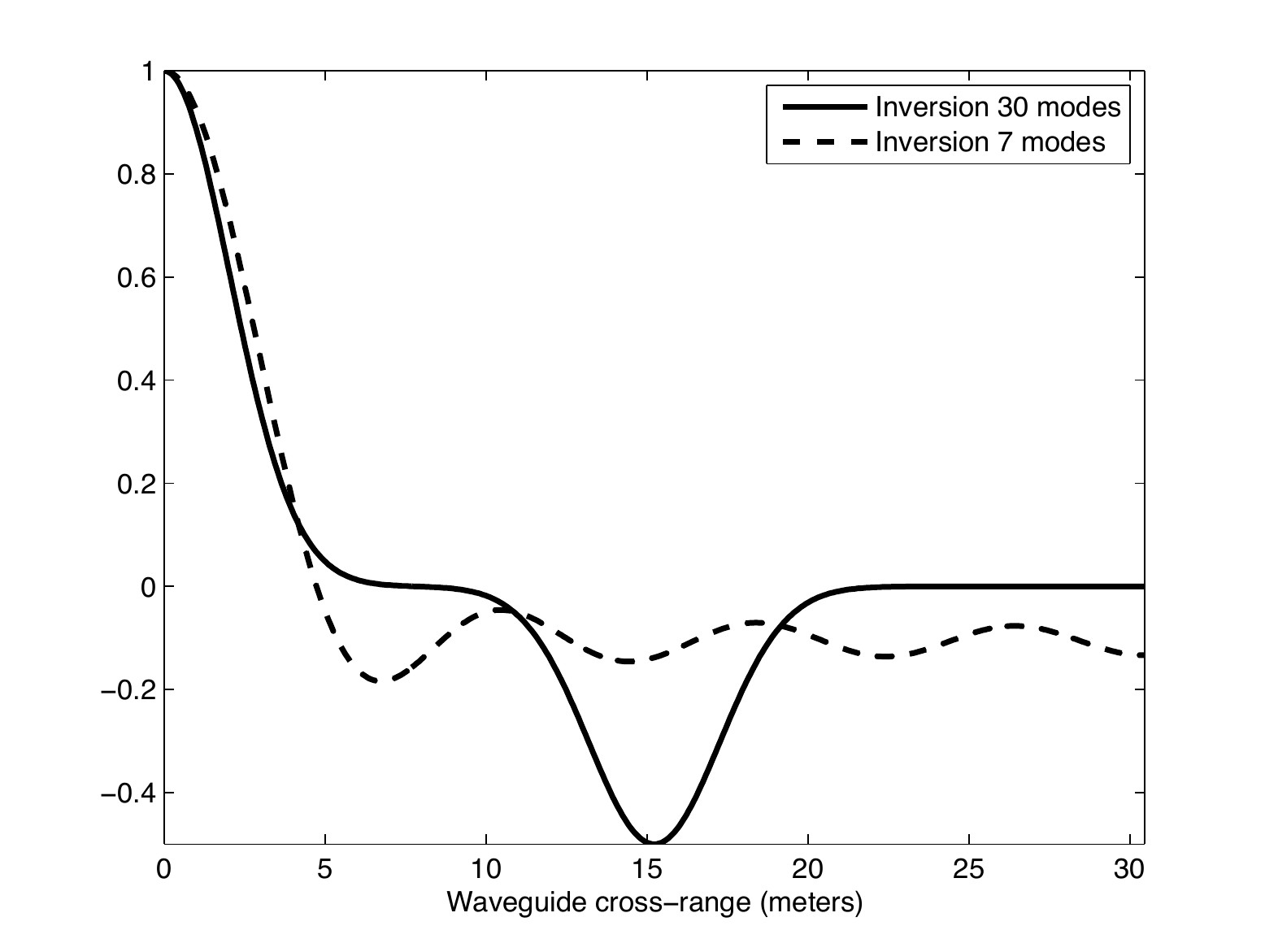}
\end{center}
\vspace{-0.1in}
\caption{Left  $\xi(x) = \mathcal{N}(x_o= X/4, \sigma =
  X/30)$. Right  $\xi(x) = \mathcal{N}(x_o= X/4, \sigma =
  X/15)$. Top line the source cross-range profile $\xi(x)$ (full line)
  and the autocorrelation $\cR_\xi(x)$ (dotted line). Middle plots
  show the exact $|\xi_j|^2$ and the recovered one for cut-off at $J =
  30$ (circle) and $J = 7$ (cross). Bottom plots show the recovered
  $\cR_\xi$ for cut-off at $J = 30$ (full line) and $J = 7$ (dotted
  line). }
\label{fig:inversion}
\end{figure}

In Figure \ref{fig:inversion} we present inversion results for $\xi(x)
= \mathcal{N}(X/4,X/30)$ (left plots) and $\xi(x) =
\mathcal{N}(X/4,X/15)$ (right plots).  The plots on the top line show
the source cross-range profile and the autocorrelation $\cR_\xi(x)$.
The plots in the middle line show the exact values $|\hat \xi_j|^2$
and the estimated ones for cut-off at $J = 30$ and $J = 7$,
respectively. The estimates are calculated using (\ref{eq:inve9}) with
regularization (\ref{eq:XJ}). In the waveguide filled with a random
medium for the cut-off at $J = 30$ the array is at $Z_\cA =
\cL_{eq}/40$, and for $J = 7$ we have $Z_\cA = \cL_{eq}/10.$ The
regularization is chosen so that the exponentials in (\ref{eq:XJ}) are
bounded by $ e^{|\Lambda_j| Z_\cA} \lesssim 10,$ for $j = 1, \ldots,
J.$ The bottom plots show the estimated autocorrelation calculated
using equation (\ref{autocorrelation}), with the series truncated at
$j = J$ and $|\hat \xi_j|^2$ replaced by the estimates.  The results
show that the regularization with $J = 30$ gives a good approximation
of the (first) largest Fourier coefficients and therefore of the
autocorrelation. However, the estimates for $J = 7$ are poor and give
no information about the location of the source (the minimum of the
autocorrelation is not evident in the estimates). The standard
deviation of the Gaussian centered at zero (the peak of the
autocorrelation), which determines the width of the support of the
source, is related to the rate of decay of the Fourier coefficients.
Thus, we can estimate it even for $J = 7$ in the case of the broader
source (bottom right plot) but not for the narrower source (bottom
left plot).

\begin{figure}[h]
\begin{center}
\includegraphics[width=0.31 \textwidth]{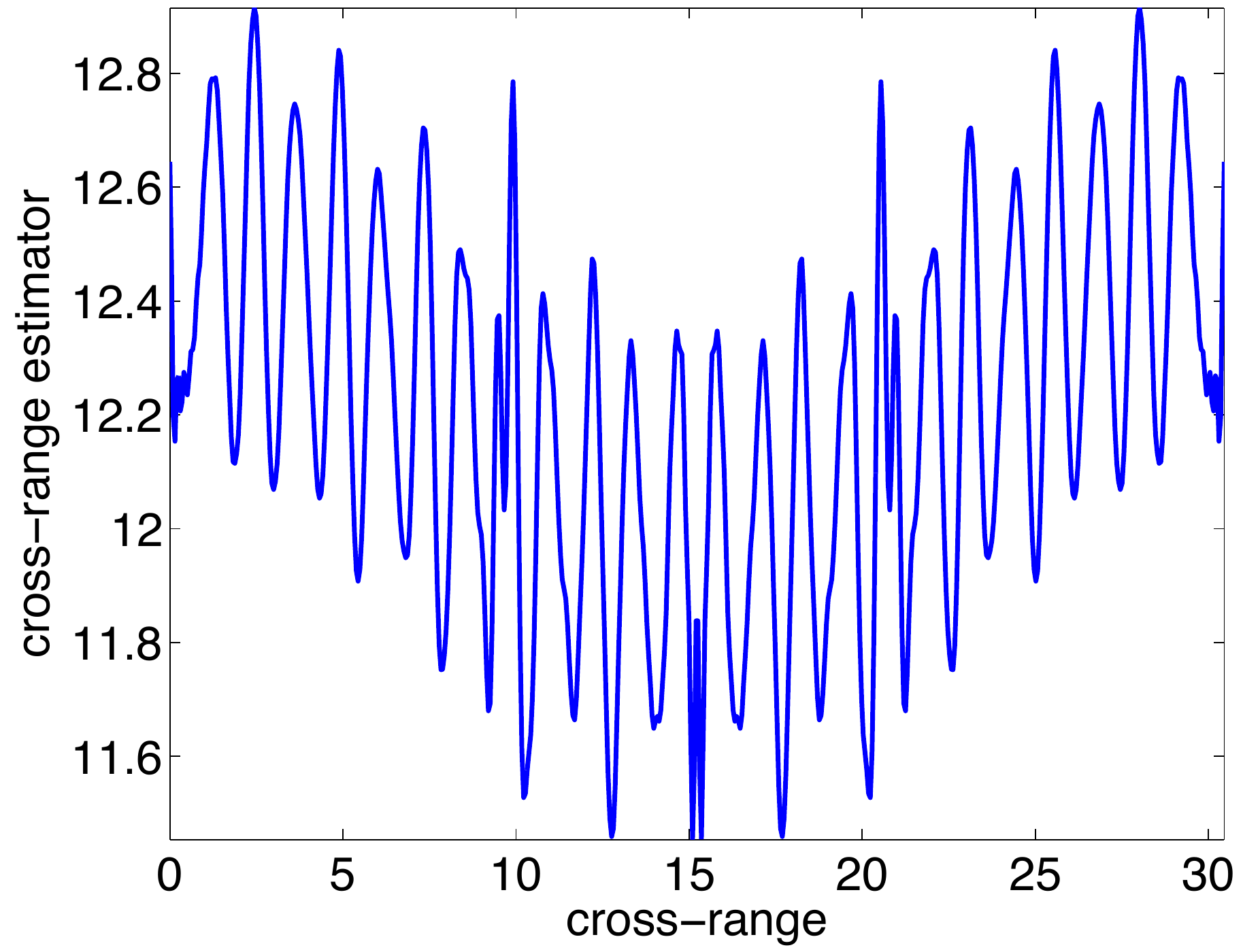}
\includegraphics[width=0.31 \textwidth]{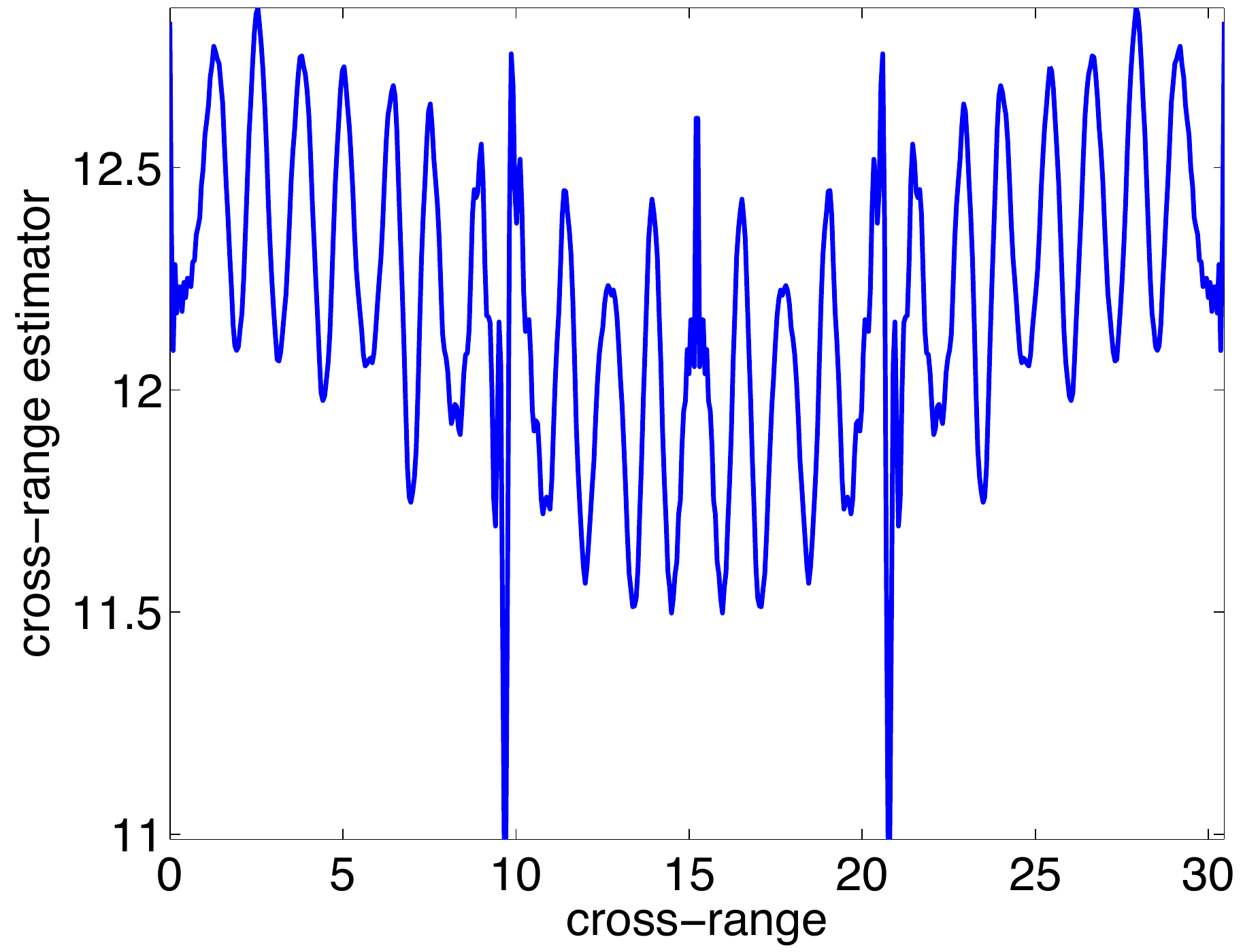}
\includegraphics[width=0.31 \textwidth]{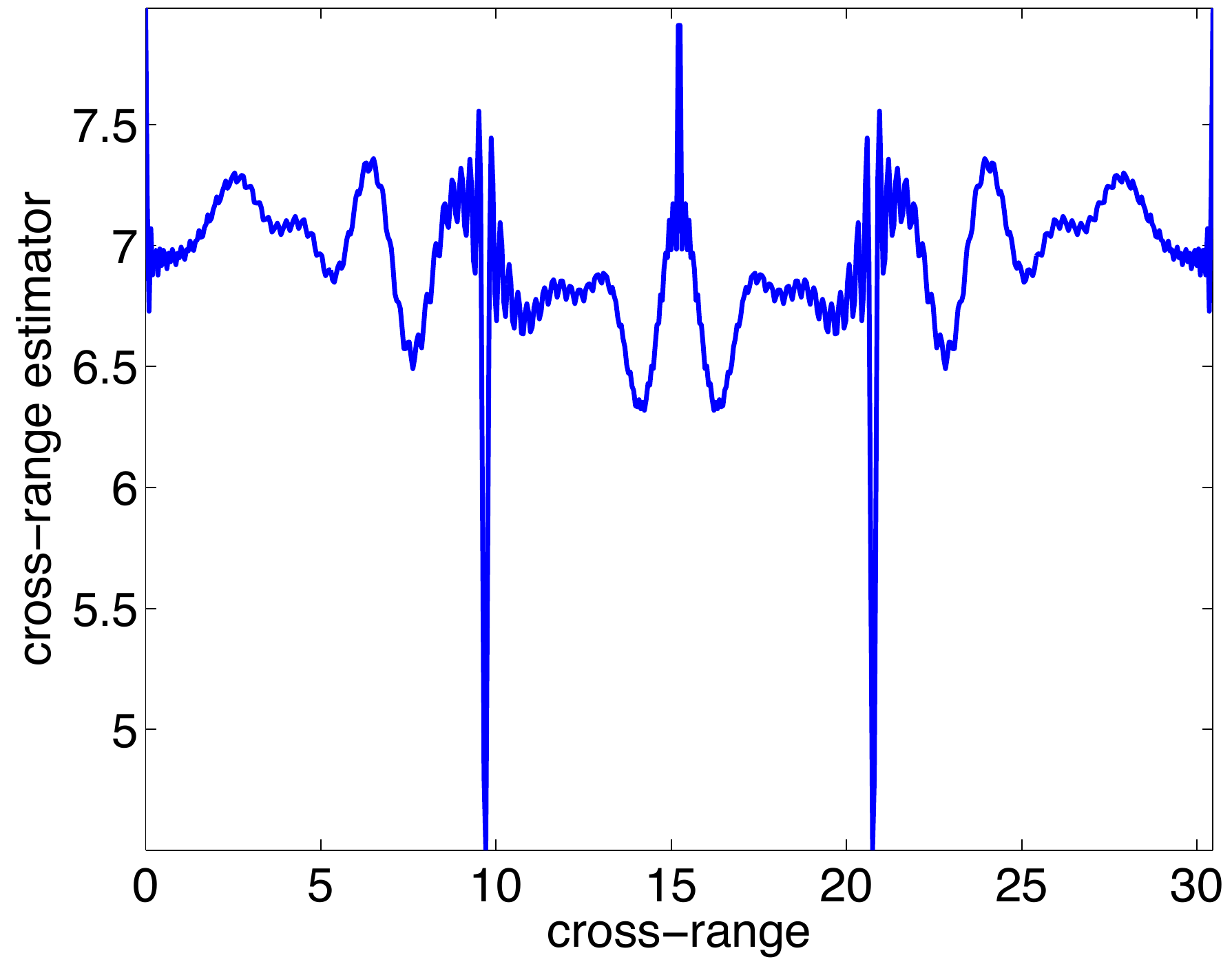}
\end{center}
\vspace{-0.1in}
\caption{ Estimation of the cross-range location of a point-like
  source at $X/\pi \approx 9.697 \lambda_o$ using wideband
  measurements of the waves in the equipartition regime.  Left $\om
  \in (\om_o, 2 \om_o)$, middle $\om \in (\om_o, 3 \om_o)$ and right
  $\om \in (0.5 \om_o,3 \om_o)$. Abscissa is cross-range in $\la_o$
  and ordinate is the objective function (\ref{eq:LSQ}). The minima
  indicate the location of the source and its mirror image with
  respect to the center of the waveguide.}
\label{fig:estim}
\end{figure}
The last numerical illustration in Figure \ref{fig:estim} demonstrates
the benefit of a large bandwidth in the estimation of the cross-range
profile of the source at very long ranges, where the measured waves
are in the equipartition regime, as discussed in section
\ref{sect:I6}. We use the prior knowledge that $\xi(x) \approx
\delta(x-x_o)$, and calculate the objective function
\begin{equation}
\label{eq:LSQ}
\mbox{Obj}(x) = \left\{ \sum_{j=1}^{N_{_M}} \left[\left|\sin
  \left(\frac{\pi j x}{X}\right)\right|^2 - {\boldsymbol{\gamma}}
  \right]^2 \right\}^{1/2},
\end{equation}
where ${\boldsymbol{\gamma}}$ is the solution of the least squares problem 
\begin{equation}
\label{eq:LSQP}
\mbox{arg min} \|\mathbb{B} {\boldsymbol{\gamma}} -
\left(N_1\theta_1, \ldots N_{_M} \theta_{_M})^T \right\|^2 \quad
\mbox{such that}~ {\boldsymbol{\gamma}} \ge {\bf 0},
\end{equation}
where the inequality is understood component-wise. We solve
(\ref{eq:LSQP}) with the MATLAB function \emph{lsqnonneg}. As
indicated in the caption of Figure \ref{fig:estim} we consider three
frequency bands, sampled in steps of $0.02\om_o$. In the first case
$\om \in (\om_o, 2 \om_o)$, so $\mathbb{B} \in \mathbb{R}^{50 \times
  81}$, with rank $50$, and the number of modes ranges from $N_1 = 40$
to $N_{50} = 81$. In the second case $\om \in (\om_o, 3 \om_o)$, so
$\mathbb{B} \in \mathbb{R}^{100 \times 121}$, with rank $92$, and the
number of modes ranges from $N_1 = 40$ to $N_{100} = 121$.  In the
last case $\om \in (0.5\om_o, 3 \om_o)$, so $\mathbb{B} \in
\mathbb{R}^{125 \times 121}$, with rank $110$, and the number of modes
ranges from $N_1 = 20$ to $N_{125} = 121$. Note how in the last two
simulations the minima of the objective function indicate the
cross-range location $x_o = X/\pi$ of the source and its mirror image
with respect to the axis of the waveguide. In the first simulation the
bandwidth is not wide enough and the estimation is ambiguous.

\subsection{Analysis of the structure of the matrix of eigenvectors} 
\label{sect:I4}
Here we use a simplified model $\Upsilon$ of $\Gamma$ to show
with analysis that ${\bf U}$ has a nearly vanishing block in the upper
right corner.  The model neglects the energy transfer between modes
that are not immediate neighbors, meaning that $\Upsilon$ is
tridiagonal. It applies to a different regime than that considered in
the numerical simulations, so the results complement the previous
ones.  The regime is for large correlation lengths satisfying
\[
k_o \ell = O(N), \qquad N \gg 1,
\]
so that $\Upsilon$ is a good approximation (up to a multiplicative
factor) of $\Gamma$.

The definition of $\Upsilon$ is
\begin{equation}
k_o \Upsilon_{jq} = \left\{ \begin{array}{ll} \Gamma_{jq}, \quad
  &|j-q| = 1, \\ 0, &|j-q| >1,
\end{array} \right. 
\label{eq:ana1}
\end{equation}
for $j \ne q$, and
\begin{equation}
k_o \Upsilon_{jj} = \left\{ \begin{array}{ll}
-\Gamma_{jj-1} - \Gamma_{jj+1}, \quad & 
2 \le j \le N-1, \\
-\Gamma_{12}, & j = 1, \\
-\Gamma_{N-1N}, & j = N.
\end{array} \right. 
\label{eq:ana1p}
\end{equation}
We factor out $k_o$ for convenience of the calculations, and use the
expression (\ref{eq:st22}) of $\Gamma$ with the assumption that
the fluctuations in the medium play the dominant role\footnote{That
  the matrix of eigenvectors has negligible entries in the upper right
  corner is pertinent to the estimation of the cross-range profile
  $\xi(x)$. As mentioned in the previous section this is more useful
  in waveguides filled with random media.}. Then, the diagonal of
$\Upsilon$ scales as
\begin{equation}
\Upsilon_{jj} \sim \frac{N^2}{N - j +1}, \qquad j = 1, \ldots, N.
\label{eq:ana2}
\end{equation}
We summarize the properties of the spectrum of $\Upsilon$ in the next
proposition proved in appendix \ref{ap:proof}. We denote its
eigenvectors and eigenvalues with the same symbols ${\bf u}_j$ and
$\Lambda_j$. This is an abuse of notation, but the spectrum of
$\Gamma$ can be related to that of $\Upsilon$ using known
perturbation theory \cite{parlett}.

\vspace{0.1in}
\begin{proposition}
\label{thm1}
The tridiagonal matrix $\Upsilon$ has the following properties:
\begin{enumerate}
\item The eigenvectors form an orthonormal basis of $\mathbb{R}^N$ and
  $\Lambda_j \le 0.$
\item The null space is one dimensional.
\item The norm is $\|\Upsilon\| = O(N^2).$
\item $|\Lambda_j|= O(N^2)$ for indices $j$ satisfying $N-j = O(1)$.
\item If $\Lambda_j$ is a ``large eigenvalue'', meaning that $\delta =
  N/|\Lambda_j| \ll 1$, and $J$ is a spectral cut-off satisfying $J
  \le N/2$, we have $ \displaystyle \sum_{q=1}^J u_{qj}^2 \le O(\delta^2).  $
\end{enumerate}
\end{proposition}
The first two properties are the same as those stated earlier for
$\Gamma$, under the assumption that its off-diagonal entries are
strictly positive. The last property confirms our expectation that the
matrix ${\bf U}$ has a nearly vanishing upper right corner.

\section{Summary}
\label{sect:summary}
We presented an analysis of the inverse source problem in perturbed
two dimensional acoustic waveguides, with data given by time resolved
measurements of the pressure field $p(t,\vx)$ at a remote array of
sensors. The waves are trapped by pressure release boundaries and are
guided along the the range direction, the axis of the waveguide.  The
perturbations consist of small scale fluctuations of the boundaries
and the sound speed in the medium that fills the waveguide. Such
fluctuations cannot be known in detail in practice and are thus
modeled with random processes.  This places the problem in a
stochastic framework. The inversion is carried in a single waveguide,
one realization of the random model, and the goal is to obtain robust
estimates of the source density $\rho(\vx)$. Robust means insensitive
(statistically stable) with respect to the particular realization of
the random perturbations of the waveguide.

Typical imaging methods are based on the assumption that the field
$p(t,\vx)$ is coherent, equal to its statistical expectation plus some
small additive noise. This holds approximately in weak scattering
regimes i.e., when the array is not too far from the source. We
consider strong scattering regimes where $p(t,\vx)$ is incoherent, it
is essentially a random, mean zero field.

Our inversion methodology is based on the theory of wave propagation
in random waveguides
\cite{kohler77,dozier,garnier_papa,garnier2008effective, ABG-12}. This
theory decomposes the wave field in a countable set of modes, which
are time harmonic propagating and evanescent waves. It models the
cumulative wave scattering effects of the perturbations in the
waveguide by the mode amplitudes, which are complex valued random
fields. We use their statistical description to obtain the following
results: (1) We show how to get high fidelity estimates of the energy
carried by the modes to the array from cross-correlations of the
incoherent data. We explain which cross-correlations are useful and
how to calculate them.  (2) As the waves propagate and scatter in the
random waveguide, they interchange energy. This is described by a
system of transport equations with initial condition that depends on
the unknown source density $\rho$. We analyze the invertibility of
this system.  (3) We quantify what can be recovered about the source
in terms of the range to the array. The cumulative scattering effects
impede the inversion process, and the longer the range, the more
pronounced the impediment.

The energies of the propagating modes encode the source information in
terms of a matrix of absolute values of Fourier coefficients of
$\rho$. It is impossible to determine this matrix uniquely from the
estimated energies (the problem is under-determined), unless there is
additional information about $\rho$. We assume that it is a separable
function $\rho(\vx) = \xi(x) \zeta(z)$, where $x$ is the cross-range
component of $\vx$ and $z$ is the range, along the axis of the
waveguide.  We study in detail two cases: (1) The estimation of the
range profile $\zeta(z)$ when the cross-range $\xi(x)$ is known, and
(2) The estimation of the cross-range profile $\xi(x)$ when the source
has point-like support in range $\zeta(z) = \delta(z)$. Other known
range profiles $\zeta(z)$ may be considered as well, but they do not
bring new insight to the inversion process.  In both cases there is
ambiguity about the source, because only the absolute value of the
Fourier coefficients of $\zeta(z)$ or $\xi(x)$ can be determined. We
can expect only limited information about $\rho(\vx)$, such as the
size of its support in range or cross-range. This can be estimated
from the autocorrelation functions of $\zeta(z)$ or $\xi(x)$, which
can be approximated using the absolute values of their Fourier
coefficients.

The range profile estimation turns out to be the easier of the two
cases. We can determine the vector $(|\hat \zeta(\beta_j)|)_{1\le j
  \le N}$ of absolute values of the Fourier transform of $\zeta(z)$
evaluated at the wavenumbers $\beta_j$ of the $N$ propagating modes,
and the calculation is well posed no matter how far the array is from
the source. The wavenumbers $\beta_j$ sample the interval
$(0,\om_o/c_o)$, in steps that decrease monotonically with $N$. Here
$\om_o$ is the central frequency of the signal emitted by the source
and $c_o$ is the reference wave speed in the medium that fills the
waveguide. Thus, we can obtain good approximations of the
autocorrelation of the range profile $\zeta(z)$, specially in high
frequency regimes.

The cross-range estimation entails the calculation of the vector
$(|\hat \xi_j|)_{1 \le j \le N}$ of Fourier coefficients of $\xi(x)$.
The Fourier basis is defined by the eigenfunctions of the second
derivative operator in $x$, which are sin functions in our case.
Although the mode energies define uniquely the vector $(|\hat
\xi_j|)_{1 \le j \le N}$, the calculation is ill posed and the problem
becomes worse as the range separation between the source and the array
increases. Cumulative scattering transfers energy between the modes,
and the longer the waves travel, the harder it is to determine the
initial energy distribution, which is defined by $(|\hat \xi_j|)_{q
  \le j \le N}$. There is a range scale, called the equipartition
distance $\cL_{eq}$, beyond which the energy becomes uniformly
distributed between the modes, independent of the initial state. The
waves lose all information about the cross-range profile at such
ranges, and the inversion for $\xi(x)$ becomes impossible. This is for
a narrow frequency band. If a wide frequency band is available, then
the estimation of the cross-range profile may be improved.

The analysis in this paper is for two dimensional waveguides with
reflecting boundaries. It extends to leaky waveguides where energy is
lost by radiation through a boundary, such as the ocean floor. The
system of transport equations that models the propagation of energy in
such waveguides is derived in \cite[Equation (4.3)]{kohler77}. It is
almost the same as the system analyzed in this paper, expect that
there is damping of energy due to the radiation. This damping adds to
the ill posedness of the inversion.

Extensions to three dimensional acoustic waveguides with reflecting
boundaries are straightforward, and do not introduce anything new if
there are no degeneracies (multiplicity) of the eigenvalues of the
Laplacian in the cross-range.  It is difficult to quantify such
degeneracies for arbitrary cross-sections of the waveguide. But in
certain cases like rectangular cross-sections with sides $L_1$ and
$L_2$, degeneracies occur if and only if $L_1/L_2$ is a rational
number.  In vectorial problems, such as electromagnetic waveguides,
degeneracies are unavoidable for any cross-range profile, because of
different states of polarization of the waves
\cite{AB_EM,marcuse91}. Degeneracies are interesting because they
introduce statistical correlations between the amplitudes of the modes
that correspond to degenerate eigenvalues. We no longer have scalar
valued energies carried by each mode, but Hermitian matrices that
describe the propagation of energy by the set of degenerate modes
\cite{AB_EM}. The transport equations are more complicated
\cite{AB_EM}, but they may lead to extra information about the
cross-range profile of the source, not just the absolute value of its
Fourier coefficients. However, there is no gain in the stability of
the inverse problem. The transfer of energy between the modes occurs
in any type of random waveguide, and the estimation of the initial
energy state, which determines the cross-range of the source, remains
exponentially ill-posed for narrow bandwidths.

\section*{Acknowledgements}
This work was partially supported
by the AFOSR Grant FA9550-12-1-0117 and the ONR Grant N00014-12-1-0256.

\appendix 
\section{The model of the cross-correlations}
\label{ap:A}
We obtain from (\ref{eq:inv2}), (\ref{eq:inv4}) and definition
(\ref{eq:st13}) that relates the mode amplitudes to the propagator
that
\begin{align}
\hat \cC_j (h) \approx & \, \hat
\psi\left(\frac{h}{H}\right)\int_{-\infty}^\infty \frac{d \om}{B^2}
\left| \hat f \left( \frac{\om-\om_o}{B} \right) \right|^2 \sum_{q, q'
  = 1}^N Q_{jq} Q_{jq'} \sum_{l,l'=1}^N \frac{1}{4
  \beta_l(\om_o)\beta_{l'}(\om_o)} \times \nonumber
\\ &\int_{-\infty}^\infty \frac{d u}{2 \pi } \, \hat \chi (u) e^{i u
  {\left[\ep^2 t_o-\beta_q'(\om_o) Z_{_\cA}\right]}/{T}}
\int_{-\infty}^\infty \frac{d u'}{2 \pi } \, \overline{\hat \chi (u')}
e^{-i u' \left[\ep^2 t_o-\beta_{q'}'(\om_o) Z_{_\cA}\right]/{T}}
\times \nonumber \\ & \qquad e^{ i \frac{Z_{_\cA}}{\ep^2} \left\{
  \beta_q(\om_o) - \beta_{q'}(\om_o) + (\om-\om_o) [\beta_{q}' (\om_o)
    -\beta'_{q'}(\om_o)]\right\} + i h \beta_{q'}'(\om_o) Z_{_\cA}}
\times \nonumber \\& \int_0^X \hspace{-0.02in} dx \hspace{-0.02in}
\int_{-\infty}^\infty \hspace{-0.02in} dz \, \rho(x,z) \phi_l(x) e^{-i
  \beta_q(\om)z} \int_0^X \hspace{-0.02in} dx' \hspace{-0.02in}
\int_{-\infty}^\infty \hspace{-0.02in} dz' \, \overline{\rho(x',z')}
\phi_{l'}(x') e^{i \beta_{q'}(\om)z'} \times \nonumber \\ & \qquad
\mathbb{P}_{ql}^\ep \left(\om-\frac{\ep^2 u}{T},Z_{_\cA},z'\right)
\overline{ \mathbb{P}_{q'l'}^\ep \left(\om-\ep^2 h - \frac{\ep^2
    u'}{T},Z_{_\cA},z'\right)}. \label{eq:A1}
\end{align}
When we calculate the expectation of (\ref{eq:A1}) using the moment
formula (\ref{eq:st18}), we see that only the terms with $q = q'$ and
$l = l'$ survive in the sum. The coherent terms for $q = l$, 
$q'=l'$ and $q \ne q'$ in the second moment (\ref{eq:st18}) do not
appear at full aperture, where $Q_{jq} = \delta_{jq}$. It the array
has partial aperture but $\cA$ is large enough to have a diagonally
dominant matrix $Q$, the coherent terms are small because of the small
weights $Q_{jq}$ for $q \ne j$, and specially because of the assumption
that $Z_{_\cA} >\cS_1$. The result is 
\begin{align}
\EE \left[\hat \cC_j(h)\right] \approx \, \hat \psi
  \left(\frac{h}{H}\right) \int_{-\infty}^\infty \frac{d \om}{B^2}
  \left| \hat f \left( \frac{\om-\om_o}{B} \right) \right|^2
  \sum_{q=1}^N Q^2_{jq} \sum_{l=1}^N \frac{\left|\hat
    \rho_l\left[\beta_q(\om)\right]\right|^2 }{4 \beta_l(\om_o)
    \beta_q(\om_o)} \times \nonumber \nonumber \\ 
  ~\int_{-\infty}^\infty \frac{d u}{2 \pi } \, \hat \chi (u)
  \int_{-\infty}^\infty \frac{d u'}{2 \pi } \, \overline{\hat \chi
    (u')} e^{i (u-u') \ep^2 t_o/T} \hat \cW_q^{(l)}\left(\om, h +
  \frac{u'-u}{T},Z_{_\cA}\right) ,
\label{eq:A3}
\end{align}
where $\hat \rho_k(\beta)$ are the Fourier coefficients of the source
density defined by (\ref{eq:A2}). Taking the inverse Fourier transform
of (\ref{eq:A3}) and changing variables $ h' = h +
{(u'-u)}/{T} $ 
\begin{align}
\EE \left[\cC_j(\tau)\right] \approx &\, \left| \hat
\chi\left(\frac{\tau-\ep^2
  t_o}{T}\right)\right|^2\int_{-\infty}^\infty \frac{d \om}{B^2}
\left| \hat f \left( \frac{\om-\om_o}{B} \right) \right|^2
\sum_{q=1}^N Q^2_{jq} \times \nonumber \\ &\sum_{l=1}^N
\frac{\left|\hat \rho_l\left[\beta_q(\om)\right]\right|^2 }{4
  \beta_l(\om_o) \beta_q(\om_o)} \int_{-\infty}^\infty\frac{d h'}{2
  \pi} \hat \cW_q^{(l)}\left(\om, h',Z_{_\cA}\right) \hat \psi
\left(\frac{h'}{H}\right) e^{-i h' \tau}.
\label{eq:A3N}
\end{align}
Here we used that $\psi$ is smooth and $HT \gg 1$, to approximate
\[
\hat \psi \left(\frac{h'}{H}- \frac{u'-u}{TH}\right) \approx \hat \psi
\left(\frac{h'}{H}\right).
\]
Equation (\ref{eq:A4}) follows from (\ref{eq:A3N}) and definition
(\ref{eq:st20FT}) of the Wigner transform. We also use that the
bandwidth is small and that the Fourier transform $\hat f$ of the
pulse and $\hat \rho_l(\beta)$ are smooth in $\om$. In fact the latter
is analytic because $\rho$ has compact support.

To assess the statistical stability of $\cC_j(\tau)$ we need the
fourth order moments of the propagator. These are given in
\cite[Appendix D]{BIT-10}, and the variance of $\cC_j(\tau)$ follows
after a long calculation which we explain briefly. Since it is defined
by
\[
{\rm var}\left[\cC_j(\tau)\right] = \EE
\left[\left|\cC_j(\tau)\right|^2\right] - \left|\EE
\left[\cC_j(\tau)\right]\right|^2,
\]
we need fourth order moments like
\[
\EE\left[\mathbb{P}_{q_1l_1}^\ep
  \left(\om,Z_{_\cA},z_1\right) \overline{ \mathbb{P}_{q_2l_2}^\ep
    \left(\om,Z_{_\cA},z_2\right)}\mathbb{P}_{q'_1l'_1}^\ep
  \left(\om',Z_{_\cA},z'_1\right) \overline{ \mathbb{P}_{q'_2l'_2}^\ep
    \left(\om',Z_{_\cA},z'_2\right)}\right]
\] 
where we neglect the order $\ep^2$ offsets in the arguments, because
they do not play any role. These moments factorize in the product of
two second moments at frequencies $\om$ and $\om'$ when $|\om'-\om|
\gg \ep^2 \om_o$, so in the calculation of the variance we are left
with the integration over the small strip $\{\om, \om' : ~ ~
|\om-\om'| \ll \ep^2 \om_o\}$. This makes the variance smaller than
the square of the mean (\ref{eq:A4}), by a factor of $\ep^2 \om_o/B =
\ep^{2-\alpha} \ll 1$, as long as the mean is large. This happens for
example when the matrix $Q$ is diagonally dominant and we evaluate the
cross-correlation at a time $\tau$ for which $\cW_j^{(l)}$ is large.

\section{Justification of the perturbative analysis of arrival times}
\label{ap:B}
To show that $i h \cB' + \Gamma$ is a perturbation of $\Gamma$ for
$|h| \le H$, let us calculate the ratio $|\Gamma_{jj}|/(H \beta'_j)$
for $j = 1, \ldots, N$.

Using (\ref{eq:st22}) and the definition of $\beta_j$ we have
\begin{align}
\frac{|\Gamma_{jj}|}{H \beta'_j} = \frac{\om_o}{H} \frac{\pi^4
  \ell j^2}{k_o^2 X^4} \sum_{q \ne j} \frac{q^2}{\beta_q}\left\{\hat
\cR_{_B} \left[\ell (\beta_j-\beta_q)\right] + \hat \cR_{_T}\left[\ell
  (\beta_j-\beta_q)\right] \right\}+ \nonumber \\ \frac{\om_o}{H}
\frac{k_o^2 \ell }{4} \sum_{q \ne j} \frac{1}{\beta_q}\hat
\cR_{\nu_{jq}}\left[\ell(\beta_j-\beta_q)\right],
\label{eq:B1}
\end{align}
and we estimate next each term.  For the first term
\begin{equation}
T_1 = \frac{\om_o}{H} \frac{\pi^4 \ell j^2}{k_o^2 X^4} \sum_{q \ne j}
\frac{q^2}{\beta_q}\hat \cR_{_B} \left[\ell (\beta_j-\beta_q)\right]
\label{eq:B2}
\end{equation}
we use that $ \pi/X \approx k/N$ and $ \beta_q \approx k
\sqrt{1-(q/N)^2}, $ to write
\begin{align*}
T_1 &\approx \frac{\om_o}{H} \frac{k_o \ell j^2}{N} \frac{1}{N} \sum_{q
  \ne j} \frac{(q/N)^2}{\sqrt{1-(q/N)^2}} \hat \cR_{_B} \left[k_o \ell
  (\sqrt{1-(j/N)^2}-\sqrt{1-(q/l)^2})\right] \nonumber \\ & \approx
\frac{\om_o}{H} \frac{k_o \ell j^2}{N} \int_0^1 du \,
\frac{u^2}{\sqrt{1-u^2}} \hat \cR_{_B} \left[k_o \ell
  (\sqrt{1-(j/N)^2}-\sqrt{1-u^2})\right].
\end{align*}
Moreover, changing variables $ s = \sqrt{1-(j/N)^2}-\sqrt{1-u^2}, $ we
obtain
\begin{align}
T_1 \approx \frac{\om_o}{H} \frac{k_o \ell j^2}{N} \int ds
\left[1-(s-\sqrt{1-(j/N)^2})^2\right]^{1/2} \hat \cR_{_B} (k_o \ell s)
= O\left(\frac{\om_o}{H} \frac{j^3}{N^2}\right).
\end{align}
The second term in (\ref{eq:B1}) is similar to $T_1$ and to estimate
the third term we need
\begin{align*}
\hat \cR_{\nu_{jq}}\left[\ell(\beta_j-\beta_l)\right] = \frac{4}{X^4}
\int d \eta \, e^{i \ell (\beta_j-\beta_q) \eta} \int_0^X d x_1
\int_0^X dx_2 \, \sin \left(\frac{\pi j x_1}{X}\right) \times \\ \sin
\left(\frac{\pi j x_2}{X}\right)\sin \left(\frac{\pi q x_1}{X}\right)
\sin \left(\frac{\pi q x_2}{X}\right)
\cR_\nu\left(\frac{x_2-x_1}{\ell},\eta\right),
\end{align*}
where we assume that the fluctuations $\nu$ are stationary in both
range and cross-range. Changing variables to $ \bar x = (x_1 + x_2)/2$
and $\tilde x = x_2-x_1,$ we have using basic trigonometry
\begin{align*}
\hat \cR_{\nu_{jq}}\left[\ell(\beta_j-\beta_l)\right] \approx
\frac{\ell}{X} \int d \eta \, e^{i \ell (\beta_j-\beta_q) \eta}
\int_{-X/\ell}^{X/\ell} \frac{d \tilde x}{\ell} \cos \left(\frac{\pi j
  \tilde x}{X}\right)\cos \left(\frac{\pi q \tilde x}{X}\right)
\cR_\nu\left(\frac{\tilde x}{\ell},\eta\right).
\end{align*}
Moreover, assuming that $X \gg \ell$ we can approximate the last
integral by the Fourier transform of $\cR_\nu$ in the first argument,
denoted by $\breve \cR_\nu$, and get
\begin{align*}
\hat \cR_{\nu_{jq}}\left[\ell(\beta_j-\beta_l)\right] &\approx
\frac{\ell}{2 X} \int d \eta \, e^{i \ell (\beta_j-\beta_q) \eta}
\left[ \breve \cR_\nu \left(\frac{k \ell (j-q)}{N},\eta\right) +\breve
  \cR_\nu\left(\frac{k \ell (j-q)}{N},\eta\right)\right] \nonumber
\\ &= \frac{\ell}{2 X} \left[ \hat \cR_\nu \left(\frac{k \ell
    (j-q)}{N},\ell (\beta_j-\beta_q)\right) +\hat
  \cR_\nu\left(\frac{k \ell
    (j-q)}{N},\ell(\beta_j-\beta_q)\right)\right].
\end{align*}
The third term in (\ref{eq:B1}) becomes 
\begin{align}
T_3 &= \frac{\om_o}{H} \frac{k_o^2 \ell }{4} \sum_{q \ne j}
\frac{1}{\beta_q}\hat \cR_{\nu_{jq}}\left[\ell(\beta_j-\beta_q)\right]
\nonumber \\ &= \frac{\om_o}{H} \frac{(k_o \ell)^2}{8 \pi} \frac{1}{N}
\sum_{q \ne j} \left[ \hat \cR_\nu \left(\frac{k \ell (j-q)}{N},\ell
  (\beta_j-\beta_q)\right) +\hat \cR_\nu\left(\frac{k \ell
    (j-q)}{N},\ell(\beta_j-\beta_q)\right)\right] \nonumber
\\ &\approx \frac{\om_o}{H} \frac{(k_o \ell)^2}{8 \pi} \int_0^1
\frac{du}{\sqrt{1-u^2}} \hat \cR_\nu \left(k \ell (j/N-u),k_o \ell
(\sqrt{1-(j/N)^2}-\sqrt{1-u^2})\right) \nonumber \\ &= O
\left(\frac{\om_o}{H} k_o \ell\right), \label{eq:B4}
\end{align}
where we used that $k_o \ell \gg 1$ in forward scattering
approximation regimes. Indeed, the forward scattering approximation 
requires that \cite{kohler77,garnier_papa,ABG-12} 
\[
\hat \cR_{\nu_jq}[\ell(\beta_j + \beta_q)] \ll 1, \qquad \forall ~j, q
= 1, \ldots, N,
\]
which implies $k_o \ell \gg 1$.

\begin{figure}[t]
\begin{center}
\includegraphics[width=0.335 \textwidth]{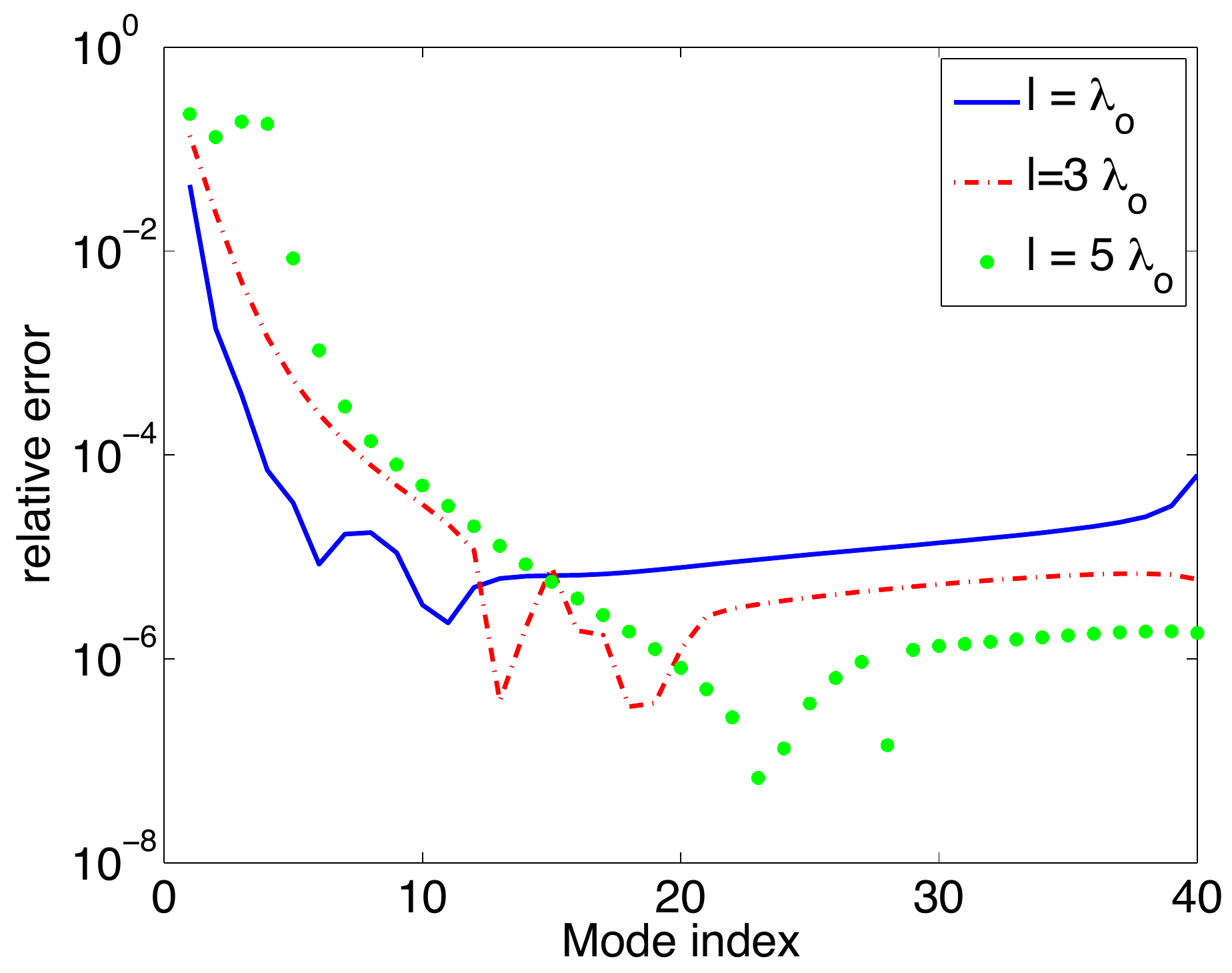}
\includegraphics[width=0.32 \textwidth]{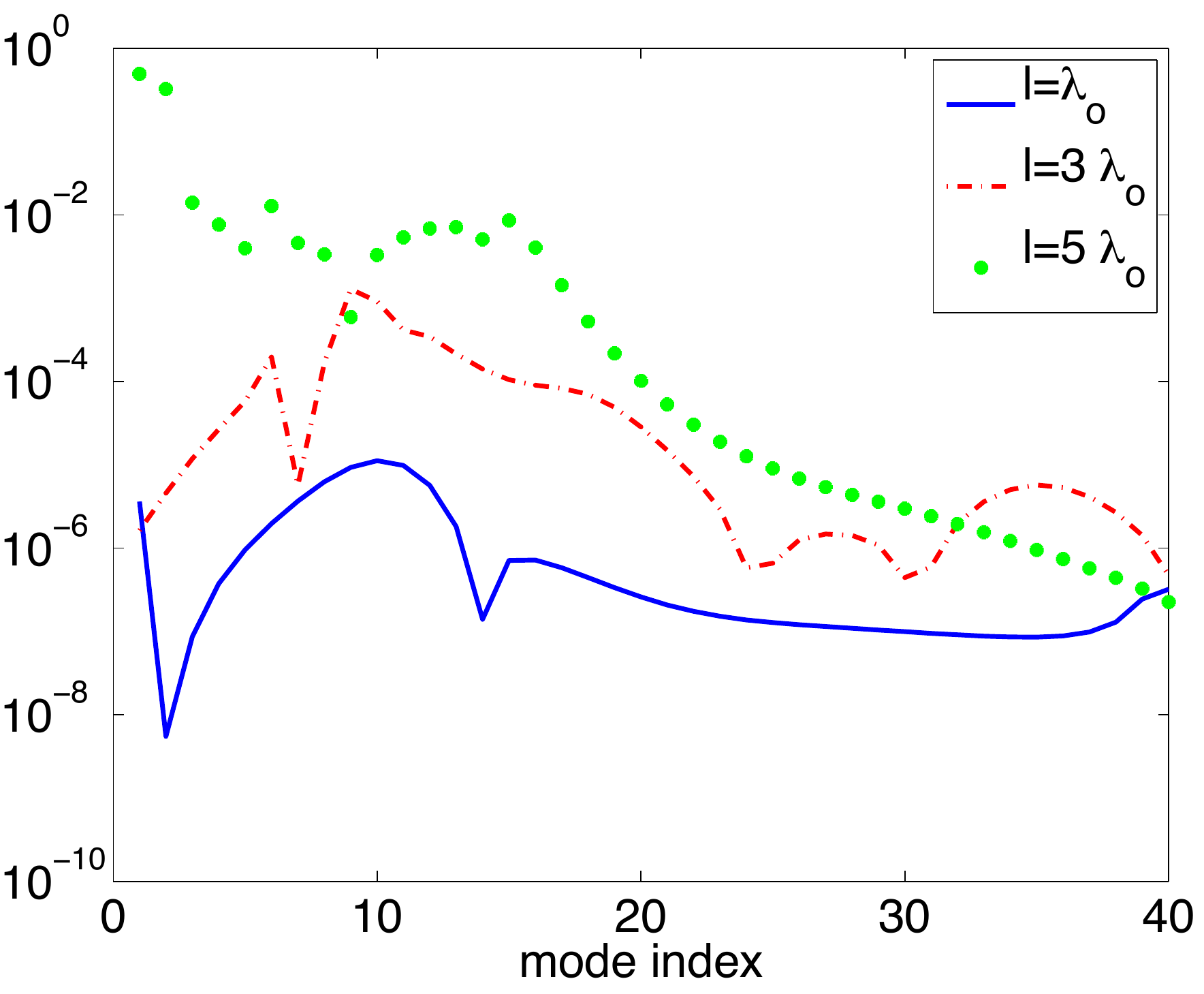}
\end{center}
\vspace{-0.1in}
\caption{Relative error of the predicted eigenvalues by the regular
  perturbation theory equation (\ref{eq:PertL}). Waveguide filled
  with a random medium (left plot) and with a random top boundary
  (right plot). The abscissa is mode index. }
\label{fig:lambda_pert}
\end{figure}

Gathering the  results (\ref{eq:B1})-(\ref{eq:B4}) we see that 
\begin{equation}
\frac{|\Gamma_{jj}|}{H \beta'_j} = \frac{\om_o}{H} \left[ O
  \left(\frac{j^3}{N^2}\right) + O(k_o \ell)\right].
\end{equation}
The second term is due to the the fluctuations in the medium and is
large when $H \sim \om_o$ because $k_o \ell \gg 1$. The first term is
due to the fluctuations of the boundary and is large for $j \gg
N^{2/3}$. In either case, the Frobenius norm of $\Gamma$ is much
larger than that of $H \cB'$, so we have a matrix perturbation
problem. To illustrate the accuracy of the perturbation analysis, we
display in Figure \ref{fig:lambda_pert} the relative error of the
approximation of the eigenvalues for the simulations in section
\ref{sect:I3}.

\section{Proof of Proposition \ref{thm1}}
\label{ap:proof}
That the eigenvectors form an orthonormal basis follows from the
symmetry of $\Upsilon$. We also obtain from
(\ref{eq:ana1})-(\ref{eq:ana1p}) that the quadratic forms of
$\Upsilon$ are
\[
{\bf v}^T \Upsilon {\bf v} = -\sum_{j=2}^N
\Upsilon_{jj-1}\left(v_j-v_{j-1}\right)^2 \le 0, \qquad \forall ~ {\bf
  v} = (v_1, \ldots, v_N)^T\in \mathbb{R}^N,
\]
so the eigenvalues must satisfy $\Lambda_j \le 0.$ We order them as
$0 = \Lambda_1 \ge \Lambda_2 \ge \ldots \Lambda_N$.

We have $(1,1, \ldots, 1)^T \in {\rm Null}(\Upsilon)$ by
construction. To prove property 2, we take a large enough $\gamma$ so
that all the entries in the matrix $ \Upsilon_{_\gamma} = \Upsilon +
\gamma I $ are positive. This matrix is of Perron-Frobenius type, and
its eigenvalues are equal to $\Lambda_j + \gamma$. The largest
eigenvalue $\Lambda_o + \gamma$ is simple, and therefore the null
space of $\Upsilon$ is one dimensional.

The variational definition of $|\Lambda_N|$ as the maximum of the
Rayleigh quotient of $-\Upsilon$ gives that $|\Lambda_N|$ is larger
than $|\Upsilon_{jj}|$, for any $j = 1, \ldots, N$. But $\Upsilon_{NN}
= O(N^2)$ by (\ref{eq:ana2}), and property 3 follows from
$\|\Upsilon\| = |\Lambda_N|.$

Consider square blocks $\Upsilon_m$ of $\Upsilon$, containing the last
$m = O(1)$ elements on its diagonal. By (\ref{eq:ana2}) they scale
like $ \Upsilon_m = N^2 \widetilde{\Upsilon}_m, $ where
$\widetilde{\Upsilon}_m$ have entries of order one. Cauchy's
interlacing theorem gives that $ |\Lambda_{N-m+j}| \ge N^2
|\widetilde{\lambda}_{j}|$, for $ j = 1, \ldots, m,$  where $\widetilde
\lambda_j \le 0 $ are the eigenvalues of $\widetilde{\Upsilon}_m$ in
decreasing order. To prove Property 4 it remains to show that these
are all $O(1)$. First, let us see that $\widetilde{\Upsilon}_m$ has a
trivial null space. Indeed, suppose that ${\bf v} \in {\rm
  Null}(\widetilde{\Upsilon}_m)$ and write equation $
\widetilde{\Upsilon}_m {\bf v} = {\bf 0} $ row by row. Starting from
the last row to the second, and using definitions
(\ref{eq:ana1})-(\ref{eq:ana1p}), we obtain that all entries in ${\bf
  v}$ must be equal to say $v$. However, the first equation gives that
$v = 0$, because the elements in the first row of
$\widetilde{\Upsilon}_m$ do not add to zero. Thus, the null space is
trivial. The smallest in magnitude eigenvalue equals the minimum of
the Rayleigh quotient
\[
\frac{{\bf v}^T (-\widetilde{\Upsilon}_m) {\bf v}}{{\bf v}^T {\bf v}} =
\left[ \widetilde \beta_{N-m} v_1^2 + \displaystyle \sum_{j =
    1}^{m-1} \widetilde \beta_{N-m+j} (v_{j+1}-v_j)^2\right]/{ \displaystyle
  \sum_{j=1}^N v_j^2},
\]
where $\widetilde \beta_j = \beta_j/N^2 = O(1)$ and the right hand
side is obtained by direct calculation. All the terms in this
expression are non-negative and at least one of them must be $O(1)$.
Thus, we see that $|\widetilde \lambda_j| \ge O(1)$ and property 4
follows.

To prove the last property, let $\Lambda$ be a large eigenvalue of
$\Upsilon$ and ${\bf u}$ its associated eigenvector. We see from
definition (\ref{eq:ana1})-(\ref{eq:ana1p}) that $|\Upsilon_{jj}| \ge
\Upsilon_{jj\pm1}$, and using that
\[
|\Lambda| u_j = - \Upsilon_{jj-1} u_{j-1} + |\Upsilon_{jj}| u_j -
\Upsilon_{jj+1} u_{j+1},
\]
we obtain the bound 
\[
|\Upsilon_{jj}| \left(|u_{j-1}| + |u_{j}| + |u_{j+1}|\right) \ge |\Lambda| |u_j|.
\]
Moreover, multiplying by $|u_j|$ and summing over $j = 2, \ldots, J$ we
get the estimate 
\begin{align*}
\sum_{j=2}^J u_j^2 &\le |\Lambda|^{-1} \sum_{j=2}^J |\Upsilon_{jj}|
\left(|u_{j-1} u_j| + u_{j}^2 + |u_{j+1}u_j|\right) \\ &\le
\frac{C\delta N}{N+1-J} \sum_{j=2}^J \left(|u_{j-1} u_j| + u_{j}^2 +
|u_{j+1}u_j|\right),
 \end{align*}
with the second inequality implied by (\ref{eq:ana2}) and $1/|\Lambda|
= \de/N$. Since $J \le N/2$, we have $N/(N+1-J) \le 2$.  Now use
Young's inequality
\[
|u_j u_{j\pm 1}| \le \frac{\widetilde \de u_{j\pm 1}^2}{2} +
\frac{u_j^2}{2 \widetilde \de},
\]
which holds for any $\widetilde \de > 0$. We let $\widetilde \de = C
\de/4$ and obtain that
\begin{align*}
\sum_{j=2}^J u_j^2 &\le C \de \sum_{j=2}^J \left[ \de \left(u_{j-1}^2
  + u_{j+1}^2\right) + 2 \left(1 + \frac{1}{ 4 C \de}
  \right) u_{j}^2 \right]
 \end{align*}
or, equivalently, 
\begin{align*}
\sum_{j=2}^J u_j^2 &\le \frac{4C \de^2}{1- 4 C \de - 4 c \de^2} (u_1^2
+ u_{J+1}^2) \le \frac{4C \de^2}{1- 4 C \de - 4 c \de^2}.
 \end{align*}
The last inequality is because $\|{\bf u}\| = 1$. It remains to show
that $|u_1| \sim \de $. This follows from $ |\Upsilon_{11}| (u_1-u_2)
= |\Lambda| u_1, $ the estimate (\ref{eq:ana2}) that gives
$|\Upsilon_{11}| = O(N)$, and the assumption $|\Lambda| = \de/N$.

\bibliographystyle{plain} \bibliography{BIBLIO}

\end{document}